\documentclass[a4paper,headlines=2.1]{scrartcl}

\usepackage[utf8]{inputenc}
\usepackage{amsmath,amsthm,amssymb}
\usepackage{amsrefs}
\usepackage{tikz}
\usetikzlibrary{arrows,shapes,automata,petri,positioning,calc,arrows.meta,fit}
\usepackage{hyperref}
\usepackage{bbm}
\usepackage{xcolor}
\usepackage{tabto}
\usepackage{bm}
\usepackage{enumitem}
\usepackage{authblk}
\usepackage{mathtools}

\usepackage{xr}

\newtheorem{theorem}{Theorem}

\allowdisplaybreaks

\newcommand{\R}{\mathbb{R}}

\newcommand{\N}{\mathbb{N}}
\newcommand{\E}{\mathbb{E}}

\newcommand{\cW}{\mathcal{W}}
\newcommand{\Cov}{\operatorname{Cov}}

\newcommand{\ol}{\overline}
\newcommand{\tr}{\operatorname{tr}}

\newcommand{\tX}{\tilde{X}}

\newsavebox\thesmashminipage

\newcommand\s{1}
\renewcommand\d{\s cm}
\tikzset{
	black/.style={
		circle,
		thick,
		draw=black,
		fill=gray!50,
		minimum size=0.05*\d,
		scale = \s,
	},
	graynode/.style={
		circle,
		thick,
		draw=black,
		fill=gray!20,
		minimum size=0.05*\d,
		scale = \s,
	},
	white/.style={
		circle,
		thick,
		draw=black,
		minimum size=0.05*\d,
		scale = \s,
	},
	noborder/.style={
		circle,
		thick,
		minimum size=0.05*\d,
		scale = \s,
	},
	el/.style = {inner sep=2pt, align=left, sloped, font=\tiny},
	every label/.append style = {font=\tiny},
	position/.style args={#1:#2 from #3}{
		at=(#3.#1), anchor=#1+180, shift=(#1:#2)
	}
}

\numberwithin{equation}{section}

\AtBeginEnvironment{biblist}{\catcode`\#=12 }

\title{Trace Moments of the Sample Covariance Matrix with Graph-Coloring}
\date{}
\author{Ben Deitmar}
\affil{\small \textit{Department of Mathematical Stochastics, ALU Freiburg \protect\\ Ernst-Zermelo-Str. 1, 79104 Freiburg, Germany \protect\\ E-mail: ben.deitmar@stochastik.uni-freiburg.de}}

\begin{document}
	
\thispagestyle{empty}
\maketitle
\vspace{-1.2cm}
	
\begin{abstract}
	\begin{center}\textbf{Abstract}\end{center}
	\noindent
	Let $\bm{S}_{p,n}$ denote the sample covariance matrix based on $n$ independent identically distributed $p$-dimensional random vectors in the null-case. The main result of this paper is an explicit expansion of trace moments and power-trace covariances of $\bm{S}_{p,n}$ simultaneously for both high- and low-dimensional data. To this end we expand a well-known ansatz of describing trace moments as weighted sums over routes or graphs. The novelty to our approach is an inherent coloring of the examined graphs and a decomposition of graphs into their tree-structure and their \textit{seed graphs}, which allows for some elegant formulas explaining the effect of the tree structures on the number of Euler-tours. The weighted sums over graphs become weighted sums over the possible seed graphs, which in turn are much easier to analyze.
\end{abstract}


\vspace{0.0cm}

\section{Introduction}\label{SectionIntroduction}
Let $\bm{S}_{p,n} \in \R^{p \times p}$ denote the sample covariance matrix for a data set of $n$ independent $p$-dimensional random vectors $(X_{1\,j},...,X_{p\,j})^T, \, 1 \leq j \leq n$, where each random vector consists of iid components with mean zero and variance one. If $\bm{X}_{p,n}=(X_{i \, j})$ is the corresponding $(p \times n)$ data-matrix, then the sample covariance matrix is given by $\bm{S}_{p,n} = \frac{1}{n} \bm{X}_{p,n} \bm{X}_{p,n}^T$.

\hspace{-0.7cm}\rule[-0.5cm]{8cm}{0.2mm}\\
\\
\hspace{1cm}
Supported by the DFG Research Unit 5381

\newpage
\noindent
Throughout we will assume $p \leq n$. We can do this almost without loss of generality, since the cyclic property of the trace implies $\tr(\bm{S}_{p,n}^l) = \frac{p^l}{n^l} \tr(\bm{S}_{n,p}^l)$.\\
For any two numbers $a,b \in \R$ use the notation $a \land b = \min(a,b)$ and $a \lor b = \max(a,b)$.

\subsection{Main results}
The main contribution of this paper is the development of combinatorial methods for analyzing the effect of tree-structures of graphs on the number of Euler-tours. This allows for trace moment expansions of the form of the following two theorems.
\begin{theorem}\label{Thm_MeanExpansion}
	For any $l \in \N$, assume that $\E[X_{1 \, 1}^{2l}] < \infty$. Then for $p \leq n$ we have
	\begin{align*}
		\E[\tr(\bm{S}_{p,n}^l)]&\\
		= \sum\limits_{b=1}^{l\land p} \Bigg[ & \frac{{p \choose b}{n-b \choose l+1-b}}{n^l} \, l! {l-1 \choose b-1}\\
		& \hspace{0cm} + \frac{{p \choose b}{n-b \choose l-b}}{n^l} \Big\{ A_1(l,b) + (\E[X_{1 \, 1}^4]-3) \, A_2(l,b) \Big\} +  \mathcal{O}\left(\frac{p^b}{n^{b+1}}\right) \Bigg] \ ,
	\end{align*}
	where
	\begin{align*}
		& A_1(l,b) := \frac{b! (l-b)!}{2} \left( {2l \choose 2b} + (2b-1) {l \choose b}^2 \right)
	\end{align*}
	and the 'fourth moment correction term' is $A_2(l,b) := b! \, (l-b)! {l \choose b-1}{l \choose b+1}$.
\end{theorem}\
Proof in Section \ref{MainResult1}
\begin{theorem}\label{Thm_CovExpansion}
	For any $l_1,l_2 \in \N$, assume that $\E[X_{1 \, 1}^{2l_1+2l_2}] < \infty$. Then for $p \leq n$ we have
	\begin{align*}
		& \hspace{-0.5cm} \operatorname{Cov}[\tr(\bm{S}_{p,n}^{l_1}), \tr(\bm{S}_{p,n}^{l_2})]\\
		& = \sum\limits_{b=1}^{l_1+l_2 \land p} \Bigg[ \frac{{p \choose b} {n-b \choose l_1+l_2-b}}{n^{l_1+l_2}} \Big\{ C_1(l_1,l_2,b) + \big(\E[X^{4}_{1 \, 1}]-3\big) C_2(l_1,l_2,b) \Big\} + \mathcal{O}\left(\frac{p^{b}}{n^{b+1}}\right) \Bigg] \ ,
	\end{align*}
	where
	\begin{align*}
		& C_1(l_1,l_2,b) := 2 \, b! \, (l_1+l_2-b)! \sum\limits_{k=0}^{b} {l_1 \choose k} {l_2 \choose b - k} \sum\limits_{m=0}^{b-k} m {l_1 \choose k+m} {l_2 \choose b-m-k}\\
		& C_2(l_1,l_2,b) := b! \, (l_1+l_2-b)! \sum\limits_{k=0}^{b - 1} {l_1 \choose k} {l_1 \choose k+1} {l_2 \choose b-1-k} {l_2 \choose b - k} \ .
	\end{align*}
\end{theorem}\
Proof in Section \ref{MainResult2}\\
\\
We would like to emphasize the fact that these expansions - using $\tr(\bm{S}_{p,n}^l) = \frac{p^l}{n^l} \tr(\bm{S}_{n,p}^l)$ - hold in any asymptotic regime for the behavior of $\frac{p}{n}$, whereas most results on trace moments of sample covariance matrices are either in the regime $\frac{p}{n} \rightarrow y > 0$ (\cite{BaiCLT}, \cite{Jonsson}, \cite{NajimYao}, \cite{Peche}) or in the regime $\frac{p}{n} \searrow 0$ (\cite{AndersonCLT} by Anderson and Zeitouni). Also the methods are highly adaptable, as we showcase in Section \ref{Section_Adaptability} of the appendix by generalizing our main results to the setting of complex entries and to a re-sampling setting.

\subsection{Related works}
Bai and Silverstein in \cite{BaiCLT} were able to show a CLT for $\tr(f(\bm{S}_{p,n}))$, where $f$ is an analytic functional. In their paper they also derived the limits of the mean $\E[\tr(\bm{S}_{p,n}^l)]$ and covariance $\operatorname{Cov}[\tr(\bm{S}_{p,n}^{l_1}),\tr(\bm{S}_{p,n}^{l_2})]$ for the asymptotic regime $\frac{p}{n} \rightarrow y > 0$ (see Equations (1.23) and (1.24) of \cite{BaiCLT}) based on previous work by Jonsson \cite{Jonsson}. They however needed to assume $\E[X^4_{i \, j}]=3$ and their results are strictly asymptotic in nature.\\
\\
A more modern result in the same direction - but for arbitrary complex entries with existing fourth moment - was shown by Najim and Yao in \cite{NajimYao}. They showed first a uniform CLT for Stieltjes transforms
\begin{align}
	& s(z) = \tr\big( (\bm{S}_{p,n} - \operatorname{Id}_{p \times p})^{-1} \big) = \frac{1}{p} \sum\limits_{j=1}^p \frac{1}{\lambda_j(\bm{S}_{p,n})-z}
\end{align}
and with Helffer-Sjöstrand calculus derived a CLT for linear spectral statistics for non-analytic functionals. In Section 3.5 of \cite{NajimYao} they give the asymptotic mean- and covariance-structure of the limiting Gaussian distribution for their Stieltjes transform CLT in terms of solutions to self-consistent equations. With their formulas it might be possible to recover the explicit limits of $\E[\tr(\bm{S}_{p,n}^l)]$ and $\operatorname{Cov}[\tr(\bm{S}_{p,n}^{l_1}), \tr(\bm{S}_{p,n}^{l_2})]$ in the regime $\frac{p}{n} \rightarrow y > 0$.\\
Our Theorem \ref{ThmComplexEntries} provides these limits directly - making the application of the strong results of \cite{NajimYao} by Najim and Yao more straight forward - and our methods work equally well in settings where the Stieltjes transform method is not yet able to formulate CLTs, such as the setting of Theorem \ref{Thm_ReSample}. Additionally the ideas developed here may be used to find even more terms in the expansions of the mean and covariance, allowing for even higher accuracy of such CLTs.\\
\\
In the regime $\frac{p}{n} \searrow 0$ Anderson and Zeitouni in \cite{AndersonCLT} showed a general CLT for regularized sample covariance matrices and were able to explicitly give the mean and covariance of the limiting Gaussian process. Unfortunately their joint cumulant summability condition (see Assumption 2.2 of \cite{AndersonCLT}) is not satisfied in the null case, which we are considering, since for example when $r=2$ we have
\begin{align*}
	& \sum\limits_{i_1,...,i_r=0}^\infty \big| \bm{C}(X_{i_1\,1},...,X_{i_r\,1}) \big| = \sum\limits_{i=0}^\infty \big| \bm{C}(X_{i\,1},...,X_{i\,1}) \big| \overset{r=2}{=} \sum\limits_{i=0}^\infty \operatorname{Var}[X_{i\,1}] = \infty \ .
\end{align*}\
\\
\\
Another asymptotic regime of interest was studied by P\'ech\'e in \cite{Peche} and earlier in a series of works by Sinai and Soshnikov (\cite{Soshnikov1}, \cite{Soshnikov2}, \cite{Soshnikov3}, \cite{Soshnikov4}). Assume $\frac{p}{n} \rightarrow y >0$ and let $l=l_n$ grow with $n$. In the case $l_n \sim n^{\frac{2}{3}}$ the limiting behavior of the trace moments $\E\big[ \tr(\bm{S}_{p,n}^{l_n}) \big]$ can be shown to determine the limiting behavior of the eigenvalues of $\bm{S}_{p,n}$ at the correct scaling to achieve Tracy-Widom results. P\'ech\'e uses this to show universality of the Tracy-Widom law. In the paper \cite{Peche} P\'ech\'e also develops a CLT for the regime $l_n << \sqrt{n}$, where she only focuses on the first-order expansion
\begin{align*}
	& \E\big[ \tr(\bm{S}_{p,n}^{l_n}) \big] = \sum\limits_{k=1}^{l_n} \Big( \frac{n}{p} \Big)^k \frac{1}{l_n} {l_n \choose k} {l_n-1 \choose k-1} \big( 1 + o(1) \big) \ .
\end{align*}
This expansion is consistent with our Theorem \ref{Thm_MeanExpansion}, but our bounds to get the $\mathcal{O}\left(\frac{p^b}{n^{b+1}}\right)$ term assume constant exponent $l$. Luckily, the exact same arguments as in the first two paragraphs of the proof of Proposition 2.4 in \cite{Peche} may also be used to show that the result of our Theorem 1 - or in higher generality the equality (\ref{Eq_ComplexResult1}) - still holds for sub-Gaussian entries in the regime $l_n << n^{\frac{1}{4}}$, where the $\mathcal{O}\left(\frac{p^b}{n^{b+1}}\right)$ term must be replaced with $\mathcal{O}\left( \frac{l_n^4}{n} \right)$. 
\\
\\
For sample covariance matrices $\bm{S}_{p,n}$ the exact formulas for the trace moments $\E\big[ \tr(\bm{S}_{p,n}^l) \big]$ are only known in the setting where the entries $(X_{i \, j})_{i,j \in \N}$ are iid complex standard normal, in other words when $\bm{S}_{p,n}$ is an isotropic complex Wishart matrix. The formulas can be found with some generalizations of Harer-Zagier recursion and are
\begin{align*}
	& \E\big[ \tr(\bm{S}_{p,n}^l) \big] = \frac{l!}{n^l} \sum\limits_{b,w \in \N} {p \choose b} {n \choose w} {l-1 \choose b-1, w-1} \ .
\end{align*}
This is for example shown in Corollary 1.9 of \cite{Vassilieva} where Vassilieva derives a new representation of $\E\big[ \tr\big((A \bm{X}_{p,n} B \bm{X}_{p,n}^*)^l \big) \big]$ for fixed matrices $A,B$ and the above exact formula arises as a corollary of her main theorem.\\
Exact formulas for the covariances $\Cov\big[ \tr(\bm{S}_{p,n}^{l_1}), \tr(\bm{S}_{p,n}^{l_2}) \big]$ are not known in any setting.

\subsection{Overview of our method}\label{MethodOverview}
We can write the mean of the trace $\tr(\bm{S}_{p,n}^{l})$ as
\begin{align}\label{Eq_TraceRepresentation1}
	& \E[\tr(\bm{S}^l)] = \frac{1}{n^{l}} \sum\limits_{\substack{\bm{i} \in [p]^l}} \sum\limits_{\bm{j} \in [n]^l} \underbrace{\E[(X_{i_1 \, j_1} X_{i_2 \, j_1}) \, (X_{i_2 \, j_2} X_{i_3 \, j_2}) \, \dots \, (X_{i_l \, j_l} X_{i_1\, j_l})]}_{=: W_{\bm{i},\bm{j}}} \ .
\end{align}
If for each entry $X_{a \, b}$ in the product $(X_{i_1 \, j_1} X_{i_2 \, j_1}) \, (X_{i_2 \, j_2} X_{i_3 \, j_2}) \, \dots \, (X_{i_l \, j_l} X_{i_1\, j_l})$ one draws an edge from $a$ to $b$, one arrives at a graph of the form\\
\\
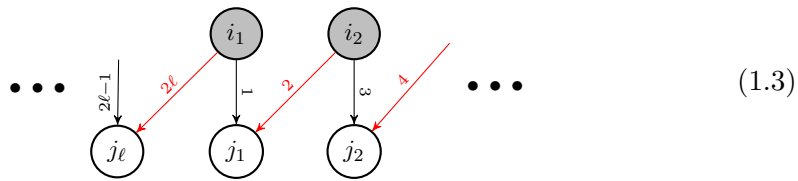
\begin{minipage}{0.95\textwidth}
	\centering
	\renewcommand\s{0.85}
	\begin{align}\label{Pic_ReversedCircuitGraph}
		&
	\end{align}
	\vspace{-2cm}\\
	\begin{tikzpicture}[node distance=\d and \d,>=stealth',auto, every place/.style={draw}]
		\node [black] (i1) {$i_1$};
		\node [black] (i2) [right=of i1] {$i_2$};
		\node [white] (j1) [below=of i1] {$j_1$};
		\node [white] (j2) [below=of i2] {$j_2$};
		\node [noborder] (i3) [right=of i2] {};
		\node [noborder] (i0) [left=of i1] {\phantom{$i_{\ell}$}};
		\node [white] (j0) [left=of j1] {$j_{\ell}$};

		\node [noborder] (t0) [above left =0.25*\d and 0.5*\d of j0] {$\bullet \, \bullet \, \bullet$};
		\node [noborder] (t1) [below right =0.25*\d and 0*\d of i3] {$\bullet \, \bullet \, \bullet$};
		
		\path[->]
		(i0) edge node[el,above] {$2\ell\hspace{-0.08cm}-\hspace{-0.08cm}1$} (j0)
		(i1) edge [color=red] node[el,above] {$2\ell$} (j0)
		(i1) edge node[el,above] {$1$} (j1)
		(i2) edge [color=red] node[el,above] {$2$} (j1)
		(i2) edge node[el,above] {$3$} (j2)
		(i3) edge [color=red] node[el,above] {$4$} (j2);
		
	\end{tikzpicture}
\end{minipage}\\
\\
By changing the direction of every second edge, the graph\\
\\
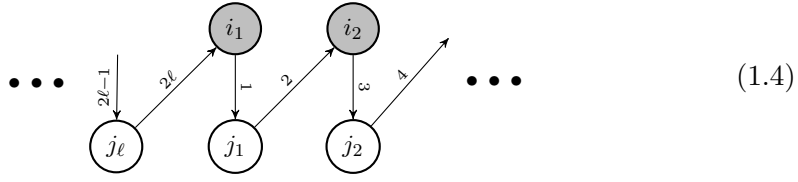
\begin{minipage}{0.95\textwidth}
	\centering
	\renewcommand\s{0.85}
	\begin{align}\label{Pic_CircuitGraph1}
		&
	\end{align}
	\vspace{-2cm}\\
	\begin{tikzpicture}[node distance=\d and \d,>=stealth',auto, every place/.style={draw}]
		\node [black] (i1) {$i_1$};
		\node [black] (i2) [right=of i1] {$i_2$};
		\node [white] (j1) [below=of i1] {$j_1$};
		\node [white] (j2) [below=of i2] {$j_2$};
		\node [noborder] (i3) [right=of i2] {};
		\node [noborder] (i0) [left=of i1] {\phantom{$i_{\ell}$}};
		\node [white] (j0) [left=of j1] {$j_{\ell}$};

		\node [noborder] (t0) [above left =0.25*\d and 0.5*\d of j0] {$\bullet \, \bullet \, \bullet$};
		\node [noborder] (t1) [below right =0.25*\d and 0*\d of i3] {$\bullet \, \bullet \, \bullet$};
		
		\path[->]
		(i0) edge node[el,above] {$2\ell\hspace{-0.08cm}-\hspace{-0.08cm}1$} (j0)
		(j0) edge node[el,above] {$2\ell$} (i1)
		(i1) edge node[el,above] {$1$} (j1)
		(j1) edge node[el,above] {$2$} (i2)
		(i2) edge node[el,above] {$3$} (j2)
		(j2) edge node[el,above] {$4$} (i3);
		
	\end{tikzpicture}
\end{minipage}\\
\\
describes a walk through the vertices $\{i_1,...,i_l\} \cup \{j_1,...,j_l\}$. We call a directed multigraph of the above form a \textit{circuit (multi-)graph}. A circuit graph has an inherent ordering to its edges and is uniquely defined by $\bm{i}$ and $\bm{j}$. The exact order in which the vertices are traversed is given by the zipped sequence $\langle \bm{i},\bm{j} \rangle := (i_1,j_1,i_2,j_2,...,i_l,j_l)$.\\
\\
In the above pictures we have given each $i_\bullet$ and $j_\bullet$ a separate vertex. If $i_{s} = j_t$, then in the above pictures the vertices representing $i_{s}$ and $j_t$ will be merged. As an example let $\bm{i} = (v_2,v_3,v_1)$ and $\bm{j} = (v_3,v_4,v_4)$, then the resulting circuit graph $G_{\langle \bm{i},\bm{j} \rangle} = G_{(v_2,v_3,v_3,v_4,v_1,v_4)}$ is:\\
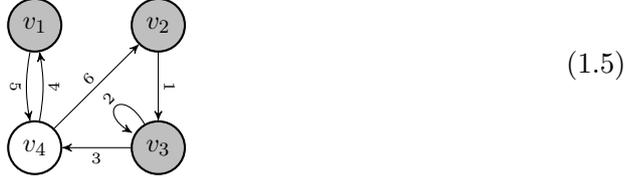
\begin{minipage}{0.95\textwidth}
	\centering
	\renewcommand\s{0.9}
	\begin{align}\label{Pic_ExampleCircuit}
		& 
	\end{align}
	\vspace{-1.9cm}\\
	\begin{tikzpicture}[node distance=\d and \d,>=stealth',auto, every place/.style={draw}]
		\node [black] (v1) {$v_1$};
		\node [black] (v2) [right=of v1] {$v_2$};
		\node [white] (v3) [below=of v1] {$v_4$};
		\node [black] (v4) [below=of v2] {$v_3$};
		
		\path[->]
		(v2) edge [bend right=0] node[el,above] {$1$} (v4)
		(v4) edge [in=150,out=120,loop] node[el,above] {$2$} (v4)
		(v4) edge [bend right=0] node[el,below] {$3$} (v3)
		(v3) edge [bend right=10] node[el,below] {$4$} (v1)
		(v1) edge [bend right=10] node[el,below] {$5$} (v3)
		(v3) edge [bend right=0] node[el,above] {$6$} (v2);
		
	\end{tikzpicture}
\end{minipage}\\
\\
A vertex in $G_{\langle \bm{i},\bm{j} \rangle}$ is colored black, if it is in the set $\bm{\{i\}} = \{i_1,...,i_l\}$. This corresponds to the heuristic, that the color black is dominant when merging two vertices from the picture (\ref{Pic_CircuitGraph1}). The reversed graph $\operatorname{R}(G_{\langle \bm{i},\bm{j} \rangle})$ is defined as a copy of $G_{\langle \bm{i},\bm{j} \rangle}$ where the direction of each even numbered edge is reversed, which now corresponds to the picture (\ref{Pic_ReversedCircuitGraph}) with merged vertices. In tune with our example (\ref{Pic_ExampleCircuit}) we now have $\operatorname{R}\big(G_{(v_2,v_3,v_3,v_4,v_1,v_4)}\big)$:\\
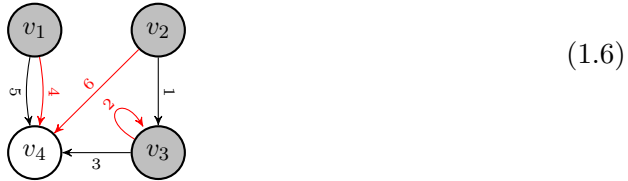
\begin{minipage}{0.95\textwidth}
	\centering
	\renewcommand\s{0.9}
	\begin{align}\label{Pic_ExampleCircuitReversed}
		& 
	\end{align}
	\vspace{-1.7cm}\\
	\begin{tikzpicture}[node distance=\d and \d,>=stealth',auto, every place/.style={draw}]
		\node [black] (v1) {$v_1$};
		\node [black] (v2) [right=of v1] {$v_2$};
		\node [white] (v3) [below=of v1] {$v_4$};
		\node [black] (v4) [below=of v2] {$v_3$};
		
		\path[->]
		(v2) edge [bend right=0] node[el,above] {$1$} (v4)
		(v4) edge [in=120,out=150,loop, color=red] node[el,above] {$2$} (v4)
		(v4) edge [bend right=0] node[el,below] {$3$} (v3)
		(v1) edge [bend right=-10, color=red] node[el,above] {$4$} (v3)
		(v1) edge [bend right=10] node[el,below] {$5$} (v3)
		(v2) edge [bend right=0, color=red] node[el,above] {$6$} (v3);
		
	\end{tikzpicture}
\end{minipage}\\
\\
Since the entries $(X_{i \, j})_{i,j \in \N}$ are assumed to be independent, the mean
\begin{align}\label{Eq_Weight1}
	& W_{\bm{i},\bm{j}} = \E[(X_{i_1 \, j_1} X_{i_2 \, j_1}) \, (X_{i_2 \, j_2} X_{i_3 \, j_2}) \, \dots \, (X_{i_l \, j_l} X_{i_1\, j_l})]
\end{align}
has the product form $\prod\limits_{(s,t) \in [p] \times [n]} \E\big[X_{s,t}^{A_{s,t}}\big]$, where $A_{s,t}$ is the number of occurrences of $X_{s \, t}$ in (\ref{Eq_Weight1}). By construction of $\operatorname{R}(G_{\langle \bm{i},\bm{j} \rangle})$, each entry $X_{s,t}$ in (\ref{Eq_Weight1}) is represented by one edge in the directed multigraph $\operatorname{R}(G_{\langle \bm{i},\bm{j} \rangle})$ and $(A_{s,t})_{s \in [p], t \in [n]}$ must be its adjacency matrix. We thus have
\begin{align}\label{Eq_Weight2}
	& W_{\bm{i},\bm{j}} = \prod\limits_{(s,t) \in [p] \times [n]} \E\big[X_{s,t}^{A_{s,t}(\operatorname{R}(G_{\langle \bm{i},\bm{j} \rangle}))}\big] \ .
\end{align}
\\
Define $W(G_{\langle \bm{i},\bm{j} \rangle}) := W_{\bm{i},\bm{j}}$, then the formula (\ref{Eq_TraceRepresentation1}) becomes a weighted sum over graphs, where our first step will be to split the sum by number of total and black vertices:
\begin{align}\label{Eq_TraceRepresentation2}
	& \E[\tr(\bm{S}^l)] = \frac{1}{n^{l}} \sum\limits_{\substack{\bm{i} \in [p]^l}} \sum\limits_{\bm{j} \in [n]^l} W(G_{\langle \bm{i},\bm{j} \rangle}) = \frac{1}{n^{l}} \sum\limits_{r=1}^n \sum\limits_{\substack{\bm{i} \in [p]^l, \, \bm{j} \in [n]^l \\ \#\bm{\{i\}} \cup \bm{\{j\}} = r}} W(G_{\langle \bm{i},\bm{j} \rangle}) \nonumber\\
	& = \frac{1}{n^{l}} \sum\limits_{r=1}^n \sum\limits_{b=1}^{r \land p} \sum\limits_{\substack{\bm{i} \in [p]^l, \, \bm{j} \in [n]^l \\ \#\bm{\{i\}} \cup \bm{\{j\}} = r \\ \#\bm{\{i\}} = b}} W(G_{\langle \bm{i},\bm{j} \rangle})
\end{align}
As the entries $Y_{i \, j}$ are assumed to be iid, we can exchange the set $\bm{\{i\}} \subset [p]$ with $\{1,...,b\} = [b]$ and the set $\bm{\{j\}}\setminus\bm{\{j\}} \subset [n]$ with $\{b+1,...,r\}$ without changing the weight $W(G_{\langle \bm{i},\bm{j} \rangle})$. This yields
\begin{align}\label{Eq_TraceRepresentation3}
	& \E[\tr(\bm{S}^l)] = \frac{1}{n^{l}} \sum\limits_{r=1}^n \sum\limits_{b=1}^{r \land p} {p \choose b} {n-b \choose r-b} \sum\limits_{\substack{\bm{i} \in [b]^l, \, \bm{j} \in [r]^l \\ \bm{\{i\}} \cup \bm{\{j\}} = [r] \\ \bm{\{i\}} = [b]}} W(G_{\langle \bm{i},\bm{j} \rangle}) \ .
\end{align}
The product form (\ref{Eq_Weight2}) tells us that the weight $W(G_{\langle \bm{i},\bm{j} \rangle})$ must already be zero, if there is ever only a single edge (regardless of direction) between two vertices. For the graph to have non-zero weight, each connection between vertices must have at least two edges, meaning we can have at most $l = \frac{\#\text{edges}}{2}$ many connections. As the graph $G_{\langle \bm{i},\bm{j} \rangle}$ is by construction connected (excluding vertices from $[n]$ which do not occur in $\bm{\{i\}} \cup \bm{\{j\}}$), there can be at most $l+1$ many vertices in $\bm{\{i\}} \cup \bm{\{j\}}$. We can adjust the above formula to
\begin{align}\label{Eq_TraceRepresentation4}
	& \E[\tr(\bm{S}^l)] = \sum\limits_{r=1}^{(l+1) \land n} \sum\limits_{b=1}^{r \land p} \underbrace{\frac{{p \choose b} {n-b \choose r-b}}{n^l}}_{\mathcal{O}\big( \frac{p^b n^{r-b}}{n^l} \big)} \overbrace{\sum\limits_{\substack{\bm{i} \in [b]^l, \, \bm{j} \in [r]^l \\ \bm{\{i\}} \cup \bm{\{j\}} = [r] \\ \bm{\{i\}} = [b]}} W(G_{\langle \bm{i},\bm{j} \rangle})}^{\text{independent of }p,n} \nonumber\\
	& = \sum\limits_{r=l}^{l+1} \sum\limits_{b=1}^{r \land p} \Bigg[ \frac{{p \choose b} {n-b \choose r-b}}{n^l} \sum\limits_{\substack{\bm{i} \in [b]^l, \, \bm{j} \in [r]^l \\ \bm{\{i\}} \cup \bm{\{j\}} = [r] \\ \bm{\{i\}} = [b]}} W(G_{\langle \bm{i},\bm{j} \rangle}) + \mathcal{O}\Big( \frac{p^b n^{l-1-b}}{n^l} \Big) \Bigg] \nonumber\\
	& = \sum\limits_{b=1}^{(l+1) \land p} \Bigg[ \sum\limits_{r=l \lor b}^{l+1} \frac{{p \choose b} {n-b \choose r-b}}{n^l} \sum\limits_{\substack{\bm{i} \in [b]^l, \, \bm{j} \in [r]^l \\ \bm{\{i\}} \cup \bm{\{j\}} = [r] \\ \bm{\{i\}} = [b]}} W(G_{\langle \bm{i},\bm{j} \rangle}) + \mathcal{O}\Big( \frac{p^b}{n^{b+1}} \Big) \Bigg] \ .
\end{align}
We are now only interested in the cases $r=l+1$ and $r=l$, where $r$ describes the total number of vertices in the graph. We will see in Lemma \ref{WeightOfTrees}, that $G_{\bm{i},\bm{j}}$ with $l+1$ many vertices and non-zero weight must always have a tree structure, which makes them easy to count, and must always have weight $1$. In Lemma \ref{CountingColoredTrees} we show that there are $l! {l-1 \choose b-1}$ many such graphs with $b$ many black vertices, which yields
\begin{align}
	& \sum\limits_{\substack{\bm{i} \in [b]^l, \, \bm{j} \in [l+1]^l \\ \bm{\{i\}} \cup \bm{\{j\}} = [l+1] \\ \bm{\{i\}} = [b]}} W(G_{\langle \bm{i},\bm{j} \rangle}) = l! {l-1 \choose b-1} \ .
\end{align}
In Proposition \ref{WeightOfFusedTrees} we will see that all $G_{\bm{i},\bm{j}}$ with $l$ many vertices and non-zero weight fall into one of the three categories, for which we have depicted an example each here:\\
\\
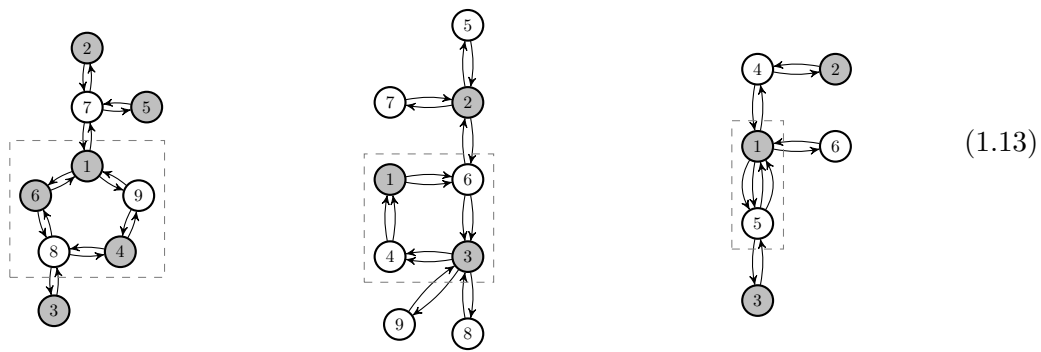
\begin{minipage}{0.3\textwidth} 
	\centering
	\renewcommand\s{0.6}
	\begin{tikzpicture}[node distance=\d and \d,>=stealth',auto, every place/.style={draw}]
		\node [black] (v1) {$1$};
		\node [white] (v9) [position=-30:{\d} from v1] {$9$};
		\node [black] (v4) [position=-108:{\d} from v9] {$4$};
		\node [black] (v6) [position=-150:{\d} from v1] {$6$};
		\node [white] (v8) [position=-72:{\d} from v6] {$8$};
		
		\node [white] (v7) [position=90:{\d} from v1] {$7$};
		
		\node [black] (v2) [position=90:{\d} from v7] {$2$};
		
		\node [black] (v5) [position=0:{\d} from v7] {$5$};
		
		\node [black] (v3) [position=-90:{\d} from v8] {$3$};
		
		\path[->]
		(v1) edge [bend right=10] (v9)
		(v9) edge [bend right=10] (v1)
		
		(v9) edge [bend right=10] (v4)
		(v4) edge [bend right=10] (v9)
		
		(v4) edge [bend right=10] (v8)
		(v8) edge [bend right=10] (v4)
		
		(v8) edge [bend right=10] (v6)
		(v6) edge [bend right=10] (v8)
		
		(v6) edge [bend right=10] (v1)
		(v1) edge [bend right=10] (v6)
		
		(v3) edge [bend right=10] (v8)
		(v8) edge [bend right=10] (v3)
		
		(v1) edge [bend right=10] (v7)
		(v7) edge [bend right=10] (v1)
		
		(v7) edge [bend right=10] (v5)
		(v5) edge [bend right=10] (v7)
		
		(v7) edge [bend right=10] (v2)
		(v2) edge [bend right=10] (v7)
		;
		
		\node[fit=(v1)(v6)(v9)(v4), draw, dashed,gray] {};
		
	\end{tikzpicture}
\end{minipage}
\begin{minipage}{0.3\textwidth}
	\centering
	\renewcommand\s{0.6}
	\begin{tikzpicture}[node distance=\d and \d,>=stealth',auto, every place/.style={draw}]
		\node [black] (v1) {$1$};
		\node [white] (v6) [right=of v1] {$6$};
		\node [white] (v4) [below=of v1] {$4$};
		\node [black] (v3) [below=of v6] {$3$};
		
		\node [black] (v2) [above=of v6] {$2$};
		\node [white] (v5) [above=of v2] {$5$};
		\node [white] (v7) [left=of v2] {$7$};
		\node [white] (v8) [below=of v3] {$8$};
		\node [white] (v9) [below left=of v3] {$9$};
		
		\path[->]
		(v1) edge [bend right=10] (v6)
		(v6) edge [bend right=10] (v3)
		(v3) edge [bend right=10] (v4)
		(v4) edge [bend right=10] (v1)
		(v1) edge [bend right=-10] (v6)
		(v6) edge [bend right=-10] (v3)
		(v3) edge [bend right=-10] (v4)
		(v4) edge [bend right=-10] (v1)
		
		(v6) edge [bend right=-10] (v2)
		(v2) edge [bend right=-10] (v6)
		
		(v7) edge [bend right=-10] (v2)
		(v2) edge [bend right=-10] (v7)
		
		(v5) edge [bend right=-10] (v2)
		(v2) edge [bend right=-10] (v5)
		
		(v8) edge [bend right=-10] (v3)
		(v3) edge [bend right=-10] (v8)
		
		(v9) edge [bend right=-10] (v3)
		(v3) edge [bend right=-10] (v9)
		;
		
		\node[fit=(v1)(v3), draw, dashed,gray] {};
	\end{tikzpicture}
\end{minipage}
\begin{minipage}{0.39\textwidth}
	\centering
	\renewcommand\s{0.6}
	\begin{align}\label{Pic_SeedExamples}
		& \phantom{h} \nonumber\\
		& \phantom{h}
	\end{align}
	\vspace{-2.2cm}\\
	\hspace{-0.9cm}
	\begin{tikzpicture}[node distance=\d and \d,>=stealth',auto, every place/.style={draw}]
		\node [black] (v1) {$1$};
		\node [white] (v5) [below=of v1] {$5$};
		\node [white] (v6) [right=of v1] {$6$};
		\node [white] (v4) [above=of v1] {$4$};
		\node [black] (v2) [right=of v4] {$2$};
		\node [black] (v3) [below=of v5] {$3$};
		
		\path[->]
		(v1) edge [bend right=10] (v4)
		(v4) edge [bend right=10] (v1)
		
		(v4) edge [bend right=10] (v2)
		(v2) edge [bend right=10] (v4)
		
		(v1) edge [bend right=10] (v5)
		(v5) edge [bend right=10] (v1)
		
		(v5) edge [bend right=10] (v3)
		(v3) edge [bend right=10] (v5)
		
		(v1) edge [bend right=10] (v6)
		(v6) edge [bend right=10] (v1)
		
		(v1) edge [bend right=30] (v5)
		(v5) edge [bend right=30] (v1)
		;
		
		\node[fit=(v1)(v5), draw, dashed,gray] {};
	\end{tikzpicture}
\end{minipage}\\
\\
Each of these examples has a tree structure surrounding a central part of the graph, which we have marked with a dashed box. The central part will be called the \textit{seed graph} and can always be defined by iteratively removing leaves of the graph until there are no leaves left. (We use a way of removing leaves, which preserves color in the remaining vertices and ensures the graph remains a circuit graph. See Definition \ref{DefRemovingLeaves} for details.) The examples are then called \textit{sprouts} of their respective seed graphs. The three relevant categories are defined by the structure of their seed graph.
\begin{itemize}
	\item Category 1:\\
	The seed graph of the left most graph in (\ref{Pic_SeedExamples}) has a ring-like structure of length $l_0 = 5$ and is traversed in both directions. We thus call the graph a sprout of an element of $\operatorname{2-d-Ring}_{5}$. Category 1 contains all sprouts of seed graphs from $\operatorname{2-d-Ring}_{l_0}$ with $l_0 \in \N\setminus \{2\}$.\\
	(For $l_0=1$ the seed graph would consist of one vertex with two self-loops.)
	
	\item Category 2:\\
	The seed graph of the middle graph in (\ref{Pic_SeedExamples}) has a ring-like structure of length $l_0=4$ and is only traversed in one direction. We thus call the graph a sprout of an element of $\operatorname{1-d-Ring}_{4}$. Category 2 contains all sprouts of seed graphs from $\operatorname{2-d-Ring}_{l_0}$ with even $l_0 \geq 4$.
	
	\item Category 3:\\
	The seed graph of the right most graph in (\ref{Pic_SeedExamples}) has exactly two vertices and two edges in each direction between the two. We call the graph a sprout of an element of $\operatorname{2-d-Ring}_{2}$. Category 3 contains all sprouts of seed graphs from $\operatorname{2-d-Ring}_{2}$.
\end{itemize}
It is easily seen (and proven in Proposition \ref{WeightOfFusedTrees}), that graphs from the first two categories have weight $1$ and graphs from the final category have weight $\E[X_{i \, j}^4]$. For the $r=l$ part of (\ref{Eq_TraceRepresentation4}) we must then count the number of graphs in each of the three categories. Counting number of seed graphs for given lengths $l_0$ and number of black vertices $b$ is simple and done in the proofs of Lemmas \ref{SproutingCorollary1} and \ref{SproutingCorollary2}. We then however still need to know the number of sprouts for each seed graph, which requires a delicate understanding of how the tree structures may be split from the specific seed graph and how the different tree structures may be counted.\\
\\
We will see in Lemma \ref{SproutingBijection} that the number of sprouts to a given seed graph depends only on the number of edges in the seed graph, in the sense, that we describe a simple bijective algorithm $\Phi$ to swap seed graphs, which is color preserving on the tree structures.\\
\\
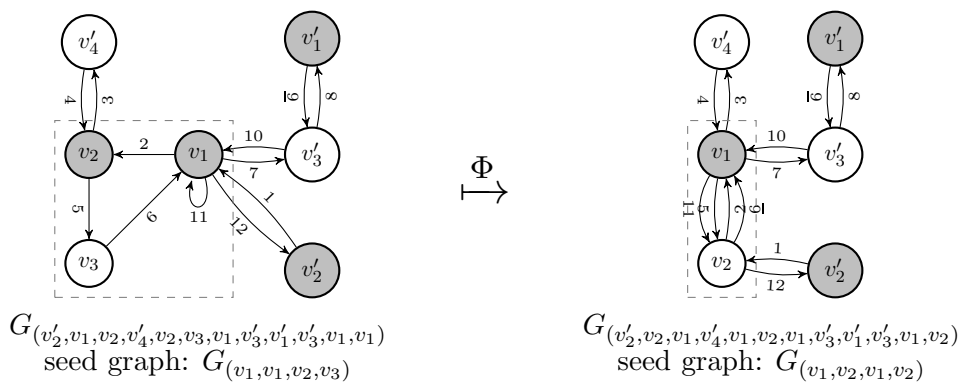
\begin{minipage}{0.45\textwidth} 
	\centering
	\renewcommand\s{0.8}
	\begin{tikzpicture}[node distance=\d and \d,>=stealth',auto, every place/.style={draw}]
		\node [black] (v1) {$v_1$};
		\node [white] (v2) [right=of v1] {$v'_3$};
		\node [black] (v3) [below=of v2] {$v'_2$};
		\node [black] (v4) [above=of v2] {$v'_1$};
		\node [black] (v5) [left=of v1] {$v_2$};
		\node [white] (v6) [below=of v5] {$v_3$};
		\node [white] (v7) [above=of v5] {$v'_4$};
		
		\path[->]
		(v1) edge node[el,above] {$2$} (v5)
		(v5) edge node[el,below] {$5$} (v6)
		(v6) edge node[el,below] {$6$} (v1)
		(v1) edge [in=-105,out=-75,loop] node[el,below] {$11$} (v1)

		(v5) edge [bend right=10] node[el,below] {$3$} (v7)
		(v7) edge [bend right=10] node[el,below] {$4$} (v5)
		
		(v1) edge [bend right=10] node[el,below] {$12$} (v3)
		(v3) edge [bend right=10] node[el,above] {$1$} (v1)
		
		(v1) edge [bend right=10] node[el,below] {$7$} (v2)
		(v2) edge [bend right=10] node[el,above] {$10$} (v1)
		
		(v2) edge [bend right=10] node[el,below] {$8$} (v4)
		(v4) edge [bend right=10] node[el,below] {$\underline{9}$} (v2)
		;
		
		\node[fit=(v1)(v5)(v6), draw, dashed,gray] {};
	\end{tikzpicture}\\
	$G_{(v'_2,v_1,v_2,v'_4,v_2,v_3,v_1,v'_3,v'_1,v'_3,v_1,v_1)}$\\
	seed graph: $G_{(v_1,v_1,v_2,v_3)}$
\end{minipage}
\scalebox{1.5}{$\xmapsto{\Phi}$}
\begin{minipage}{0.45\textwidth} 
	\centering
	\renewcommand\s{0.8}
	\begin{tikzpicture}[node distance=\d and \d,>=stealth',auto, every place/.style={draw}]
		\node [black] (v1) {$v_1$};
		\node [white] (v5) [below=of v1] {$v_2$};
		\node [white] (v2) [right=of v1] {$v'_3$};
		\node [black] (v3) [below=of v2] {$v'_2$};
		\node [black] (v4) [above=of v2] {$v'_1$};
		\node [white] (v6) [above=of v1] {$v'_4$};
		
		\path[->]
		(v1) edge [bend right=10] node[el,below] {$5$} (v5)
		(v5) edge [bend right=10] node[el,below] {$2$} (v1)
		(v1) edge [bend right=30] node[el,below] {$11$} (v5)
		(v5) edge [bend right=30] node[el,below] {$\underline{6}$} (v1)
		
		(v5) edge [bend right=10] node[el,below] {$12$} (v3)
		(v3) edge [bend right=10] node[el,above] {$1$} (v5)
		
		(v1) edge [bend right=10] node[el,below] {$3$} (v6)
		(v6) edge [bend right=10] node[el,below] {$4$} (v1)
		
		(v1) edge [bend right=10] node[el,below] {$7$} (v2)
		(v2) edge [bend right=10] node[el,above] {$10$} (v1)
		
		(v2) edge [bend right=10] node[el,below] {$8$} (v4)
		(v4) edge [bend right=10] node[el,below] {$\underline{9}$} (v2)
		;
		
		\node[fit=(v1)(v5), draw, dashed,gray] {};
	\end{tikzpicture}\\
	$G_{(v'_2,v_2,v_1,v'_4,v_1,v_2,v_1,v'_3,v'_1,v'_3,v_1,v_2)}$\\
	seed graph: $G_{(v_1,v_2,v_1,v_2)}$
\end{minipage}\\
\\
The existence of such a bijection means we do not need to count the number of sprouts for every possible seed-graph from the three aforementioned categories, but only need to check one seed graph with $2l_0$ many edges for every $l_0 \geq 1$. The best choice for a seed graph is $G_0 := G_{\underbrace{(v_1,v_2,v_1,v_2,...,v_1,v_2)}_{l_0 \text{ many pairs}}}$, since every sprout of such a seed graph will automatically be bipartite in the sense that each edge is between one black and one white vertex. We can then apply the B.E.S.T. theorem together with properties of bipartite trees to count the number of sprouts of $G_0$. In Proposition \ref{CountingSprouting} we see that there are $\frac{(l_0+b'+w')!^2}{(l_0+b')!(l_0+w')!}$ many sprouts to a seed graph with $2l_0$ many edges, $b'$ many sprouted black vertices and $w'$ many sprouted white vertices. As we are now able to count the number of graphs in each of the three categories, we can in Section \ref{MainResult1} show
\begin{align*}
	& \sum\limits_{\substack{\bm{i} \in [b]^l, \, \bm{j} \in [l]^l \\ \bm{\{i\}} \cup \bm{\{j\}} = [l] \\ \bm{\{i\}} = [b]}} W(G_{\langle \bm{i},\bm{j} \rangle}) = \frac{b! (l-b)!}{2} \Bigg( {2l \choose 2b} + (2b-1) {l \choose b}^2 \Bigg)\\[-0.5cm]
	& \hspace{4cm} + \big(\E[X_{1 \, 1}^4]-3\big) \, \mathbbm{1}_{b<l} \, b! (l-b)! {l \choose b-1}{l \choose b+1}
\end{align*}
to prove Theorem \ref{Thm_MeanExpansion}.\\
\\
The methods for Theorem \ref{Thm_CovExpansion} are similar and require little additional theory.

\subsection*{Contents}
\begin{enumerate}[label=\arabic*,leftmargin=*,labelsep=2ex,ref=\arabic*]
	\item[\ref{SectionIntroduction}] Introduction \dotfill \pageref{SectionIntroduction}
	\item[\ref{SectionGraphTheory}] Graph theoretical notation and lemmas \dotfill \pageref{SectionGraphTheory}
	\item[\ref{Section_WeightAndColoring}] Graph weight and coloring \dotfill \pageref{Section_WeightAndColoring}
	\item[\ref{SectionCountingSprouting}] Counting sprouts \dotfill \pageref{SectionCountingSprouting}
	\item[\ref{SectionMainRes1}] Proving Theorem \ref{Thm_MeanExpansion} \dotfill \pageref{SectionMainRes1}
	\item[\ref{SectionMainRes2}] Proving Theorem \ref{Thm_CovExpansion} \dotfill \pageref{SectionMainRes2}
	\item[] \hspace{-0.8cm}
	\item[\ref{Section_Adaptability}] Appendix (Adaptability of the methods) \dotfill \pageref{Section_Adaptability}
	\item[\ref{Section_TechnicalLemma}] Appendix (A technical lemma on bipartite trees) \dotfill \pageref{Section_TechnicalLemma}
	\item[] List of symbols \dotfill \pageref{ListOfSymbols}
	\item[] Acknowledgments \& Declarations \dotfill \pageref{Declarations}
	\item[] References \dotfill \pageref{References}
\end{enumerate}

\section{Graph theoretical notation and lemmas}\label{SectionGraphTheory}

\subsection{Definition (Directed multigraph)}\label{DefDirectedMultigraph}
A \textit{directed multigraph} $G$ is a triple $(V,E,f_G)$ consisting of a finite \textit{vertex set} $V$, a finite \textit{edge set} $E$ and a map $f_G: E \rightarrow V \times V$. We say the edge $e \in E$ has \textit{tail} in $\operatorname{tail}(e) = v \in V$ and has \textit{head} in $\operatorname{head}(e) = v' \in V$, if $f_G(e) = (v,v')$.\\
\\
This definition of a multigraph is commonly known as a \textit{multigraph with edges with own identity}, since edges $e,e' \in E$ can be distinct even if they both originate and terminate at the same vertices.\\
\\
A directed multigraph $G = (V,E,f_G)$ will be called \textit{labeled}, if $V$ has a canonical ordering. (We will mostly just be interested in $V \subset \N$.) We can then label these vertices by their ordering, meaning $V$ can be written as $V_r = \{v_1,...,v_r\}$ for some $r \in \N$.\\
\\
A \textit{linearly ordered} directed multigraph $G$ is a labeled directed multigraph $(V_r,E_N,f_G)$, where $E$ also has a canonical ordering. We again label the edges by ordering and write $E$ as $E_N = \{e_1,...,e_N\}$ for some $N \in \N$. We will often use $[N] := \{1,...,N\}$ as $E_N$.\\
\\
For fixed $V_r$ and $E_N$ let $\mathcal{G}_{V_r,E_N}$ denote the set of all linearly ordered directed multigraphs $G=(V_r,E_N,f_G)$.

\subsection{Definition (Visited vertices and exhaustive graphs)}\label{DefVisitedVertex}
For a directed multigraph $G = (V,E,f_G)$ a vertex $v \in V$ is called \textit{visited}, if there exists an edge $e \in E$, for which $v$ is head or tail. Let $V(G)$ denote the set of all visited vertices in $G$. The graph $G$ is called \textit{exhaustive}, if every vertex $v \in V$ is visited, i.e. $V(G) = V$.\\
\\
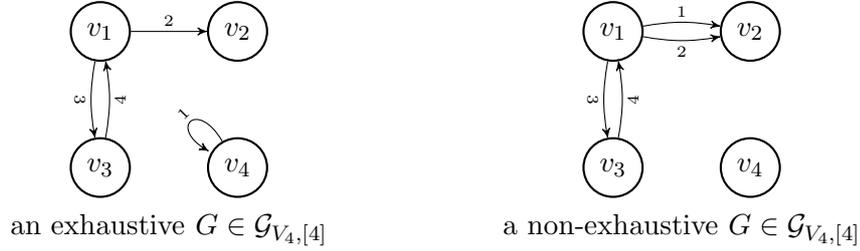
\begin{minipage}{0.45\textwidth}
	\centering
	\begin{tikzpicture}[node distance=\d and \d,>=stealth',auto, every place/.style={draw}]
		\node [white] (v1) {$v_1$};
		\node [white] (v2) [right=of v1] {$v_2$};
		\node [white] (v3) [below=of v1] {$v_3$};
		\node [white] (v4) [below=of v2] {$v_4$};
		
		\path[->]
		(v1) edge [bend right=0] node[el,above] {$2$} (v2)
		(v1) edge [bend right=10] node[el,below] {$3$} (v3)
		(v3) edge [bend right=10] node[el,below] {$4$} (v1)
		
		(v4) edge [in=150,out=120,loop] node[el,above] {$1$} (v4);
		
	\end{tikzpicture}\\
	an exhaustive $G \in \mathcal{G}_{V_4,[4]}$
\end{minipage}
\begin{minipage}{0.45\textwidth}
	\centering
	\begin{tikzpicture}[node distance=\d and \d,>=stealth',auto, every place/.style={draw}]
		\node [white] (v1) {$v_1$};
		\node [white] (v2) [right=of v1] {$v_2$};
		\node [white] (v3) [below=of v1] {$v_3$};
		\node [white] (v4) [below=of v2] {$v_4$};
		
		\path[->]
		(v1) edge [bend right=-10] node[el,above] {$1$} (v2)
		(v1) edge [bend right=10] node[el,below] {$2$} (v2)
		(v1) edge [bend right=10] node[el,below] {$3$} (v3)
		(v3) edge [bend right=10] node[el,below] {$4$} (v1);
		
	\end{tikzpicture}\\
	a non-exhaustive $G \in \mathcal{G}_{V_4,[4]}$
\end{minipage}

\subsection{Definition (Undirected connection)}\label{DefUndirectedConnection}
Let $G = (V,E,f_G)$ be a directed multigraph with more than one vertex. For two different vertices $v, w \in V$ we say there is an \textit{undirected connection} between $v$ and $w$, if there exists an edge $e \in E$ with $f_G(e) \in \{(v,w),(w,v)\}$. The number of such edges does not play a role.\\
Let $\operatorname{U}(G)$ describe the undirected simple graph (possibly with self-loops), which we get by replacing all undirected connections of $G$ with undirected edges.\\
The directed multigraph $G$ is called \textit{undirectedly connected}, if $\operatorname{U}(G)$ is connected.\\
\\
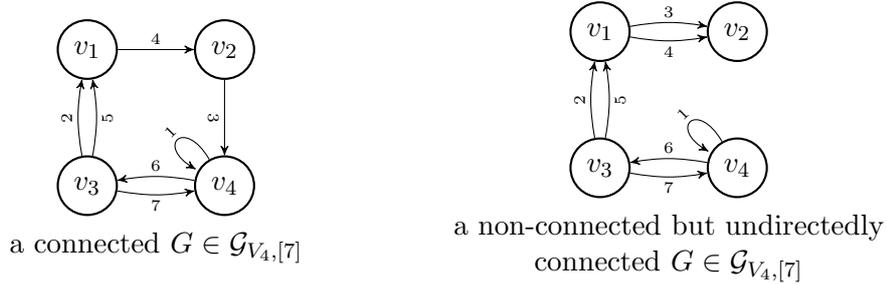
\begin{minipage}{0.45\textwidth}
	\centering
	\begin{tikzpicture}[node distance=\d and \d,>=stealth',auto, every place/.style={draw}]
		\node [white] (v1) {$v_1$};
		\node [white] (v2) [right=of v1] {$v_2$};
		\node [white] (v3) [below=of v1] {$v_3$};
		\node [white] (v4) [below=of v2] {$v_4$};
		
		\path[->]
		(v1) edge [bend right=0] node[el,above] {$4$} (v2)
		(v2) edge [bend right=0] node[el,below] {$3$} (v4)
		(v4) edge [bend right=10] node[el,above] {$6$} (v3)
		(v3) edge [bend right=10] node[el,below] {$7$} (v4)
		(v3) edge [bend right=-10] node[el,above] {$2$} (v1)
		(v3) edge [bend right=10] node[el,below] {$5$} (v1)
		(v4) edge [in=150,out=120,loop] node[el,above] {$1$} (v4);
		
	\end{tikzpicture}\\
	a connected $G \in \mathcal{G}_{V_4,[7]}$
\end{minipage}
\begin{minipage}{0.45\textwidth}
	\centering
	\begin{tikzpicture}[node distance=\d and \d,>=stealth',auto, every place/.style={draw}]
		\node [white] (v1) {$v_1$};
		\node [white] (v2) [right=of v1] {$v_2$};
		\node [white] (v3) [below=of v1] {$v_3$};
		\node [white] (v4) [below=of v2] {$v_4$};
		
		\path[->]
		(v1) edge [bend right=-10] node[el,above] {$3$} (v2)
		(v1) edge [bend right=10] node[el,below] {$4$} (v2)
		(v4) edge [bend right=10] node[el,above] {$6$} (v3)
		(v3) edge [bend right=10] node[el,below] {$7$} (v4)
		(v3) edge [bend right=-10] node[el,above] {$2$} (v1)
		(v3) edge [bend right=10] node[el,below] {$5$} (v1)
		(v4) edge [in=150,out=120,loop] node[el,above] {$1$} (v4);
		
	\end{tikzpicture}\\
	a non-connected but undirectedly connected $G \in \mathcal{G}_{V_4,[7]}$
\end{minipage}

\subsection{Definition (Routes and circuit graphs)}\label{DefCircuitGraph}
For any fixed vertex set $V_r$ a \textit{route} through $V_r$ of length $N$ is defined to be a sequence $\bm{i} \subset V_r^N$ with the property
\begin{align*}
	& V_r = \bm{\{i\}} := \text{set of all vertices occurring in } \bm{i} \ .
\end{align*}
For a given route $\bm{i} \subset V_r^N$ let $E_N := [N]$, then we define a linearly ordered directed multigraph $G_{\bm{i}} = (V_r,[N],f_{\bm{i}}) \in \mathcal{G}_{V_r,E_N}$ by
\begin{align*}
	& f_{\bm{i}}(\underbrace{k}_{=e_k}) := \begin{cases}
		(i_k,i_{k+1}) & \text{, if } k < N\\
		(i_N,i_1) & \text{, if } k = N
	\end{cases} \ .
\end{align*}
The set of \textit{circuit (multi-)graphs} on $V_r$ of length $N$ is then defined as
\begin{align*}
	& \mathcal{C}_{V_r,N} := \{G_{\bm{i}} \mid \bm{i} \subset V_r^N \text{ route through $V_r$ of length $N$} \} \ .
\end{align*}
Such graphs are by construction exhaustive and connected.\\
\\
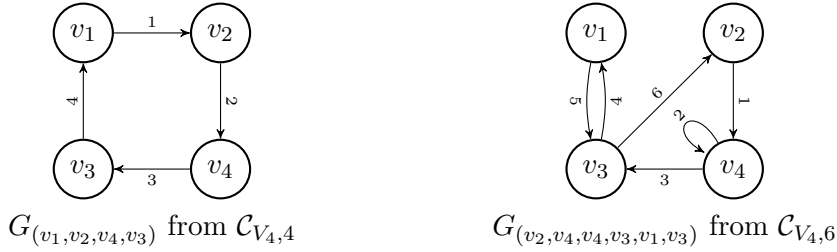
\begin{minipage}{0.45\textwidth}
	\centering
	\begin{tikzpicture}[node distance=\d and \d,>=stealth',auto, every place/.style={draw}]
		\node [white] (v1) {$v_1$};
		\node [white] (v2) [right=of v1] {$v_2$};
		\node [white] (v3) [below=of v1] {$v_3$};
		\node [white] (v4) [below=of v2] {$v_4$};
		
		\path[->]
		(v1) edge [bend right=0] node[el,above] {$1$} (v2)
		(v2) edge [bend right=0] node[el,above] {$2$} (v4)
		(v4) edge [bend right=0] node[el,below] {$3$} (v3)
		(v3) edge [bend right=0] node[el,above] {$4$} (v1);
		
	\end{tikzpicture}\\
	$G_{(v_1,v_2,v_4,v_3)}$ from $\mathcal{C}_{V_4,4}$
\end{minipage}
\begin{minipage}{0.45\textwidth}
	\centering
	\begin{tikzpicture}[node distance=\d and \d,>=stealth',auto, every place/.style={draw}]
		\node [white] (v1) {$v_1$};
		\node [white] (v2) [right=of v1] {$v_2$};
		\node [white] (v3) [below=of v1] {$v_3$};
		\node [white] (v4) [below=of v2] {$v_4$};
		
		\path[->]
		(v2) edge [bend right=0] node[el,above] {$1$} (v4)
		(v4) edge [in=150,out=120,loop] node[el,above] {$2$} (v4)
		(v4) edge [bend right=0] node[el,below] {$3$} (v3)
		(v3) edge [bend right=10] node[el,below] {$4$} (v1)
		(v1) edge [bend right=10] node[el,below] {$5$} (v3)
		(v3) edge [bend right=0] node[el,above] {$6$} (v2);
		
	\end{tikzpicture}\\
	$G_{(v_2,v_4,v_4,v_3,v_1,v_3)}$ from $\mathcal{C}_{V_4,6}$
\end{minipage}


\subsection{Definition (Reversal operator)}\label{DefReversed}
For any $r,N \in \N$ define the \textit{reversal operator}
\begin{align*}
	& \operatorname{R}: \mathcal{C}_{V_r,N} \rightarrow \mathcal{G}_{V_r,[N]}
\end{align*}
by reversing the direction of every second edge $e_k \in [N] = \{1,...,N\}$. In other words if $G=(V_r,[N],f) \in \mathcal{C}_{V_r,N}$, then $\operatorname{R}(G)=(V_r,[N],f')$, where $f'$ is given by
\begin{align*}
	& f'(k) = 
	\begin{cases}
		(v_i,v_j) & \text{ for } f(k) = (v_i,v_j) \text{, if $k$ is odd}\\
		(v_j,v_i) & \text{ for } f(k) = (v_i,v_j) \text{, if $k$ is even}
	\end{cases} \ .
\end{align*}

\subsection{Definition (Balanced directed multigraphs)}\label{DefBalanced}
We call a directed multigraph $G=(V,E,f_G)$ \textit{balanced}, if the edges $E$ can be split into \textit{(balanced) edge pairs} $(e,e')$ such that the head of edge $e$ is the tail of edge $e'$ and vice versa.\\
\\
If $G$ is labeled, an equivalent definition would be to say that its adjacency matrix $A(G)$ is symmetric and has only even entries on the diagonal.

\subsection{Definition (Balanced tree)}\label{DefBalancedTree}
For any $l \in \N$ and $V_{l+1}$ the set
$$
\mathcal{T}_{V_{l+1}} := \{G \in \mathcal{C}_{V_{l+1},2l} \mid G \text{ is balanced} \}
$$
will be called the set of \textit{balanced trees} of length $l$.\\
\\
Since elements $G \in \mathcal{C}_{V_{l+1},2l}$ are by construction exhaustive and connected, we know $\operatorname{U}(G)$ to be connected with $l+1$ vertices. As $\operatorname{U}(G)$ can have at most $l$ many edges, it is an elementary exercise in graph theory to see that $\operatorname{U}(G)$ can have no cycles and must have exactly $l$ edges, i.e. be a tree. Also each undirected edge of $\operatorname{U}(G)$ must correspond to one balanced edge pair in $G$.\\
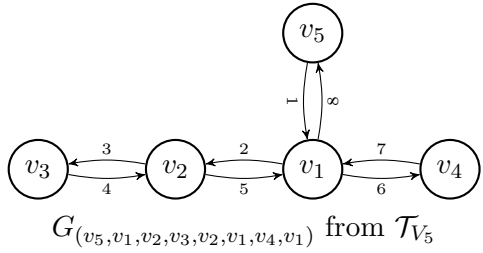
\begin{minipage}{0.45\textwidth}
	\centering
	\begin{tikzpicture}[node distance=\d and \d,>=stealth',auto, every place/.style={draw}]
		\node [white] (v1) {$v_1$};
		\node [white] (v4) [right=of v1] {$v_4$};
		\node [white] (v5) [above=of v1] {$v_5$};
		\node [white] (v2) [left=of v1] {$v_2$};
		\node [white] (v3) [left=of v2] {$v_3$};
		
		\path[->]
		(v1) edge [bend right=10] node[el,above] {$2$} (v2)
		(v2) edge [bend right=10] node[el,below] {$5$} (v1)
		
		(v1) edge [bend right=10] node[el,below] {$8$} (v5)
		(v5) edge [bend right=10] node[el,below] {$1$} (v1)
		
		(v1) edge [bend right=10] node[el,below] {$6$} (v4)
		(v4) edge [bend right=10] node[el,above] {$7$} (v1)
		
		(v2) edge [bend right=10] node[el,above] {$3$} (v3)
		(v3) edge [bend right=10] node[el,below] {$4$} (v2);
		
	\end{tikzpicture}\\
	$G_{(v_5,v_1,v_2,v_3,v_2,v_1,v_4,v_1)}$ from $\mathcal{T}_{V_{5}}$
\end{minipage}
\begin{minipage}{0.5\textwidth}
	
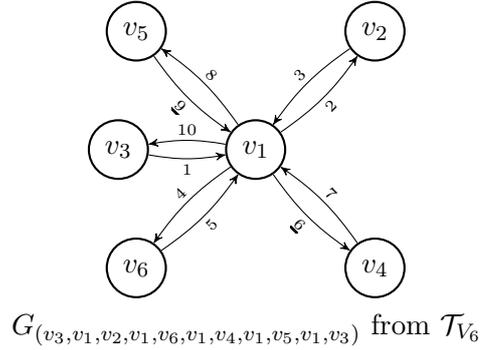
\centering
	\begin{tikzpicture}[node distance=\d and \d,>=stealth',auto, every place/.style={draw}]
		\node [white] (v1) {$v_1$};
		\node [white] (v2) [above right=of v1] {$v_2$};
		\node [white] (v3) [left=of v1] {$v_3$};
		\node [white] (v5) [above left=of v1] {$v_5$};
		\node [white] (v6) [below left=of v1] {$v_6$};
		\node [white] (v4) [below right=of v1] {$v_4$};
		
		\path[->]
		(v1) edge [bend right=10] node[el,below] {$2$} (v2)
		(v2) edge [bend right=10] node[el,above] {$3$} (v1)
		
		(v1) edge [bend right=10] node[el,above] {$10$} (v3)
		(v3) edge [bend right=10] node[el,below] {$1$} (v1)
		
		(v1) edge [bend right=10] node[el,below] {$\underline{6}$} (v4)
		(v4) edge [bend right=10] node[el,above] {$7$} (v1)
		
		(v1) edge [bend right=10] node[el,above] {$8$} (v5)
		(v5) edge [bend right=10] node[el,below] {$\underline{9}$} (v1)
		
		(v1) edge [bend right=10] node[el,above] {$4$} (v6)
		(v6) edge [bend right=10] node[el,below] {$5$} (v1);
		
	\end{tikzpicture}\\
	$G_{(v_3,v_1,v_2,v_1,v_6,v_1,v_4,v_1,v_5,v_1,v_3)}$ from $\mathcal{T}_{V_{6}}$
\end{minipage}

\subsection{Definition (Ring-type graphs)}\label{DefRing}
For $l_0 \in \N$ and a vertex set $V_{l_0}$ a circuit graph $G = G_{\bm{i}}$ from $\mathcal{C}_{V_{l_0},2l_0}$ will be called \textit{of ring-type}, if $\operatorname{U}(G)$ is a cycle graph and each undirected connection in $G$ consists of two edges. If $G$ is additionally balanced, it will be called a \textit{two-directional} ring-type graph. Otherwise it will be called a \textit{one-directional} ring type graph. Note that a one-directional ring type graph is only possible for $l_0 \geq 3$.\\
\\
For any $G \in \mathcal{C}_{V_{l_0},2l_0}$ we (with slight abuse of notation as $V_{l_0}$ is lost) write $G \in \operatorname{1-d-Ring}_{l_0}$, if $G$ is a one-directional ring-type graph with ring length $l_0$. Analogously we write $G \in \operatorname{2-d-Ring}_{l_0}$, if $G$ is a two-directional ring-type graph with ring length $l_0$.\\

\renewcommand\s{1}
\begin{minipage}{0.4\textwidth} 
	\centering
	\begin{tikzpicture}[node distance=\d and \d,>=stealth',auto, every place/.style={draw}]
		\node [white] (v1) {$v_1$};
		
		\path (v1) edge [in=150,out=120,loop] node[el,above] {$1$} (v1);
		\path (v1) edge [in=100,out=70,loop] node[el,above] {$2$} (v1);
	\end{tikzpicture}\\
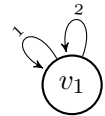
	$G_{(v_1,v_1)}$, the only\\
	element in $\operatorname{2-d-Ring}_{1}$
\end{minipage}
\begin{minipage}{0.4\textwidth} 
	\centering
	\begin{tikzpicture}[node distance=\d and \d,>=stealth',auto, every place/.style={draw}]
		\node [white] (v1) {$v_1$};
		\node [white] (v2) [below=of v1] {$v_2$};
		
		\path[->]
		(v1) edge [bend right=20] node[el,above] {$1$} (v2)
		(v2) edge [bend right=20] node[el,above] {$2$} (v1)
		(v1) edge [bend right=30] node[el,below] {$3$} (v2)
		(v2) edge [bend right=30] node[el,below] {$4$} (v1);
		
		\node [white] (v3) [right=of v1] {$v_1$};
		\node [white] (v4) [right=of v2] {$v_2$};
		
		\path[->]
		(v4) edge [bend right=20] node[el,above] {$1$} (v3)
		(v3) edge [bend right=20] node[el,above] {$2$} (v4)
		(v4) edge [bend right=30] node[el,below] {$3$} (v3)
		(v3) edge [bend right=30] node[el,below] {$4$} (v4);
	\end{tikzpicture}\\
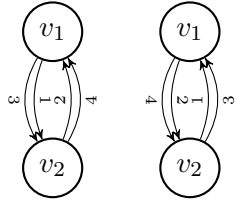
	$G_{(v_1,v_2,v_1,v_2)}$ and $G_{(v_2,v_1,v_2,v_1)}$,\\
	the only elements in $\operatorname{2-d-Ring}_{2}$
\end{minipage}
\\
\\
\begin{minipage}{0.5\textwidth} 
	\centering
	\begin{tikzpicture}[node distance=\d and \d,>=stealth',auto, every place/.style={draw}]
		\node [white] (v1) {$v_1$};
		\node [white] (v2) [right=of v1] {$v_4$};
		\node [white] (v3) [below right=of v2] {$v_3$};
		\node [white] (v4) [below left=of v3] {$v_2$};
		\node [white] (v5) [left=of v4] {$v_6$};
		\node [white] (v6) [above left=of v5] {$v_5$};
		\path[->] (v1) edge [bend right=10] node[el,below] {$2$} (v2);
		\path[->] (v1) edge [bend right=-10] node[el,above] {$8$} (v2);
		
		\path[->] (v2) edge [bend right=10] node[el,below] {$3$} (v3);
		\path[->] (v2) edge [bend right=-10] node[el,above] {$9$} (v3);
		
		\path[->] (v3) edge [bend right=10] node[el,above] {$4$} (v4);
		\path[->] (v3) edge [bend right=-10] node[el,below] {$10$} (v4);
		
		\path[->] (v4) edge [bend right=10] node[el,above] {$5$} (v5);
		\path[->] (v4) edge [bend right=-10] node[el,below] {$11$} (v5);
		
		\path[->] (v5) edge [bend right=10] node[el,above] {$6$} (v6);
		\path[->] (v5) edge [bend left=10] node[el,below] {$12$} (v6);
		
		\path[->] (v6) edge [bend right=10] node[el,below] {$1$} (v1);
		\path[->] (v6) edge [bend right=-10] node[el,above] {$7$} (v1);
		
	\end{tikzpicture}\\
	$G_{(v_5,v_1,v_4,v_3,v_2,v_6,v_5,v_1,v_4,v_3,v_2,v_6)}$\\
	from $\operatorname{1-d-Ring}_{6}$
\end{minipage}
\begin{minipage}{0.3\textwidth} 
	\centering
	\begin{tikzpicture}[node distance=\d and \d,>=stealth',auto, every place/.style={draw}]
		\node [white] (v1) {$v_2$};
		\node [white] (v2) [right=of v1] {$v_4$};
		\node [white] (v3) [below=of v2] {$v_1$};
		\node [white] (v4) [left=of v3] {$v_3$};
		
		\path[->] (v1) edge [bend right=10] node[el,below] {$1$} (v2);
		\path[->] (v2) edge [bend right=10] node[el,above] {$6$} (v1);
		
		\path[->] (v2) edge [bend right=10] node[el,below] {$2$} (v3);
		\path[->] (v3) edge [bend right=10] node[el,below] {$5$} (v2);
		
		\path[->] (v3) edge [bend right=10] node[el,above] {$3$} (v4);
		\path[->] (v4) edge [bend right=10] node[el,below] {$4$} (v3);
		
		\path[->] (v4) edge [bend right=10] node[el,below] {$8$} (v1);
		\path[->] (v1) edge [bend right=10] node[el,below] {$7$} (v4);
		
	\end{tikzpicture}\\
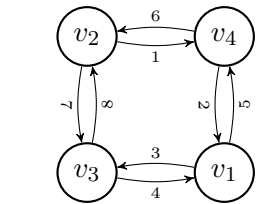
	$G_{(v_2,v_4,v_1,v_3,v_1,v_4,v_2,v_3,v_2)}$\\
	from $\operatorname{2-d-Ring}_{4}$
\end{minipage}

\subsection{Definition (Balanced leaves)}
Let $G = (V,E,f)$ be a directed multigraph with more than one vertex. We call a vertex $v \in V$ a \textit{balanced leaf} of $G$, if there is exactly one edge with head $v$, one other edge with tail $v$ and both these edges are between $v$ and one other vertex $w \in V$. Accordingly, in circuit graphs $G=G_{\bm{i}} \in \mathcal{C}_{V_r,N}$ a vertex $v_j$ is a balanced leaf, if and only if $v_j$ occurs only once in $\bm{i}$ and this occurrence's left- and right-hand neighbors in $\bm{i}$ are equal. If the sequence $\bm{i}$ starts or ends in $v_j$, we loop around the ends of $\bm{i}$ to find these neighbors.

\subsection{Definition (Removing balanced leaves)}\label{DefRemovingLeaves}
For any $G_{\bm{i}} \in \mathcal{C}_{V_r,N}$ with a balanced leaf $v_j \in V_{r}$ and $N>2$ let $v_s$ be the only neighbor of $v_j$ in the graph $G_{\bm{i}}$, then the sub-sequence $(v_s,v_j,v_s)$ must occur in $\bm{i}$, though it might happen that the sub-sequence is interrupted by the end of the route, in which case the sub-sequence will continue at the beginning of the route. We define the modified route
$$
\bm{i}' = \left.
\begin{cases}
	(v_{i_1},...,v_{s},\widehat{v_j},\widehat{v_s},...,v_{i_{2l}}) & \text{, if } v_{2l} \neq v_j\\
	(v_{s},v_{i_2},...,v_{i_{2l-2}},\widehat{v_{s}},\widehat{v_j}) & \text{, if } v_{2l} = v_j
\end{cases} \right\}
\in \{v_1,...,\widehat{v_j},...,v_r\}^{N-2}
$$
by ignoring the singular occurrence of $v_j$ in $\bm{i}$ together with the next entry $v_s$ of $\bm{i}$, if it exists. If the next entry does not exists, then $v_j$ must be the last entry in $\bm{i}$ and we instead ignore the entry $v_s$ previous to the occurrence of $v_j$. In both cases we have ignored the occurrence of $v_j$ and an occurrence of $v_s$. (This definition guarantees that in both cases the positions of the remaining entries stay the same modulo $2$.) Let $e,e' \in E_{N}$ be the two edges between $v_j$ and $v_s$ in $G_{\bm{i}}$, then we call the graph
\begin{align*}
	& \widetilde{G} := G_{\bm{i}'} \in \mathcal{C}_{\widetilde{V}_{r-1} := V_{r}\setminus\{v_j\},N-2}
\end{align*}
the \textit{version of $G_{\bm{i}}$ with $v_j$ removed}.\\
\\
We had assumed $N>2$ in order to guarantee that $\widetilde{G}$ still has edges (and can thus have visited vertices). For $N \leq 2$, we say that $G_{\bm{i}}$ has no balanced leaves.

\subsection{Definition (Seed graph)}\label{DefSeed}
For any $G \in \mathcal{C}_{V_r,2l}$ let the \textit{seed graph} $\operatorname{S}(G) \in \mathcal{C}_{V_{r_0},2l_0}$ (for certain $r_0,l_0$ with $l_0-r_0 = l-r$) be given by the following recursive definition. If $G$ has no balanced leaves (includes $l=1$), we define $\operatorname{S}(G) := G$. Otherwise $G$ let $v_j$ be the 'smallest' (equivalently 'lowest indexed') balanced leaf of $G$ and $\widetilde{G}$ be the version of $G$ with $v_j$ removed. We recursively define $\operatorname{S}(G) := \operatorname{S}(\widetilde{G})$ to be the seed graph of $\widetilde{G}$.\\
\\
In the opposite direction we say $G$ is a \textit{sprout} of $\operatorname{S}(G)$. Also the vertices in $\operatorname{S}(G)$ are called \textit{seed vertices} and vertices from $G \setminus \operatorname{S}(G)$ are called \textit{sprouted vertices}.\\

\begin{minipage}{0.45\textwidth}
	\centering
	\begin{tikzpicture}[node distance=\d and \d,>=stealth',auto, every place/.style={draw}]
		\node [white] (v1) {$v_1$};
		\node [white] (v3) [right=of v1] {$v_3$};
		\node [white] (v4) [right=of v3] {$v_4$};
		\node [white] (v2) [above=of v4] {$v_2$};
		
		\path[->]
		(v2) edge [bend right=0] node[el,above] {$1$} (v4)
		(v4) edge [in=150,out=120,loop] node[el,above] {$2$} (v4)
		(v4) edge [bend right=0] node[el,below] {$3$} (v3)
		(v3) edge [bend right=10] node[el,above] {$4$} (v1)
		(v1) edge [bend right=10] node[el,below] {$5$} (v3)
		(v3) edge [bend right=0] node[el,above] {$6$} (v2);
		
		\node[fit=(v2)(v3)(v4), draw, dashed,blue] {};
		
	\end{tikzpicture}\\
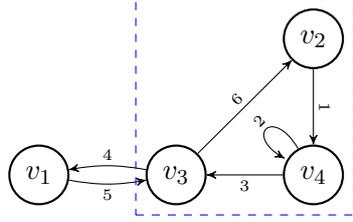
	$G_{(v_2,v_4,v_4,v_3,v_1,v_3)}$ from $\mathcal{C}_{V_4,6}$\\
	with \textcolor{blue}{$\operatorname{S}(G_{...}) = G_{(v_2,v_4,v_4,v_3)}$}
\end{minipage}
\begin{minipage}{0.45\textwidth}
	\centering
	\begin{tikzpicture}[node distance=\d and \d,>=stealth',auto, every place/.style={draw}]
		\node [white] (v1) {$v_1$};
		\node [white] (v2) [below=of v1] {$v_2$};
		\node [white] (v3) [left=of v1] {$v_3$};
		\node [white] (v4) [right=of v2] {$v_4$};
		\node [white] (v5) [right=of v1] {$v_5$};

		\path[->]
		(v1) edge [bend right=20] node[el,above] {$2$} (v2)
		(v2) edge [bend right=20] node[el,above] {$3$} (v1)
		(v1) edge [bend right=30] node[el,below] {$\underline{6}$} (v2)
		(v2) edge [bend right=30] node[el,below] {$\underline{9}$} (v1)

		(v3) edge [bend right=10] node[el,below] {$1$} (v1)
		(v1) edge [bend right=10] node[el,above] {$10$} (v3)
		
		(v1) edge [bend right=10] node[el,below] {$4$} (v5)
		(v5) edge [bend right=10] node[el,above] {$5$} (v1)
		
		(v2) edge [bend right=10] node[el,below] {$7$} (v4)
		(v4) edge [bend right=10] node[el,above] {$8$} (v2)
		;
		
		\node[fit=(v1)(v2), draw, dashed,blue] {};
		
	\end{tikzpicture}\\
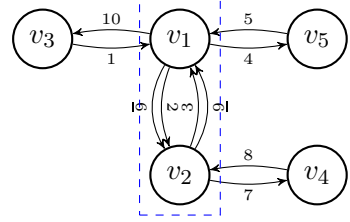
	$G_{(v_3,v_1,v_2,v_1,v_5,v_1,v_2,v_4,v_2,v_1)}$ from $\mathcal{C}_{V_5,10}$\\
	with \textcolor{blue}{$\operatorname{S}(G_{...}) = G_{(v_2,v_1,v_2,v_1)}$}
\end{minipage}

\subsection{Remark (Seeds of trees)}\label{SeedsOfTrees}
By properties of the well known Prüfer-code algorithm, balanced trees are precisely the elements $G \in \mathcal{C}_{V_{l+1},2l}$ where the seed graph consist of two vertices connected by a balanced edge pair. Balanced trees are also the only type of circuit graph for which the seed graph depends on the ordering of $V_r$. To avoid problems, that may stem from this, we will in Proposition \ref{CountingSprouting} only look at balanced trees with an edge pair between the two largest (highest indexed) vertices. This guarantees that their seed graph will consist of these two largest vertices.\\
In later applications we will not be examining balanced trees and we will not need to address these problems by making requirements to the order of the vertices.\\

\begin{minipage}{0.45\textwidth}
	\centering
	\begin{tikzpicture}[node distance=\d and \d,>=stealth',auto, every place/.style={draw}]
		\node [white] (v1) {$v_1$};
		\node [white] (v4) [right=of v1] {$v_4$};
		\node [white] (v5) [above=of v1] {$v_5$};
		\node [white] (v2) [left=of v1] {$v_2$};
		\node [white] (v3) [left=of v2] {$v_3$};
		
		\path[->]
		(v1) edge [bend right=10] node[el,above] {$2$} (v2)
		(v2) edge [bend right=10] node[el,below] {$5$} (v1)
		
		(v1) edge [bend right=10] node[el,below] {$8$} (v5)
		(v5) edge [bend right=10] node[el,below] {$1$} (v1)
		
		(v1) edge [bend right=10] node[el,below] {$6$} (v4)
		(v4) edge [bend right=10] node[el,above] {$7$} (v1)
		
		(v2) edge [bend right=10] node[el,above] {$3$} (v3)
		(v3) edge [bend right=10] node[el,below] {$4$} (v2);
		
		\node[fit=(v1)(v5), draw, dashed,blue] {};
		
	\end{tikzpicture}\\
	$G_{(v_5,v_1,v_2,v_3,v_2,v_1,v_4,v_1)}$ from $\mathcal{T}_{V_{5}}$\\
	with $\textcolor{blue}{\operatorname{S}(G_{...}) = G_{(v_5,v_1)}}$
\end{minipage}
\begin{minipage}{0.45\textwidth}
	\centering
	\begin{tikzpicture}[node distance=\d and \d,>=stealth',auto, every place/.style={draw}]
		\node [white] (v1) {$v_1$};
		\node [white] (v4) [right=of v1] {$v_4$};
		\node [white] (v5) [above=of v1] {$v_3$};
		\node [white] (v2) [left=of v1] {$v_2$};
		\node [white] (v3) [left=of v2] {$v_5$};
		
		\path[->]
		(v1) edge [bend right=10] node[el,above] {$2$} (v2)
		(v2) edge [bend right=10] node[el,below] {$5$} (v1)
		
		(v1) edge [bend right=10] node[el,below] {$8$} (v5)
		(v5) edge [bend right=10] node[el,below] {$1$} (v1)
		
		(v1) edge [bend right=10] node[el,below] {$6$} (v4)
		(v4) edge [bend right=10] node[el,above] {$7$} (v1)
		
		(v2) edge [bend right=10] node[el,above] {$3$} (v3)
		(v3) edge [bend right=10] node[el,below] {$4$} (v2);
		
		\node[fit=(v2)(v3), draw, dashed,blue] {};
		
	\end{tikzpicture}\\
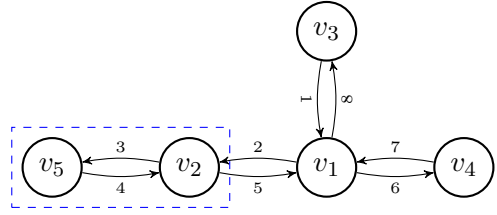
	$G_{(v_3,v_1,v_2,v_5,v_2,v_1,v_4,v_1)}$ from $\mathcal{T}_{V_{5}}$\\
	with $\textcolor{blue}{\operatorname{S}(G_{...}) = G_{(v_2,v_5)}}$
\end{minipage}


\section{Graph weight and coloring}\label{Section_WeightAndColoring}

\subsection{Definition (Black-white coloring of circuit graphs)}\label{DefBlackWhiteCol}
For any $G \in \mathcal{C}_{V_r,N}$ we call a vertex $v_i \in V_r$ \textit{white}, if it is only tail of even numbered edges $e_k=k \in [N]$. Due to $G$ being a circuit graph, this is equivalent to $v_i$ only being head of odd numbered edges. We call a vertex of $G$ \textit{black}, if it is not white.\\
\\
In terms of the route $\bm{i}_G$ this is equivalent to a vertex $v_j \in V_r$ being white, iff $j$ only appears in even numbered entries of $\bm{i}_G \in V_r^N$.\\
\\
Let $B(G) \subset V_r$ denote the set of all black vertices. By construction we have
\begin{align*}
	& B(G_{\bm{i}}) := \{v_{i_1},v_{i_3},...\} = \{v_{i_t} \mid t \leq N \text{ odd}\} \ .
\end{align*}\
\\
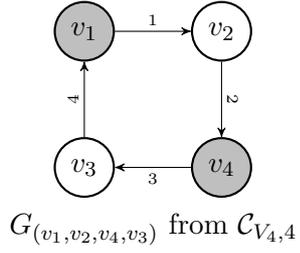
\begin{minipage}{0.45\textwidth}
	\centering
	\begin{tikzpicture}[node distance=\d and \d,>=stealth',auto, every place/.style={draw}]
		\node [black] (v1) {$v_1$};
		\node [white] (v2) [right=of v1] {$v_2$};
		\node [white] (v3) [below=of v1] {$v_3$};
		\node [black] (v4) [below=of v2] {$v_4$};
		
		\path[->]
		(v1) edge [bend right=0] node[el,above] {$1$} (v2)
		(v2) edge [bend right=0] node[el,above] {$2$} (v4)
		(v4) edge [bend right=0] node[el,below] {$3$} (v3)
		(v3) edge [bend right=0] node[el,above] {$4$} (v1);
		
	\end{tikzpicture}\\
	$G_{(v_1,v_2,v_4,v_3)}$ from $\mathcal{C}_{V_4,4}$
\end{minipage}
\begin{minipage}{0.45\textwidth}
	
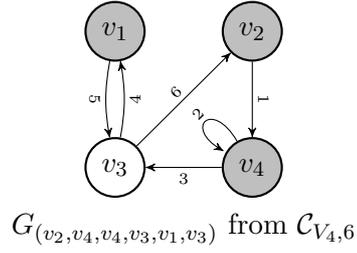
\centering
	\begin{tikzpicture}[node distance=\d and \d,>=stealth',auto, every place/.style={draw}]
		\node [black] (v1) {$v_1$};
		\node [black] (v2) [right=of v1] {$v_2$};
		\node [white] (v3) [below=of v1] {$v_3$};
		\node [black] (v4) [below=of v2] {$v_4$};
		
		\path[->]
		(v2) edge [bend right=0] node[el,above] {$1$} (v4)
		(v4) edge [in=150,out=120,loop] node[el,above] {$2$} (v4)
		(v4) edge [bend right=0] node[el,below] {$3$} (v3)
		(v3) edge [bend right=10] node[el,below] {$4$} (v1)
		(v1) edge [bend right=10] node[el,below] {$5$} (v3)
		(v3) edge [bend right=0] node[el,above] {$6$} (v2);
		
	\end{tikzpicture}\\
	$G_{(v_2,v_4,v_4,v_3,v_1,v_3)}$ from $\mathcal{C}_{V_4,6}$
\end{minipage}
\\
\\
\\
The above definition is easily extended to non-exhaustive circuit graphs by calling all unvisited vertices white. This will only become necessary in Definition \ref{DefDoubleCircuit}.

\subsection{Lemma (Removing balanced leaves does not change the coloring)}\label{RemovingLeaves_Coloring}
For any $G_{\bm{i}} \in \mathcal{C}_{V_r,E_{2l}}$ with balanced leaf $v_j$ (of arbitrary coloring) let $\widetilde{G} \in \mathcal{C}_{\widetilde{V}_{r-1},2l-2}$ be the version of $G_{\bm{i}}$ with $v_j$ removed as in Definition \ref{DefRemovingLeaves}. The property
\begin{align*}
	& B(G) \setminus \{v_j\} = B(\widetilde{G})
\end{align*}
holds.
\begin{proof}\
	\\
	By construction of the route $\bm{i}'$ in Definition \ref{DefRemovingLeaves} the positions of the entries in $\bm{i}'$ are the same as their positions in $\bm{i}$ modulo $2$. We thus have
	\begin{align*}
		& B(G_{\bm{i}}) \setminus \{v_j\} = \{v_{i_1},v_{i_3},...,v_{i_{2l-1}}\} \setminus \{v_j\} = \{v_{i_1'},v_{i_3'},...,v_{i'_{2l-3}}\} = B(\widetilde{G}) \ . \qedhere
	\end{align*}
\end{proof}

\subsection{Lemma (Coloring of balanced trees)}\label{TreeColoring}
For any balanced tree $G_{\bm{i}} \in \mathcal{T}_{V_{l+1}}$ each edge must be between a black and a white vertex.\\
\\
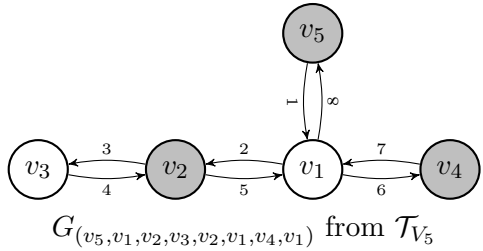
\begin{minipage}{0.45\textwidth}
	\centering
	\begin{tikzpicture}[node distance=\d and \d,>=stealth',auto, every place/.style={draw}]
		\node [white] (v1) {$v_1$};
		\node [black] (v4) [right=of v1] {$v_4$};
		\node [black] (v5) [above=of v1] {$v_5$};
		\node [black] (v2) [left=of v1] {$v_2$};
		\node [white] (v3) [left=of v2] {$v_3$};
		
		\path[->]
		(v1) edge [bend right=10] node[el,above] {$2$} (v2)
		(v2) edge [bend right=10] node[el,below] {$5$} (v1)
		
		(v1) edge [bend right=10] node[el,below] {$8$} (v5)
		(v5) edge [bend right=10] node[el,below] {$1$} (v1)
		
		(v1) edge [bend right=10] node[el,below] {$6$} (v4)
		(v4) edge [bend right=10] node[el,above] {$7$} (v1)
		
		(v2) edge [bend right=10] node[el,above] {$3$} (v3)
		(v3) edge [bend right=10] node[el,below] {$4$} (v2);
		
	\end{tikzpicture}\\
	$G_{(v_5,v_1,v_2,v_3,v_2,v_1,v_4,v_1)}$ from $\mathcal{T}_{V_{5}}$
\end{minipage}
\begin{minipage}{0.5\textwidth}
	
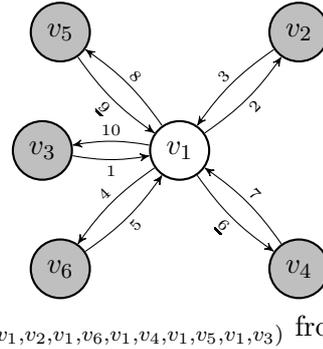
\centering
	\begin{tikzpicture}[node distance=\d and \d,>=stealth',auto, every place/.style={draw}]
		\node [white] (v1) {$v_1$};
		\node [black] (v2) [above right=of v1] {$v_2$};
		\node [black] (v3) [left=of v1] {$v_3$};
		\node [black] (v5) [above left=of v1] {$v_5$};
		\node [black] (v6) [below left=of v1] {$v_6$};
		\node [black] (v4) [below right=of v1] {$v_4$};
		
		\path[->]
		(v1) edge [bend right=10] node[el,below] {$2$} (v2)
		(v2) edge [bend right=10] node[el,above] {$3$} (v1)
		
		(v1) edge [bend right=10] node[el,above] {$10$} (v3)
		(v3) edge [bend right=10] node[el,below] {$1$} (v1)
		
		(v1) edge [bend right=10] node[el,below] {$\underline{6}$} (v4)
		(v4) edge [bend right=10] node[el,above] {$7$} (v1)
		
		(v1) edge [bend right=10] node[el,above] {$8$} (v5)
		(v5) edge [bend right=10] node[el,below] {$\underline{9}$} (v1)
		
		(v1) edge [bend right=10] node[el,above] {$4$} (v6)
		(v6) edge [bend right=10] node[el,below] {$5$} (v1);
		
	\end{tikzpicture}\\
	$G_{(v_3,v_1,v_2,v_1,v_6,v_1,v_4,v_1,v_5,v_1,v_3)}$ from $\mathcal{T}_{V_{6}}$
\end{minipage}\
\\
\\
This property is easily seen to be true. The proof writes itself by iterative removal of leaves using Lemma \ref{RemovingLeaves_Coloring}.

\subsection{Lemma (Only balanced trees have $l+1$ vertices and positive weight)}\label{WeightOfTrees}
For any vertex set $V_{l+1}$ and $\bm{i},\bm{k} \in V_{l+1}^l$ with $\bm{\{i\}} \cup \bm{\{k\}} = V_{l+1}$ let $G_{\left< \bm{i},\bm{k} \right>}$ be the circuit graph in $\mathcal{C}_{V_{l+1},2l}$ with route $\left< \bm{i},\bm{k} \right>$, then
\begin{align}\label{Eq_WeightOfTrees_1}
	\cW(G_{\left< \bm{i}, \bm{k} \right>}) := & \E\left[ (X_{i_1 \, k_1}X_{i_2 \, k_1}) ... (X_{i_{l-1} \, k_{l-1}}X_{i_l \, k_{l-1}}) \, (X_{i_l \, k_l}X_{i_1 \, k_l}) \right]\\
	 = & \mathbbm{1}_{G_{\left< \bm{i}, \bm{k} \right>} \in \mathcal{T}_{V_{l+1}}} \ .
\end{align}
\begin{proof}\
	\\
	The product form (\ref{Eq_Weight2}) yields $\cW(G_{\left< \bm{i}, \bm{k} \right>}) = \prod\limits_{i,j=1}^l \E\Big[X_{i \, j}^{A_{v_i,v_j}(\operatorname{R}(G_{\left< \bm{i}, \bm{k} \right>}))}\Big]$ and thus the weight can only be non-zero, when every edge in the reversed graph $\operatorname{R}(G_{\left< \bm{i}, \bm{k} \right>})$ occurs at least twice. There can then only be at most $l$ many connections in $G_{\left< \bm{i}, \bm{k} \right>}$. As $G_{\left< \bm{i}, \bm{k} \right>}$ is connected, the only way for it to have $l+1$ many vertices is for $U(G_{\left< \bm{i}, \bm{k} \right>})$ to be a tree. It follows that $G_{\left< \bm{i}, \bm{k} \right>}$ must be a balanced tree.
\end{proof}

\subsection{Lemma (Number of balanced trees $T$ with given $B(T)$)}\label{CountingColoredTrees}
For any given set $B \subset V_{l+1} = \{v_1,...,v_{l+1}\}$ the number of balanced trees $T \in \mathcal{T}_{V_{l+1}}$ with $B(T) = B$ is zero, if $b := \#B \in \{0,l+1\}$, and is otherwise given by
\begin{align*}
	& \#\{T \in \mathcal{T}_{V_{l+1}} \mid B(T) = B\} = l! {l-1 \choose b-1} \ .
\end{align*}
\begin{proof}\
	\\
	Without loss of generality assume $B=\{v_1,...,v_{b}\}$, then by Lemma \ref{TreeColoring} the adjacency matrices $\{A(T)\}$ of the balanced trees in $\{T \in \mathcal{T}_{V_{l+1}} \mid B(T) = B\}$ are precisely all adjacency matrices $\{A(t)\}$ of undirected bipartite trees $t$ with vertices $\{v_1,...,v_b\}$ on the left and $\{v_{b+1},...,v_{l+1}\}$ on the right. For each such adjacency matrix $A(t)$ the B.E.S.T. Theorem tells us that that there are $2l \, \prod\limits_{i=1}^{l+1} (\deg_{t}(v_i)-1)!$ many $T \in \mathcal{T}_{V_{l+1}}$ with $A(T)=A(t)$. The coloring $B(T)$ will match $B$ in those cases, where we start on the left hand side, so we have
	\begin{align*}
		& \#\{T \in \mathcal{T}_{V_{l+1}} \mid A(T) = A(t) , \, B(T) = B\} = l \, \prod\limits_{i=1}^{l+1} (\deg_{t}(v_i)-1)! \ .
	\end{align*}
	By summing over all possible choices of degrees $(d_i)_{i \leq l+1}$ and then over all choices of $t$ with $d_i = \deg_{t}(v_i)$ we can use the fact that there are ${l-b \choose d_1-1,...,d_{b}-1} {b-1 \choose d_{b+1}-1,...,d_{l+1}-1}$ many bipartite trees $t$ with $d_i = \deg_{t}(v_i)$ to calculate
	\begin{align*}
		& \#\{T \in \mathcal{T}_{V_{l+1}} \mid B(T) = B\}\\
		& = \sum\limits_{\substack{d_1,...,d_b \geq 1 \\ d_1+...+d_b=l}} \sum\limits_{\substack{d_{b+1},...,d_{l+1} \geq 1 \\ d_{b+1}+...+d_{l+1}=l}} \sum\limits_{\substack{t \text{ bipartite tree} \\ \text{with } d_i = \deg_{t}(v_i)}} \#\{T \in \mathcal{T}_{V_{l+1}} \mid A(T) = A(t) , \, B(T) = B\}\\
		& = \sum\limits_{\substack{d_1,...,d_b \geq 1 \\ d_1+...+d_b=l}} \sum\limits_{\substack{d_{b+1},...,d_{l+1} \geq 1 \\ d_{b+1}+...+d_{l+1}=l}} \underbrace{{l-b \choose d_1-1,...,d_{b}-1} {b-1 \choose d_{b+1}-1,...,d_{l+1}-1} \, l \, \prod\limits_{v=1}^{l+1} (d_v-1)!}_{= l(l-b)! (b-1)!}\\
		& = l \, (l-b)! \, (b-1)! \, {l-1 \choose b-1} \, {l-1 \choose l-b} = l! {l-1 \choose b-1} \ .
	\end{align*}
\end{proof}

\subsection{Lemma (Coloring of one-directional ring-type graphs)}\label{RingColor_OneDirectional}
For even $l_0 \geq 4$ any one-directional ring-type graph $G \in \operatorname{1-d-Ring}_{l_0}$ will have alternating black and white vertices along its ring structure.\\
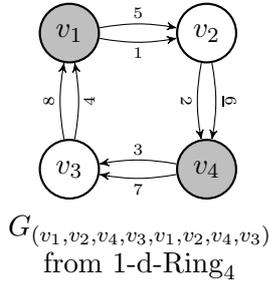
\begin{minipage}{0.45\textwidth}
	\centering
	\begin{tikzpicture}[node distance=\d and \d,>=stealth',auto, every place/.style={draw}]
		\node [black] (v1) {$v_1$};
		\node [white] (v2) [right=of v1] {$v_2$};
		\node [white] (v3) [below=of v1] {$v_3$};
		\node [black] (v4) [below=of v2] {$v_4$};
		
		\path[->]
		(v1) edge [bend right=10] node[el,below] {$1$} (v2)
		(v2) edge [bend right=10] node[el,below] {$2$} (v4)
		(v4) edge [bend right=10] node[el,above] {$3$} (v3)
		(v3) edge [bend right=10] node[el,below] {$4$} (v1)
		(v1) edge [bend right=-10] node[el,above] {$5$} (v2)
		(v2) edge [bend right=-10] node[el,above] {$\underline{6}$} (v4)
		(v4) edge [bend right=-10] node[el,below] {$7$} (v3)
		(v3) edge [bend right=-10] node[el,above] {$8$} (v1)
		;
		
	\end{tikzpicture}\\
	$G_{(v_1,v_2,v_4,v_3,v_1,v_2,v_4,v_3)}$\\
	from $\operatorname{1-d-Ring}_{4}$
\end{minipage}
\begin{minipage}{0.45\textwidth} 
	
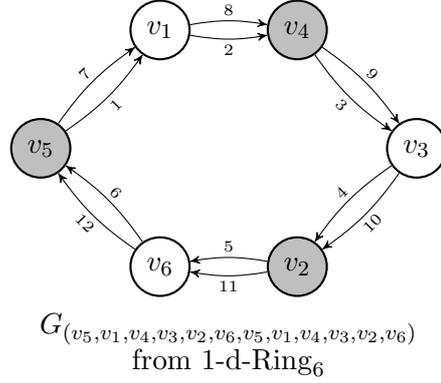
\centering
	\begin{tikzpicture}[node distance=\d and \d,>=stealth',auto, every place/.style={draw}]
		\node [white] (v1) {$v_1$};
		\node [black] (v2) [right=of v1] {$v_4$};
		\node [white] (v3) [below right=of v2] {$v_3$};
		\node [black] (v4) [below left=of v3] {$v_2$};
		\node [white] (v5) [left=of v4] {$v_6$};
		\node [black] (v6) [above left=of v5] {$v_5$};
		\path[->] (v1) edge [bend right=10] node[el,below] {$2$} (v2);
		\path[->] (v1) edge [bend right=-10] node[el,above] {$8$} (v2);
		
		\path[->] (v2) edge [bend right=10] node[el,below] {$3$} (v3);
		\path[->] (v2) edge [bend right=-10] node[el,above] {$9$} (v3);
		
		\path[->] (v3) edge [bend right=10] node[el,above] {$4$} (v4);
		\path[->] (v3) edge [bend right=-10] node[el,below] {$10$} (v4);
		
		\path[->] (v4) edge [bend right=10] node[el,above] {$5$} (v5);
		\path[->] (v4) edge [bend right=-10] node[el,below] {$11$} (v5);
		
		\path[->] (v5) edge [bend right=10] node[el,above] {$6$} (v6);
		\path[->] (v5) edge [bend left=10] node[el,below] {$12$} (v6);
		
		\path[->] (v6) edge [bend right=10] node[el,below] {$1$} (v1);
		\path[->] (v6) edge [bend right=-10] node[el,above] {$7$} (v1);
		
	\end{tikzpicture}\\
	$G_{(v_5,v_1,v_4,v_3,v_2,v_6,v_5,v_1,v_4,v_3,v_2,v_6)}$\\
	from $\operatorname{1-d-Ring}_{6}$
\end{minipage}
\begin{proof}\
	\\
	By construction two edges in the same undirected connection have the same parity and this parity must toggle along the ring structure.
\end{proof}

\subsection{Lemma (Coloring of two-directional ring-type graphs)}\label{RingColor_TwoDirectional}
For even $l_0 \geq 2$ any two-directional ring-type graph $G \in \operatorname{2-d-Ring}_{l_0}$ will have alternating black and white vertices along its ring structure. If $l_0$ is odd, then there exists one exception where two neighboring vertices are both black.\\

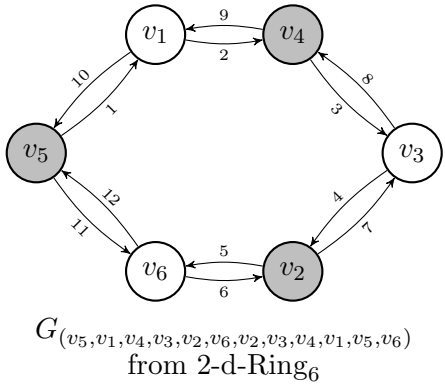
\begin{minipage}{0.45\textwidth} 
	\centering
	\begin{tikzpicture}[node distance=\d and \d,>=stealth',auto, every place/.style={draw}]
		\node [white] (v1) {$v_1$};
		\node [black] (v2) [right=of v1] {$v_4$};
		\node [white] (v3) [below right=of v2] {$v_3$};
		\node [black] (v4) [below left=of v3] {$v_2$};
		\node [white] (v5) [left=of v4] {$v_6$};
		\node [black] (v6) [above left=of v5] {$v_5$};
		
		\path[->]
		(v1) edge [bend right=10] node[el,below] {$2$} (v2)
		(v2) edge [bend right=10] node[el,above] {$9$} (v1)
		
		(v2) edge [bend right=10] node[el,below] {$3$} (v3)
		(v3) edge [bend right=10] node[el,above] {$8$} (v2)
		
		(v3) edge [bend right=10] node[el,above] {$4$} (v4)
		(v4) edge [bend right=10] node[el,below] {$7$} (v3)
		
		(v4) edge [bend right=10] node[el,above] {$5$} (v5)
		(v5) edge [bend right=10] node[el,below] {$6$} (v4)
		
		(v5) edge [bend right=10] node[el,above] {$12$} (v6)
		(v6) edge [bend right=10] node[el,below] {$11$} (v5)
		
		(v6) edge [bend right=10] node[el,below] {$1$} (v1)
		(v1) edge [bend right=10] node[el,above] {$10$} (v6)
		
		;
		
	\end{tikzpicture}\\
	$G_{(v_5,v_1,v_4,v_3,v_2,v_6,v_2,v_3,v_4,v_1,v_5,v_6)}$\\
	from $\operatorname{2-d-Ring}_{6}$
\end{minipage}
\begin{minipage}{0.45\textwidth} 
	
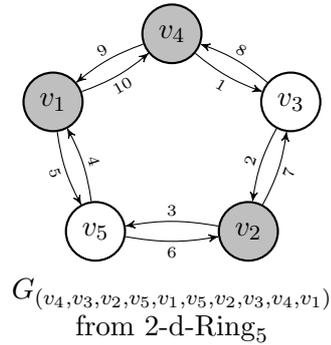
\centering
	\begin{tikzpicture}[node distance=\d and \d,>=stealth',auto, every place/.style={draw}]
		\node [black] (v1) {$v_4$};
		\node [white] (v2) [position=-30:{\d} from v1] {$v_3$};
		\node [black] (v3) [position=-108:{\d} from v2] {$v_2$};
		\node [black] (v5) [position=-150:{\d} from v1] {$v_1$};
		\node [white] (v4) [position=-72:{\d} from v5] {$v_5$};
		
		\path[->]
		(v1) edge [bend right=10] node[el,below] {$1$} (v2)
		(v2) edge [bend right=10] node[el,above] {$8$} (v1)
		
		(v2) edge [bend right=10] node[el,above] {$2$} (v3)
		(v3) edge [bend right=10] node[el,below] {$7$} (v2)
		
		(v3) edge [bend right=10] node[el,above] {$3$} (v4)
		(v4) edge [bend right=10] node[el,below] {$6$} (v3)
		
		(v4) edge [bend right=10] node[el,above] {$4$} (v5)
		(v5) edge [bend right=10] node[el,below] {$5$} (v4)
		
		(v5) edge [bend right=10] node[el,below] {$10$} (v1)
		(v1) edge [bend right=10] node[el,above] {$9$} (v5)
		;
		
	\end{tikzpicture}\\
	$G_{(v_4,v_3,v_2,v_5,v_1,v_5,v_2,v_3,v_4,v_1)}$\\
	from $\operatorname{2-d-Ring}_{5}$
\end{minipage}
\begin{proof}\
	\\
	The above pictures explain the coloring properties best. The route starts at some vertex $v_i$ and passes through the ring structure until at one point - lets say at vertex $v_j$ - it must change direction. After this first change of direction a full circuit is completed until the route arrives back at $v_j$, where a second change of direction must occur and the route walks back to the starting vertex $v_i$. It is possible for the starting vertex and the vertex, where the direction is changed, to be the same.\\
	Since every second step in the route colors its vertex black, one easily checks that for an even number of vertices the change of direction at vertex $v_j$ does not change the natural alternation of black and white vertices. Meanwhile, for an odd number of vertices we get two neighboring black vertices around the vertex $v_j$.
\end{proof}

\section{Counting sprouts}\label{SectionCountingSprouting}

\subsection{Lemma (Number of sprouts only depends on $l_0$)}\label{SproutingBijection}
For any $r_0,\overline{r_0},l_0 \in \N$ let $G_0 \in \mathcal{C}_{V_{r_0},2l_0}$ and $\overline{G_0} \in \mathcal{C}_{\overline{V}_{\overline{r_0}},2l_0}$ be circuit graphs without balanced leaves. Further let $V'_{l'}$ be a vertex set disjoint to both $V_{r_0}$ and $\overline{V}_{\overline{r_0}}$. We call $V'_{l'}$ the set of sprouted vertices, while $V_{r_0}$ and $\overline{V}_{\overline{r_0}}$ are two possible choices of seed vertices.\\
\\
If $l_0 = 2$ and $G_0$ is of the form $G_{v_i,v_j}$ for $\{v_i,v_j\} = \{v_1,v_2\}=V_2$ we assume that $v_1,v_2$ are larger than all vertices in $V'_{l'}$. The same goes for $\overline{G}_0$. We do this in order to not run into the problems addressed in Remark \ref{SeedsOfTrees}.\\
\\
Under these conditions for any $B \subset V'_{l'}$ the number of circuit graphs, which are sprouts of $G_0$, such that $B$ is the set of black sprouted vertices is, precisely the number of circuit graphs, which are sprouts of $\overline{G_0}$, such that $B$ is the set of black sprouted vertices. More precisely we have the equality
\begin{align*}
	& \#\{G \in \mathcal{C}_{V_{r_0} \sqcup V'_{l'},2l_0+2l'} \mid \operatorname{S}(G) = G_0 , \, B(G) \cap V'_{l'} = B \}\\
	& = \#\{G \in \mathcal{C}_{\overline{V}_{\overline{r_0}} \sqcup V'_{l'},2l_0+2l'} \mid \operatorname{S}(G) = \overline{G_0} , \, B(G) \cap V'_{l'} = B \} \ .
\end{align*}
We prove this by constructing a bijection $\Phi$ between the sets. Here some examples of how $\Phi$ maps sprouts:\\
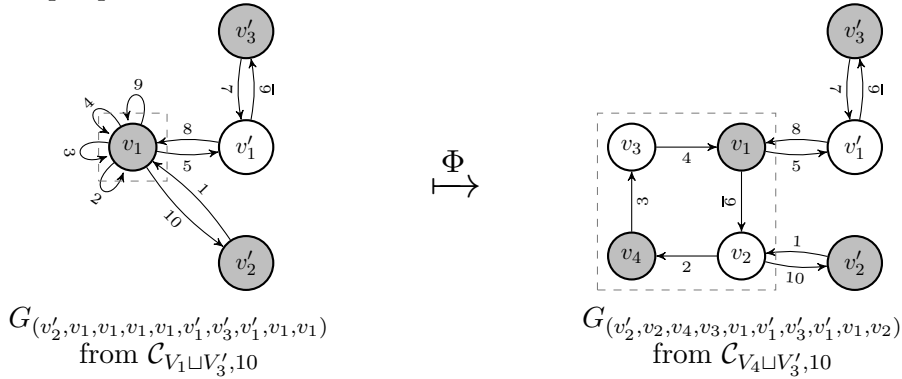
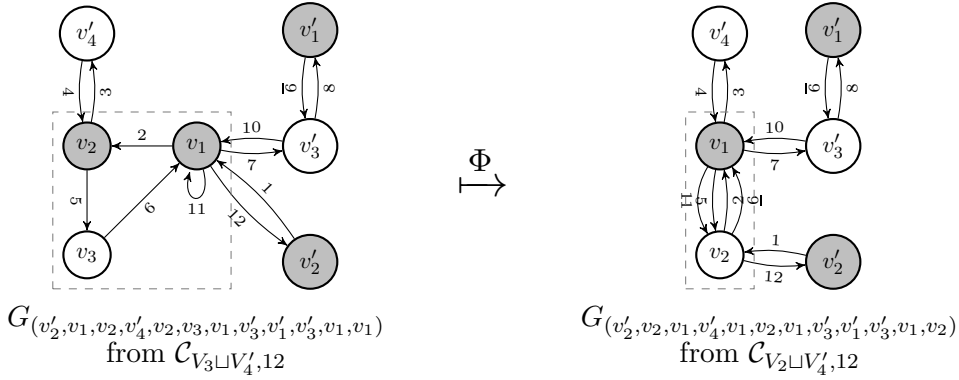
\begin{minipage}{0.45\textwidth} 
	\centering
	\renewcommand\s{0.8}
	\begin{tikzpicture}[node distance=\d and \d,>=stealth',auto, every place/.style={draw}]
		\node [black] (v1) {$v_1$};
		\node [white] (v2) [right=of v1] {$v'_1$};
		\node [black] (v3) [below=of v2] {$v'_2$};
		\node [black] (v4) [above=of v2] {$v'_3$};
		
		\path[->]
		(v1) edge [in=250,out=220,loop] node[el,below] {$2$} (v1)
		(v1) edge [in=200,out=170,loop] node[el,below] {$3$} (v1)
		(v1) edge [in=150,out=120,loop] node[el,above] {$4$} (v1)
		(v1) edge [in=100,out=70,loop] node[el,above] {$9$} (v1)
		
		(v1) edge [bend right=10] node[el,below] {$10$} (v3)
		(v3) edge [bend right=10] node[el,above] {$1$} (v1)
		
		(v1) edge [bend right=10] node[el,below] {$5$} (v2)
		(v2) edge [bend right=10] node[el,above] {$8$} (v1)
		
		(v2) edge [bend right=10] node[el,below] {$\underline{6}$} (v4)
		(v4) edge [bend right=10] node[el,below] {$7$} (v2)
		;
		
		\node[fit=(v1), draw, dashed,gray] {};
	\end{tikzpicture}\\
	$G_{(v'_2,v_1,v_1,v_1,v_1,v'_1,v'_3,v'_1,v_1,v_1)}$\\
	from $\mathcal{C}_{V_1 \sqcup V'_3,10}$
\end{minipage}
\scalebox{1.5}{$\xmapsto{\Phi}$}
\begin{minipage}{0.45\textwidth} 
	\centering
	\renewcommand\s{0.8}
	\begin{tikzpicture}[node distance=\d and \d,>=stealth',auto, every place/.style={draw}]
		\node [black] (v1) {$v_1$};
		\node [white] (v5) [below=of v1] {$v_2$};
		\node [white] (v2) [right=of v1] {$v'_1$};
		\node [black] (v3) [below=of v2] {$v'_2$};
		\node [black] (v4) [above=of v2] {$v'_3$};
		\node [white] (v6) [left=of v1] {$v_3$};
		\node [black] (v7) [left=of v5] {$v_4$};
		
		\path[->]
		(v7) edge [bend right=0] node[el,below] {$3$} (v6)
		(v5) edge [bend right=0] node[el,below] {$2$} (v7)
		(v1) edge [bend right=0] node[el,below] {$\underline{9}$} (v5)
		(v6) edge [bend right=0] node[el,below] {$4$} (v1)
		
		(v5) edge [bend right=10] node[el,below] {$10$} (v3)
		(v3) edge [bend right=10] node[el,above] {$1$} (v5)
		
		(v1) edge [bend right=10] node[el,below] {$5$} (v2)
		(v2) edge [bend right=10] node[el,above] {$8$} (v1)
		
		(v2) edge [bend right=10] node[el,below] {$\underline{6}$} (v4)
		(v4) edge [bend right=10] node[el,below] {$7$} (v2)
		;
		
		\node[fit=(v6)(v5), draw, dashed,gray] {};
	\end{tikzpicture}\\
	$G_{(v'_2,v_2,v_4,v_3,v_1,v'_1,v'_3,v'_1,v_1,v_2)}$\\
	from $\mathcal{C}_{V_4 \sqcup V'_3,10}$
\end{minipage}
\\
\vspace{0.2cm}
\\
\begin{minipage}{0.45\textwidth} 
	\centering
	\renewcommand\s{0.8}
	\begin{tikzpicture}[node distance=\d and \d,>=stealth',auto, every place/.style={draw}]
		\node [black] (v1) {$v_1$};
		\node [white] (v2) [right=of v1] {$v'_3$};
		\node [black] (v3) [below=of v2] {$v'_2$};
		\node [black] (v4) [above=of v2] {$v'_1$};
		\node [black] (v5) [left=of v1] {$v_2$};
		\node [white] (v6) [below=of v5] {$v_3$};
		\node [white] (v7) [above=of v5] {$v'_4$};
		
		\path[->]
		(v1) edge node[el,above] {$2$} (v5)
		(v5) edge node[el,below] {$5$} (v6)
		(v6) edge node[el,below] {$6$} (v1)
		(v1) edge [in=-105,out=-75,loop] node[el,below] {$11$} (v1)

		(v5) edge [bend right=10] node[el,below] {$3$} (v7)
		(v7) edge [bend right=10] node[el,below] {$4$} (v5)
		
		(v1) edge [bend right=10] node[el,below] {$12$} (v3)
		(v3) edge [bend right=10] node[el,above] {$1$} (v1)
		
		(v1) edge [bend right=10] node[el,below] {$7$} (v2)
		(v2) edge [bend right=10] node[el,above] {$10$} (v1)
		
		(v2) edge [bend right=10] node[el,below] {$8$} (v4)
		(v4) edge [bend right=10] node[el,below] {$\underline{9}$} (v2)
		;
		
		\node[fit=(v1)(v5)(v6), draw, dashed,gray] {};
	\end{tikzpicture}\\
	$G_{(v'_2,v_1,v_2,v'_4,v_2,v_3,v_1,v'_3,v'_1,v'_3,v_1,v_1)}$\\
	from $\mathcal{C}_{V_3 \sqcup V_{4}',12}$
\end{minipage}
\scalebox{1.5}{$\xmapsto{\Phi}$}
\begin{minipage}{0.45\textwidth} 
	\centering
	\renewcommand\s{0.8}
	\begin{tikzpicture}[node distance=\d and \d,>=stealth',auto, every place/.style={draw}]
		\node [black] (v1) {$v_1$};
		\node [white] (v5) [below=of v1] {$v_2$};
		\node [white] (v2) [right=of v1] {$v'_3$};
		\node [black] (v3) [below=of v2] {$v'_2$};
		\node [black] (v4) [above=of v2] {$v'_1$};
		\node [white] (v6) [above=of v1] {$v'_4$};
		
		\path[->]
		(v1) edge [bend right=10] node[el,below] {$5$} (v5)
		(v5) edge [bend right=10] node[el,below] {$2$} (v1)
		(v1) edge [bend right=30] node[el,below] {$11$} (v5)
		(v5) edge [bend right=30] node[el,below] {$\underline{6}$} (v1)
		
		(v5) edge [bend right=10] node[el,below] {$12$} (v3)
		(v3) edge [bend right=10] node[el,above] {$1$} (v5)
		
		(v1) edge [bend right=10] node[el,below] {$3$} (v6)
		(v6) edge [bend right=10] node[el,below] {$4$} (v1)
		
		(v1) edge [bend right=10] node[el,below] {$7$} (v2)
		(v2) edge [bend right=10] node[el,above] {$10$} (v1)
		
		(v2) edge [bend right=10] node[el,below] {$8$} (v4)
		(v4) edge [bend right=10] node[el,below] {$\underline{9}$} (v2)
		;
		
		\node[fit=(v1)(v5), draw, dashed,gray] {};
	\end{tikzpicture}\\
	$G_{(v'_2,v_2,v_1,v'_4,v_1,v_2,v_1,v'_3,v'_1,v'_3,v_1,v_2)}$\\
	from $\mathcal{C}_{V_2 \sqcup V_{4}',12}$
\end{minipage}
\begin{proof}\
	\\
	Note that the number of occurrences of a seed-vertex $v_\bullet$ in the route $\bm{i}$ of the sprout can be larger than the number of occurrences in the route $\bm{i}_0$ of the seed graph. For example the route of the lower left graph of our examples is $\bm{i} = (v'_2,v_1,v_2,v'_4,v_2,v_3,v_1,v'_3,v'_1,v'_3,v_1,v_1)$, while its seed route is $\bm{i}_0 = (v_2,v_3,v_1,v_1)$. Superfluous occurrences of $v_\bullet$ are added to the route $\bm{i}$ by tree structures connecting/returning to the seed graph. Consequently, if we first remove all seed-vertices $v_\bullet$, which are directly behind a sprouting vertex $v'_\bullet$ and then also remove all sprouting vertices, we are left with the seed route $\bm{i}_0$.\\
	\\
	The idea behind the construction of $\Phi$ is to use the above observation to identify the $l_0$ many critical positions in the route $\bm{i}$, which correspond to the seed route $\bm{i}_0$ and replace them with the respective entries of the other seed route $\ol{\bm{i}_0}$. The superfluous seed-vertices $v_\bullet$, which are not critical, must be changed to match the last critical $v_\bullet$ in the route. The critical positions clearly stay critical in the new route and the entire procedure is then easily seen to be reversible.
\end{proof}

\subsection{Proposition (Counting sprouts)}\label{CountingSprouting}
For any $l_0 \in \N$ let $G_0 \in \mathcal{C}_{V_2,2l_0}$ be the circuit graph with route
\begin{align*}
	& \bm{i}_0 = (\underbrace{v_1,v_2,v_1,v_2,...,v_1,v_2}_{\text{length: }2l_0}) \ .
\end{align*}
For any finite sets $B',W' \subset \N$ such that $V_2=\{v_1,v_2\}$,$B'$ and $W'$ are disjoint define $b':= \#B'$, $w' := \#W'$, $l'=b'+w'$ and $V'_{l'} = B' \sqcup W'$. Further assume that $v_1,v_2$ are the two largest vertices in $V_2 \sqcup V'_{l'}$ (see Remark \ref{SeedsOfTrees}). We then have
\begin{align*}
	& \#\underbrace{\{G \in \mathcal{C}_{V_2\sqcup V'_{l'},2l_0+2l'} \mid \operatorname{S}(G) = G_0, \, B(G) \cap V'_{l'} = B' \}}_{=: \operatorname{Spr}_{B',W'}(G_0)} = \frac{(l_0+b'+w')!^2}{(l_0+b')!(l_0+w')!} \ .
\end{align*}
By Lemma \ref{SproutingBijection} this is also true for arbitrary other $G_0 \in \mathcal{C}_{V_{r_0},2l_0}$ without balanced leaves. Since for other $G_0 \in \mathcal{C}_{V_{r_0},E_{2l_0}}$ the order of $V_{r_0} \sqcup V'_{l'}$ has no effect on the seed graph of a $G \in \mathcal{C}_{V_{r_0} \sqcup B' \sqcup W', 2l_0+2b'+2w'}$, we in this case also don't need to make any assumptions about elements of $V_{r_0}$ being larger than those of $V'_{l'}$.
\begin{proof}\
	\\
	Given the coloring of the sprouting vertices, any sprout $G$ of $G_0$ will have a bipartite structure in accordance with the coloring of the vertices and we can uniquely identify every possible adjacency matrix $A(G)$ of a sprout with a bipartite tree, where the edge $(v_1,v_2)$ is prescribed. Below we see a possible sprout of $G$ (without edge labels) and the corresponding bipartite tree $t$ with $b'+1$ vertices on the left and $w'+1$ vertices on the right.\\
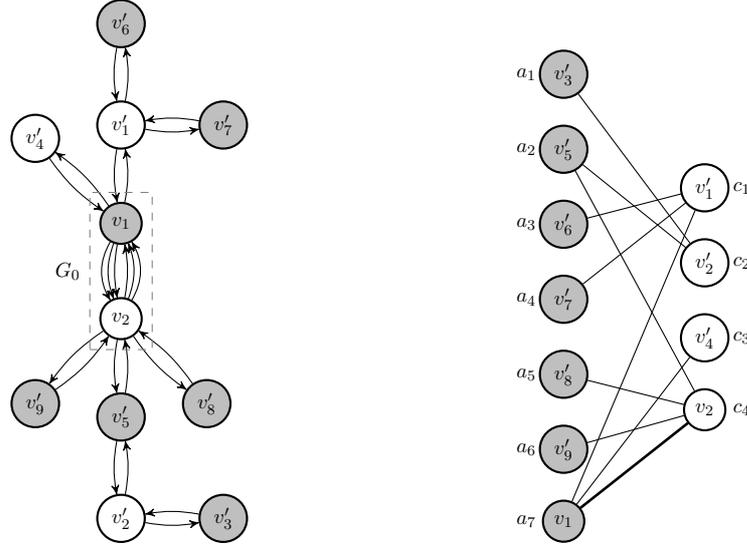
	\begin{minipage}{0.95\textwidth}
		\centering
		\renewcommand\s{0.7}
		\begin{tikzpicture}[node distance=\d and \d,>=stealth',auto, every place/.style={draw}]
			\node [black] (v1) {$v_1$};
			\node [white] (v2) [below=of v1] {$v_2$};
			\path[->] (v1) edge [bend right=10] node {} (v2);
			\path[->] (v1) edge [bend right=20] node {} (v2);
			\path[->] (v1) edge [bend right=30] node {} (v2);
			\path[->] (v2) edge [bend right=10] node {} (v1);
			\path[->] (v2) edge [bend right=20] node {} (v1);
			\path[->] (v2) edge [bend right=30] node {} (v1);
			
			\node [white] (v'1) [above=of v1] {$v'_1$};
			\node [white] (v'4) [above left=of v1] {$v'_4$};
			\node [black] (v'7) [right=of v'1] {$v'_7$};
			\node [black] (v'6) [above=of v'1] {$v'_6$};
			\path[->] (v1) edge [bend right=10] node {} (v'1);
			\path[->] (v'1) edge [bend right=10] node {} (v1);
			\path[->] (v1) edge [bend right=10] node {} (v'4);
			\path[->] (v'4) edge [bend right=10] node {} (v1);
			\path[->] (v'1) edge [bend right=10] node {} (v'7);
			\path[->] (v'7) edge [bend right=10] node {} (v'1);
			\path[->] (v'1) edge [bend right=10] node {} (v'6);
			\path[->] (v'6) edge [bend right=10] node {} (v'1);

			\node [black] (v'5) [below=of v2] {$v'_5$};
			\node [black] (v'8) [below right=of v2] {$v'_8$};
			\node [black] (v'9) [below left=of v2] {$v'_9$};
			\node [white] (v'2) [below=of v'5] {$v'_2$};
			\node [black] (v'3) [right=of v'2] {$v'_3$};
			\path[->] (v2) edge [bend right=10] node {} (v'5);
			\path[->] (v'5) edge [bend right=10] node {} (v2);
			\path[->] (v2) edge [bend right=10] node {} (v'9);
			\path[->] (v'9) edge [bend right=10] node {} (v2);
			\path[->] (v2) edge [bend right=10] node {} (v'8);
			\path[->] (v'8) edge [bend right=10] node {} (v2);
			\path[->] (v'5) edge [bend right=10] node {} (v'2);
			\path[->] (v'2) edge [bend right=10] node {} (v'5);
			\path[->] (v'2) edge [bend right=10] node {} (v'3);
			\path[->] (v'3) edge [bend right=10] node {} (v'2);

			
			
			\node[fit=(v1)(v2), draw, dashed,gray] {};
			\node[noborder] at ($(v1)+(-1.0*\d,-0.9*\d)$) {$G_0$};

			
			
			
		\end{tikzpicture}
		\ \hspace{3cm} 
		\begin{tikzpicture}[node distance=\d and \d,>=stealth',auto, every place/.style={draw}]
			\node [black] (v'3) {$v'_3$};
			\node [black] (v'5) [below =0.5*\d of v'3] {$v'_5$};
			\node [black] (v'6) [below =0.5*\d of v'5] {$v'_6$};
			\node [black] (v'7) [below =0.5*\d of v'6] {$v'_7$};
			\node [black] (v'8) [below =0.5*\d of v'7] {$v'_8$};
			\node [black] (v'9) [below =0.5*\d of v'8] {$v'_9$};
			\node [black] (v1) [below =0.5*\d of v'9] {$v_1$};
			
			\node [white] (v'1) [below right =1.5*\d and 2*\d of v'3] {$v'_1$};
			\node [white] (v'2) [below =0.5*\d of v'1] {$v'_2$};
			\node [white] (v'4) [below =0.5*\d of v'2] {$v'_4$};
			\node [white] (v2) [below =0.5*\d of v'4] {$v_2$};
			
			\draw[line width=0.05*\d] (v1) -- (v2);
			\draw (v1) -- (v'1);
			\draw (v1) -- (v'4);
			\draw (v'1) -- (v'6);
			\draw (v'1) -- (v'7);
			\draw (v2) -- (v'9);
			\draw (v2) -- (v'5);
			\draw (v2) -- (v'8);
			\draw (v'5) -- (v'2);
			\draw (v'3) -- (v'2);
			
			\node[noborder] at ($(v1)+(-0.7*\d,0)$) {$a_{7}$};
			\node[noborder] at ($(v'3)+(-0.7*\d,0)$) {$a_1$};
			\node[noborder] at ($(v'5)+(-0.7*\d,0)$) {$a_2$};
			\node[noborder] at ($(v'6)+(-0.7*\d,0)$) {$a_3$};
			\node[noborder] at ($(v'7)+(-0.7*\d,0)$) {$a_4$};
			\node[noborder] at ($(v'8)+(-0.7*\d,0)$) {$a_5$};
			\node[noborder] at ($(v'9)+(-0.7*\d,0)$) {$a_6$};
			
			\node[noborder] at ($(v2)+(0.7*\d,0)$) {$c_4$};
			\node[noborder] at ($(v'1)+(0.7*\d,0)$) {$c_1$};
			\node[noborder] at ($(v'2)+(0.7*\d,0)$) {$c_2$};
			\node[noborder] at ($(v'4)+(0.7*\d,0)$) {$c_3$};
		\end{tikzpicture}
	\end{minipage}\\
	\\
	We can now use the same idea as in Lemma \ref{CountingColoredTrees} by first summing over every possible bipartite tree $t$ with prescribed edge $(v_1,v_2)$ and then counting the number of sprouts $G$ whose adjacency matrix $A(G)$ corresponds to said bipartite tree $t$.\\
	As the connection between $v_1$ and $v_2$ in the left hand graph above is $l_0$-fold, it has $l_0$ many spanning trees and by the B.E.S.T. theorem there are $l_0 \, \prod\limits_{v \in V_2 \sqcup B' \sqcup W'} (\deg_A(v)-1)!$ many Euler-tours through such a graph. Here $A$ is the adjacency matrix of the fixed graph. Since Euler-tours are counted modulo starting edge, we can multiply this by the number of edges originating from ab black vertex $(l_0+l')$ to get the number of sprouts $G$, whose adjacency matrix $A(G)$ corresponds to $A$. We have shown
	\begin{align}
		& \#\{G \in \operatorname{Spr}_{B',W'}(G_0) \mid A(G)=A\} = (l_0+l') l_0 \, \prod\limits_{v \in V_2 \sqcup B' \sqcup W'} (\deg_A(v)-1)! \ .
	\end{align}
	Let $t$ denote the bipartite tree corresponding to $A$, then the connection $(v_1,v_2)$ is no longer $l_0$-fold in $t$ and we have $\deg_{t}(v_{1/2}) = \deg_{A}(v_{1/2})-l_0+1$, while the degrees of the other vertices remain unchanged. We so far have
	\begin{align*}
		& \#\operatorname{Spr}_{B',W'}(G_0)\\
		& = \sum\limits_{\substack{t \text{ bip. tree} \\ (v_1,v_2) \text{ edge in }t}} (l_0+l') l_0 \, (\deg_{t}(v_{1})+l_0-2)! (\deg_{t}(v_{2})+l_0-2)! \, \prod\limits_{v \in B' \sqcup W'} (\deg_t(v)-1)!\\
		& = (l_0+l') l_0 \sum\limits_{\substack{(d_{v})_{v \in V_2 \sqcup B' \sqcup W'} \subset \N \\ d_{v_1}+\sum\limits_{v \in B'} d_v = l'+1 \\ d_{v_2}+\sum\limits_{v \in W'} d_v = l'+1}} (d_{v_{1}}+l_0-2) (d_{v_{2}}+l_0-2) \, \prod\limits_{v \in B' \sqcup W'} (d_v-1)! \sum\limits_{\substack{t \text{ bip. tree} \\ (v_1,v_2) \text{ edge in }t \\ \forall v : \, \deg_t(v) = d_v}} 1 \ .
	\end{align*}
	where by bipartite tree we always mean that the vertices $\{v_1\}\sqcup B'$ are on the left hand side and $\{v_2\}\sqcup W'$ are on the right hand side. In the last step we have changed the order of summation to first go over all possible degrees of vertices on the left and right hand side of the bipartite tree. Rename the degrees of the left hand vertices $(d_v)_{v \in \{v_1\} \sqcup B'}$ into $d_1,...,d_{b'+1}$ with $d_1 = d_{v_1}$ and the degrees of the right hand vertices $(d_v)_{v \in \{v_2\} \sqcup W'}$ into $e_1,...,e_{w'+1}$ with $e_1 = d_{v_2}$, then in Lemma \ref{BipartiteSpanningTrees} we inductively show the formula
	\begin{align*}
		& \sum\limits_{\substack{t \text{ bip. tree} \\ (v_1,v_2) \text{ edge in }t \\ \forall v : \, \deg_t(v) = d_v}} 1 = \bigg( 1- \mathbbm{1}_{b'>0,w'>0} \frac{(b'-e_{1}+1)(w'-d_{1}+1)}{b'w'} \bigg)\\
		& \hspace{4cm} \times {b' \choose e_1-1,...,e_{w'+1}-1} \, {w' \choose d_1-1,...,d_{b'+1}-1} \ .
	\end{align*}
	The two above equalities with some calculations imply
	\begin{align*}
		& \#\operatorname{Spr}_{B',W'}(G_0) = \frac{(l_0+b'+w')!^2}{(l_0+b')!(l_0+w')!} \ . \qedhere
	\end{align*}
\end{proof}

\subsection{Corollary (Counting sprouts of one-directional rings)}\label{SproutingCorollary1}
For even $l_0 \geq 4$ and any $b',w' \geq 0$ define $l := l_0 +b'+w'$ and let $B$ be a subset of $V_l$ such that $b := \#B = b'+\frac{l_0}{2}$. The number of $G \in \mathcal{C}_{V_{l},2l}$ with $B(G) = B$ and the property that $\operatorname{S}(G) \in \operatorname{1-d-Ring}_{l_0}$ is given by
\begin{align*}
	& b! \, w! {l \choose b'}{l \choose w'} \ .
\end{align*}
Here $b := b'+\frac{l_0}{2}$ denotes the total number of black vertices and $w := w'+ \frac{l_0}{2} = l-b$ denotes the total number of white vertices.
\begin{proof}\
	\\
	By Lemma \ref{RingColor_OneDirectional} there will be $\frac{l_0}{2}$ many black and white vertices respectively in $\operatorname{S}(G)$, meaning there are ${ b'+\frac{l_0}{2} \choose \frac{l_0}{2}} {w'+\frac{l_0}{2} \choose \frac{l_0}{2}}$ many choices for which vertices make up $\operatorname{S}(G)$. For the route of $\operatorname{S}(G)$ there are now $\frac{l_0}{2}!$ choices for the order of the black vertices and independently $\frac{l_0}{2}!$ choices for the ordering of the white vertices of $\operatorname{S}(G)$. Since the sets $B' = B \setminus \operatorname{S}(G)$ and $W' = V_l \setminus (\operatorname{S}(G) \cup B')$ are uniquely defined and have cardinalities $b'$ and $w'$ respectively, we may now apply Proposition \ref{CountingSprouting} for fixed $\operatorname{S}(G)$ to see that there are
	\begin{align*}
		& \frac{(l_0+b'+w')!^2}{(l_0+b')!(l_0+w')!}
	\end{align*}
	many $G$ to each chosen $\operatorname{S}(G)$. In total we have counted
	\begin{align*}
		& {b'+\frac{l_0}{2} \choose \frac{l_0}{2}} {w'+\frac{l_0}{2} \choose \frac{l_0}{2}} \, \frac{l_0}{2}! \, \frac{l_0}{2}! \, \frac{(l_0+b'+w')!^2}{(l_0+b')(l_0+w')} = b! \, w! {l \choose b'}{l \choose w'}
	\end{align*}
	many $G$ with the properties described above.
\end{proof}

\subsection{Corollary (Counting sprouts of two-directional rings)}\label{SproutingCorollary2}
For $l_0 \neq 2$ and any $b',w' \geq 0$ define $l := l_0 +b'+w'$ and let $B$ be a subset of $V_l$ such that $b := \#B = b'+\lceil\frac{l_0}{2}\rceil$. The number of $G \in \mathcal{C}_{V_{l},2l}$ with $B(G) = B$ and the property that $\operatorname{S}(G) \in \operatorname{2-d-Ring}_{l_0}$ is given by $l_0 \, b! \, w! {l \choose b'}{l \choose w'}$.\\
Here $b := b'+\lceil\frac{l_0}{2}\rceil$ denotes the total number of black vertices and $w := w'+ \lfloor\frac{l_0}{2}\rfloor = l-b$ denotes the total number of white vertices.\\
\\
For $l_0 = 2$ there are only
\begin{align*}
	& b! \, w! \, {l \choose b'}{l \choose w'} = b! \, w! \, {b'+w'+2 \choose b'}{b'+w'+2 \choose w'}
\end{align*}
such $G$, which is only half as many as expected by the formula for $l_0>2$.
\begin{proof}\
	\\
	By Lemma \ref{RingColor_TwoDirectional} there will be $\lceil \frac{l_0}{2} \rceil$ black and $\lfloor {\frac{l_0}{2}} \rfloor$ white vertices in $\operatorname{S}(G)$, meaning there are ${b'+\lceil\frac{l_0}{2}\rceil \choose \lceil\frac{l_0}{2}\rceil} {w'+\lfloor\frac{l_0}{2}\rfloor \choose \lfloor\frac{l_0}{2}\rfloor}$ many choices for which vertices make up $\operatorname{S}(G)$. The form of the route of $\operatorname{S}(G)$ was discussed in the proof of Lemma \ref{RingColor_TwoDirectional}. We have $\lceil\frac{l_0}{2}\rceil!$ choices for the ordering of the black vertices in the route, $\lfloor\frac{l_0}{2}\rfloor!$ choices for the ordering of the white vertices of the route and $l_0$ choices where to initially change direction. Since the sets $B' = B \setminus \operatorname{S}(G)$ and $W' = V_l \setminus (\operatorname{S}(G) \cup B')$ are uniquely defined and have cardinalities $b'$ and $w'$ respectively, we may now apply Proposition \ref{CountingSprouting} to see that there are
	\begin{align*}
		& \frac{(l_0+b'+w')!^2}{(l_0+b')!(l_0+w')!}
	\end{align*}
	many $G$ to each chosen $\operatorname{S}(G)$. In total we have counted
	\begin{align*}
		& {b'+\lceil\frac{l_0}{2}\rceil \choose \lceil\frac{l_0}{2}\rceil} {w'+\lfloor\frac{l_0}{2}\rfloor \choose \lfloor\frac{l_0}{2}\rfloor} \, l_0 \, \Big\lceil\frac{l_0}{2}\Big\rceil! \, \Big\lfloor\frac{l_0}{2}\Big\rfloor! \, \frac{(l_0+b'+w')!^2}{(l_0+b')(l_0+w')} = l_0 \, b! \, w! {l \choose b'}{l \choose w'}
	\end{align*}
	many $G$ with the properties described above.\\
	\\
	For $l_0 = 2$ the only change in our argument is that there are not $l_0 \, \lceil\frac{l_0}{2}\rceil! \, \lfloor\frac{l_0}{2}\rfloor! = 2$ many choices for $\operatorname{S}(G)$ but only $1$ choice.
\end{proof}

\section{Proving Theorem \ref{Thm_MeanExpansion}}\label{SectionMainRes1}

\subsection{Proposition (Weight of $G \in \mathcal{C}_{V_l,E_{2l}}$)}\label{WeightOfFusedTrees}
For any $l \in \N$ suppose $V_l \subset \N$. For $\bm{i} \in V_l^l, \bm{k} \in V_l^l$ with $\bm{\{i\}} \cup \bm{\{k\}} = V_l$ we have $G_{\left< \bm{i}, \bm{k} \right>} \in \mathcal{C}_{V_l,2l}$ and
\begin{align*}
	& \cW(G_{\left< \bm{i}, \bm{k} \right>}) = \begin{cases}
		\E[X_{1 \, 1}^4] & \text{, if $\operatorname{S}(G_{\left< \bm{i}, \bm{k} \right>}) \in \operatorname{2-d-Ring}_2$}\\
		1 & \text{, if $\operatorname{S}(G_{\left< \bm{i}, \bm{k} \right>}) \in \operatorname{2-d-Ring}_{l_0}$ for some $l_0 \neq 2$}\\
		1 & \text{, if $\operatorname{S}(G_{\left< \bm{i}, \bm{k} \right>}) \in \operatorname{1-d-Ring}_{l_0}$ for even $l_0 \geq 4$}\\
		0 & \text{ else}
	\end{cases} \ .
\end{align*}
($\operatorname{1/2-d-Ring}_{l_0}$ was introduced in Definition \ref{DefRing})
\vspace{-0.3cm}\\
\renewcommand\s{0.95}
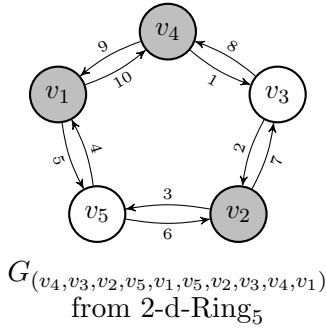
\begin{minipage}{0.45\textwidth} 
	\centering
	\begin{tikzpicture}[node distance=\d and \d,>=stealth',auto, every place/.style={draw}]
		\node [black] (v1) {$v_4$};
		\node [white] (v2) [position=-30:{\d} from v1] {$v_3$};
		\node [black] (v3) [position=-108:{\d} from v2] {$v_2$};
		\node [black] (v5) [position=-150:{\d} from v1] {$v_1$};
		\node [white] (v4) [position=-72:{\d} from v5] {$v_5$};
		
		\path[->]
		(v1) edge [bend right=10] node[el,below] {$1$} (v2)
		(v2) edge [bend right=10] node[el,above] {$8$} (v1)
		
		(v2) edge [bend right=10] node[el,above] {$2$} (v3)
		(v3) edge [bend right=10] node[el,below] {$7$} (v2)
		
		(v3) edge [bend right=10] node[el,above] {$3$} (v4)
		(v4) edge [bend right=10] node[el,below] {$6$} (v3)
		
		(v4) edge [bend right=10] node[el,above] {$4$} (v5)
		(v5) edge [bend right=10] node[el,below] {$5$} (v4)
		
		(v5) edge [bend right=10] node[el,below] {$10$} (v1)
		(v1) edge [bend right=10] node[el,above] {$9$} (v5)
		;
		
		\node[noborder] at ($(v1)+(0,0.75*\d)$) {\tiny};
		
	\end{tikzpicture}\\
	$G_{(v_4,v_3,v_2,v_5,v_1,v_5,v_2,v_3,v_4,v_1)}$\\
	from $\operatorname{2-d-Ring}_{5}$
\end{minipage}
\begin{minipage}{0.45\textwidth} 
	
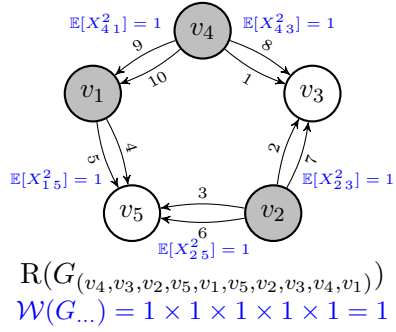
\centering
	\begin{tikzpicture}[node distance=\d and \d,>=stealth',auto, every place/.style={draw}]
		\node [black] (v1) {$v_4$};
		\node [white] (v2) [position=-30:{\d} from v1] {$v_3$};
		\node [black] (v3) [position=-108:{\d} from v2] {$v_2$};
		\node [black] (v5) [position=-150:{\d} from v1] {$v_1$};
		\node [white] (v4) [position=-72:{\d} from v5] {$v_5$};
		
		\path[->]
		(v1) edge [bend right=10] node[el,below] {$1$} (v2)
		(v1) edge [bend right=-10] node[el,above] {$8$} (v2)
		
		(v3) edge [bend right=-10] node[el,above] {$2$} (v2)
		(v3) edge [bend right=10] node[el,below] {$7$} (v2)
		
		(v3) edge [bend right=10] node[el,above] {$3$} (v4)
		(v3) edge [bend right=-10] node[el,below] {$6$} (v4)
		
		(v5) edge [bend right=-10] node[el,above] {$4$} (v4)
		(v5) edge [bend right=10] node[el,below] {$5$} (v4)
		
		(v1) edge [bend right=-10] node[el,below] {$10$} (v5)
		(v1) edge [bend right=10] node[el,above] {$9$} (v5)
		;
		
		\node[noborder] at ($(v1)+(1.2*\d,0.1*\d)$) {\tiny\textcolor{blue}{$\E[X_{4 \, 3}^2] = 1$}};
		\node[noborder] at ($(v1)+(-1.2*\d,0.1*\d)$) {\tiny\textcolor{blue}{$\E[X_{4 \, 1}^2] = 1$}};
		\node[noborder] at ($(v2)+(0.5*\d,-1.2*\d)$) {\tiny\textcolor{blue}{$\E[X_{2 \, 3}^2] = 1$}};
		\node[noborder] at ($(v5)+(-0.5*\d,-1.2*\d)$) {\tiny\textcolor{blue}{$\E[X_{1 \, 5}^2] = 1$}};
		\node[noborder] at ($(v4)+(1.0*\d,-0.5*\d)$) {\tiny\textcolor{blue}{$\E[X_{2 \, 5}^2] = 1$}};
		
	\end{tikzpicture}\\
	\vspace{-0.8cm}
	$\operatorname{R}(G_{(v_4,v_3,v_2,v_5,v_1,v_5,v_2,v_3,v_4,v_1)})$\\
	\textcolor{blue}{\small$\cW(G_{...}) = 1 \times 1 \times 1 \times 1 \times 1 = 1$}
\end{minipage}\\
\vspace{-0.5cm}\\
\begin{minipage}{0.45\textwidth}
	
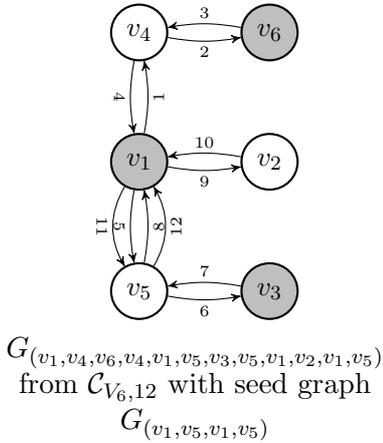
\centering
	\begin{tikzpicture}[node distance=\d and \d,>=stealth',auto, every place/.style={draw}]
		\node [black] (v1) {$v_1$};
		\node [white] (v5) [below=of v1] {$v_5$};
		\node [white] (v2) [right=of v1] {$v_2$};
		\node [white] (v4) [above=of v1] {$v_4$};
		\node [black] (v6) [right=of v4] {$v_6$};
		\node [black] (v3) [right=of v5] {$v_3$};
		
		\path[->]
		(v1) edge [bend right=10] node[el,below] {$1$} (v4)
		(v4) edge [bend right=10] node[el,below] {$4$} (v1)
		
		(v4) edge [bend right=10] node[el,below] {$2$} (v6)
		(v6) edge [bend right=10] node[el,above] {$3$} (v4)
		
		(v1) edge [bend right=10] node[el,below] {$5$} (v5)
		(v5) edge [bend right=10] node[el,below] {$8$} (v1)
		
		(v5) edge [bend right=10] node[el,below] {$6$} (v3)
		(v3) edge [bend right=10] node[el,above] {$7$} (v5)
		
		(v1) edge [bend right=10] node[el,below] {$9$} (v2)
		(v2) edge [bend right=10] node[el,above] {$10$} (v1)
		
		(v1) edge [bend right=30] node[el,below] {$11$} (v5)
		(v5) edge [bend right=30] node[el,below] {$12$} (v1)
		;
		
		\node[noborder] at ($(v4)+(0,1.6*\d)$) {\tiny};
		
	\end{tikzpicture}\\
	$G_{(v_1,v_4,v_6,v_4,v_1,v_5,v_3,v_5,v_1,v_2,v_1,v_5)}$\\
	from $\mathcal{C}_{V_6,12}$ with seed graph $G_{(v_1,v_5,v_1,v_5)}$
\end{minipage}
\begin{minipage}{0.45\textwidth}
	
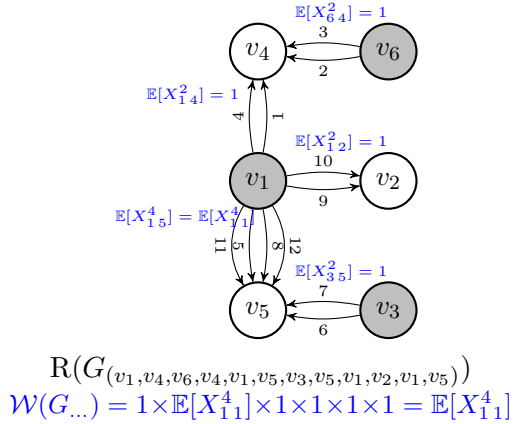
\centering
	\begin{tikzpicture}[node distance=\d and \d,>=stealth',auto, every place/.style={draw}]
		\node [black] (v1) {$v_1$};
		\node [white] (v5) [below=of v1] {$v_5$};
		\node [white] (v2) [right=of v1] {$v_2$};
		\node [white] (v4) [above=of v1] {$v_4$};
		\node [black] (v6) [right=of v4] {$v_6$};
		\node [black] (v3) [right=of v5] {$v_3$};
		
		\path[->]
		(v1) edge [bend right=10] node[el,below] {$1$} (v4)
		(v1) edge [bend right=-10] node[el,above] {$4$} (v4)
		
		(v6) edge [bend right=-10] node[el,below] {$2$} (v4)
		(v6) edge [bend right=10] node[el,above] {$3$} (v4)
		
		(v1) edge [bend right=10] node[el,below] {$5$} (v5)
		(v1) edge [bend right=-10] node[el,above] {$8$} (v5)
		
		(v3) edge [bend right=-10] node[el,below] {$6$} (v5)
		(v3) edge [bend right=10] node[el,above] {$7$} (v5)
		
		(v1) edge [bend right=10] node[el,below] {$9$} (v2)
		(v1) edge [bend right=-10] node[el,above] {$10$} (v2)
		
		(v1) edge [bend right=30] node[el,below] {$11$} (v5)
		(v1) edge [bend right=-30] node[el,above] {$12$} (v5)
		;
		
		\node[noborder] at ($(v1)+(-1.0*\d,-0.5*\d)$) {\tiny\textcolor{blue}{$\E[X_{1 \, 5}^4] = \E[X_{1 \, 1}^4]$}};
		\node[noborder] at ($(v4)+(-0.9*\d,-0.6*\d)$) {\tiny\textcolor{blue}{$\E[X_{1 \, 4}^2] = 1$}};
		\node[noborder] at ($(v6)+(-0.65*\d,0.55*\d)$) {\tiny\textcolor{blue}{$\E[X_{6 \, 4}^2] = 1$}};
		\node[noborder] at ($(v2)+(-0.65*\d,0.55*\d)$) {\tiny\textcolor{blue}{$\E[X_{1 \, 2}^2] = 1$}};
		\node[noborder] at ($(v3)+(-0.65*\d,0.55*\d)$) {\tiny\textcolor{blue}{$\E[X_{3 \, 5}^2] = 1$}};
		
	\end{tikzpicture}\\
	$\operatorname{R}(G_{(v_1,v_4,v_6,v_4,v_1,v_5,v_3,v_5,v_1,v_2,v_1,v_5)})$\\
	\textcolor{blue}{\small$\cW(G_{...}) = 1 \times \E[X_{1 \, 1}^4] \times 1 \times 1 \times 1 \times 1 = \E[X_{1 \, 1}^4]$}
\end{minipage}
\begin{proof}\
	\\
	The fact that only the seed graph is of interest for the weight is a direct consequence of the fact that removing balanced leaves from a circuit graph does not change its weight. As the original $G_{\left< \bm{i}, \bm{k} \right>}$ is assumed to have $l$ many vertices and $2l$ many edges, it is clear, that the seed graph will have $2l_0$ many edges, where $l_0$ is the number of vertices in the seed graph. This also follows from the construction of the seed graph by removing leaves, since each time a leaf is removed, we loose exactly one vertex and two edges.\\
	\\
	For the seed graph to have non-zero weight we already know by \ref{Eq_Weight2} that there can be no single edged connections between vertices and, as we know the seed graph to have $2l_0$ many edges for $l_0$ many vertices, we can see that $\operatorname{U}(S(G_{\left< \bm{i}, \bm{k} \right>}))$ (see Definition \ref{DefUndirectedConnection}) must either be
	\begin{itemize}
		\item[1)] a ring-graph
		\item[2)] a graph consisting of two vertices
		\item[3)] a path-graph with three or more vertices.
	\end{itemize}
	There are no possible seed graphs $S(G_{\left< \bm{i}, \bm{k} \right>})$ corresponding to the third option, since a seed graph must (a) have no leaves, (b) be a circuit graph and in this case (c) have only two 'extra' edges in the sense that $2l_0-2$ of its edges must be used to keep the graph $\operatorname{U}(S(G_{\left< \bm{i}, \bm{k} \right>}))$ connected. Under the assumption that $\operatorname{U}(S(G_{\left< \bm{i}, \bm{k} \right>}))$ is as in option (3) it must by (a) hold that neither of the two end-points of $S(G_{\left< \bm{i}, \bm{k} \right>})$ are leaves. By (c) this means each of the end points must have one of the two extra edges. This then makes property (b) impossible.\\
	With the properties (a-c) one similarly sees that option (2) corresponds to $\operatorname{S}(G_{\left< \bm{i}, \bm{k} \right>}) \in \operatorname{2-d-Ring}_2$. Option (1) is actually the easiest, since here there are no 'extra' edges as all $2l_0$ edges are used in the ring structure. By (b) either all edge-pairs show in the same direction or all show in opposing directions and it immediately follows that $\operatorname{S}(G_{\left< \bm{i}, \bm{k} \right>}) \in \operatorname{2-d-Ring}_{l_0}$ for some $l_0 \neq 2$ and $\operatorname{S}(G_{\left< \bm{i}, \bm{k} \right>}) \in \operatorname{1-d-Ring}_{l_0}$ for even $l_0 \geq 4$ are the only ways for option (1) to hold.
\end{proof}

\subsection{Proof of Theorem \ref{Thm_MeanExpansion}}\label{MainResult1}
\begin{proof}\
	\\
	In (\ref{Eq_TraceRepresentation4}) we had seen
	\begin{align}\label{Eq_MainResult100}
		& \E[\tr(\bm{S}^l)] = \sum\limits_{b=1}^{(l+1) \land p} \Bigg[ \sum\limits_{r=l \lor b}^{l+1} \frac{{p \choose b} {n-b \choose r-b}}{n^l} \overbrace{\sum\limits_{\substack{\bm{i} \in [b]^l, \, \bm{j} \in [r]^l \\ \bm{\{i\}} \cup \bm{\{j\}} = [r] \\ \bm{\{i\}} = [b]}} W(G_{\langle \bm{i},\bm{j} \rangle})}^{= \sum\limits_{\substack{G_{\left< \bm{i},\bm{k} \right>} \in \mathcal{C}_{[r],2l} \\ B(G_{\left< \bm{i},\bm{k} \right>}) = [b]}} \cW(G_{\left< \bm{i},\bm{k} \right>})} + \mathcal{O}\Big( \frac{p^b}{n^{b+1}} \Big) \Bigg] \ .
	\end{align}
	Lemmas \ref{WeightOfTrees} and \ref{CountingColoredTrees} yield
	\begin{align}\label{Eq_MainResult110}
		& \sum\limits_{\substack{G_{\left< \bm{i},\bm{k} \right>} \in \mathcal{C}_{[l+1],2l} \\ B(G_{\left< \bm{i},\bm{k} \right>}) = [b]}} \cW(G_{\left< \bm{i},\bm{k} \right>}) = \#\{T \in \mathcal{T}_{[l+1]} \mid B(T) = [b] \} = \mathbbm{1}_{b \leq l} \, l! \, {l-1 \choose b-1} \ .
	\end{align}
	For the $(r=l)$-term we prove
	\begin{align}\label{Eq_MainResult11}
		\sum\limits_{\substack{G_{\left< \bm{i},\bm{k} \right>} \in \mathcal{C}_{[l],2l} \\ B(G_{\left< \bm{i},\bm{k} \right>}) = \{1,...,b\}}} \cW(G_{\left< \bm{i},\bm{k} \right>}) &= \overbrace{\frac{b! (l-b)!}{2} \Bigg( {2l \choose 2b} + (2b-1) {l \choose b}^2 \Bigg)}^{=: A_1(l,b)} \nonumber\\
		& \hspace{1cm} + (\E[X_{1 \, 1}^4]-3) \underbrace{\mathbbm{1}_{b<l} \, b! (l-b)! {l \choose b-1}{l \choose b+1}}_{=: A_2(l,b)}
	\end{align}
	first for $b=l$ and then for $b<l$. By Proposition \ref{WeightOfFusedTrees} and Lemmas \ref{RingColor_OneDirectional} and \ref{RingColor_TwoDirectional} the only possible seed graph of $G_{\left< \bm{i},\bm{k} \right>}$ with positive weight for $b=l$ is $\operatorname{S}(G_{\left< \bm{i},\bm{k} \right>}) = G_{(v,v)} \in \operatorname{2-d-Ring}_1$ (for some $v \in [l]$), since there are no white vertices in $G_{\left< \bm{i},\bm{k} \right>}$ and thus no white vertices in the seed graph. By applying Corollary \ref{SproutingCorollary2} for $l_0=1$, $b'=l-1$ and $w'=0$ we get
	\begin{align*}
		& \sum\limits_{\substack{G_{\left< \bm{i},\bm{k} \right>} \in \mathcal{C}_{[l],2l} \\ B(G_{\left< \bm{i},\bm{k} \right>}) = [l]}} \cW(G_{\left< \bm{i},\bm{k} \right>}) = \#\{G \in \mathcal{C}_{[l],2l} \mid \operatorname{S}(G) \in \operatorname{2-d-Ring}_{1}, \, B(G) = [l] \}\\
		& = 1 \, l! \, 0! {l \choose l-1} {l \choose 0} = l \, l! = A(l,l) + (\E[X_{1 \, 1}^4]-3) \underbrace{B(l,l)}_{=0} \ .
	\end{align*}
	It remains to show (\ref{Eq_MainResult11}) for $b<l$. By Proposition \ref{WeightOfFusedTrees} we have
	\begin{align*}
		& \sum\limits_{\substack{G_{\left< \bm{i},\bm{k} \right>} \in \mathcal{C}_{[l],2l} \\ B(G_{\left< \bm{i},\bm{k} \right>}) = [b]}} \cW(G_{\left< \bm{i},\bm{k} \right>})\\
		& = \#\{G \in \mathcal{C}_{[l],2l} \mid \operatorname{S}(G) \in \operatorname{2-d-Ring}_{l_0}\text{ for some }l_0 \neq 2 , \, B(G) = [b] \}\\
		& \ \ \ \ + \E[X_{1 \, 1}^4] \, \#\{G \in \mathcal{C}_{[l],2l} \mid \operatorname{S}(G) \in \operatorname{2-d-Ring}_{2} , \, B(G) = [b] \}\\
		& \ \ \ \ + \#\{G \in \mathcal{C}_{[l],2l} \mid \operatorname{S}(G) \in \operatorname{1-d-Ring}_{l_0}\text{ for some even } l_0 \geq 4 , \, B(G) = [b] \} \ .
	\end{align*}
	We add zero in the form of $\pm \#\{G \in \mathcal{C}_{[l],2l} \mid \operatorname{S}(G) \in \operatorname{2-d-Ring}_{2} , \, B(G) = [b] \}$ three times. This will later make the expressions easier to handle. Heuristically this counteracts the exclusion of the case $l_0=2$ from both other cardinalities. Write
	\begin{align*}
		& \sum\limits_{\substack{G_{\left< \bm{i},\bm{k} \right>} \in \mathcal{C}_{[l],2l} \\ B(G_{\left< \bm{i},\bm{k} \right>}) = \{1,...,b\}}} \cW(G_{\left< \bm{i},\bm{k} \right>})\\
		& = \#\{G \in \mathcal{C}_{[l],2l} \mid \operatorname{S}(G) \in \operatorname{2-d-Ring}_{l_0}\text{ for some odd }l_0 , \, B(G) = [b] \}\\
		& \ \ \ \ + \#\{G \in \mathcal{C}_{[l],2l} \mid \operatorname{S}(G) \in \operatorname{2-d-Ring}_{l_0}\text{ for some even }l_0 \neq 2 , \, B(G) = [b] \}\\
		& \ \ \ \ + 2 \, \#\{G \in \mathcal{C}_{[l],2l} \mid \operatorname{S}(G) \in \operatorname{2-d-Ring}_{2} , \, B(G) = [b] \}\\
		& \ \ \ \ + \#\{G \in \mathcal{C}_{[l],2l} \mid \operatorname{S}(G) \in \operatorname{1-d-Ring}_{l_0}\text{ for some even } l_0 \geq 4 , \, B(G) = [b] \}\\
		& \ \ \ \ + \#\{G \in \mathcal{C}_{[l],2l} \mid \operatorname{S}(G) \in \operatorname{2-d-Ring}_{2} , \, B(G) = [b] \}\\
		& \ \ \ \ + (\E[X_{1 \, 1}^4]-3) \, \#\{G \in \mathcal{C}_{[l],2l} \mid \operatorname{S}(G) \in \operatorname{2-d-Ring}_{2} , \, B(G) = [b] \} \ .
	\end{align*}
	Let $w := l-b$ be the number of white vertices in $V_l$. We sum over all possible ring lengths $l_0$, which are precisely those for which we have enough black and white vertices to satisfy the coloring from Lemmas \ref{RingColor_OneDirectional} and \ref{RingColor_TwoDirectional} for $\operatorname{S}(G)$. For $b<l$, i.e. $w>0$, we get
	\begin{align*}
		\sum\limits_{\substack{G_{\left< \bm{i},\bm{k} \right>} \in \mathcal{C}_{[l],2l} \\ B(G_{\left< \bm{i},\bm{k} \right>}) = [b]}} \cW(G_{\left< \bm{i},\bm{k} \right>}) & = \sum\limits_{\substack{m=1 \\ m \leq w+1}}^{b} \#\{G \in \mathcal{C}_{[l],2l} \mid \operatorname{S}(G) \in \operatorname{2-d-Ring}_{2m-1}, \, B(G) = [b] \}\\
		& \ \ \ \ + \sum\limits_{\substack{m=2 \\ m \leq w}}^{b} \#\{G \in \mathcal{C}_{[l],2l} \mid \operatorname{S}(G) \in \operatorname{2-d-Ring}_{2m}, \, B(G) = [b] \}\\
		& \ \ \ \ + 2 \, \#\{G \in \mathcal{C}_{[l],2l} \mid \operatorname{S}(G) \in \operatorname{2-d-Ring}_{2} , \, B(G) = [b] \}\\
		& \ \ \ \ + \sum\limits_{\substack{m=2 \\ m \leq w}}^{b} \#\{G \in \mathcal{C}_{[l],2l} \mid \operatorname{S}(G) \in \operatorname{1-d-Ring}_{2m} , \, B(G) = [b] \}\\
		& \ \ \ \ + \#\{G \in \mathcal{C}_{[l],2l} \mid \operatorname{S}(G) \in \operatorname{2-d-Ring}_{2} , \, B(G) = [b] \}\\
		& \ \ \ \ + (\E[X_{1 \, 1}^4]-3) \, \#\{G \in \mathcal{C}_{[l],2l} \mid \operatorname{S}(G) \in \operatorname{2-d-Ring}_{2} , \, B(G) = [b] \} \ .
	\end{align*}
	With $b' = b-\lceil \frac{l_0}{2} \rceil$, $w' = l-b-\lfloor \frac{l_0}{2} \rfloor$ and Corollaries \ref{SproutingCorollary1} and \ref{SproutingCorollary2} this becomes
	\begin{align*}
		& \sum\limits_{\substack{G_{\left< \bm{i},\bm{k} \right>} \in \mathcal{C}_{[l],2l} \\ B(G_{\left< \bm{i},\bm{k} \right>}) = [b]}} \cW(G_{\left< \bm{i},\bm{k} \right>})\\
		& = \sum\limits_{m=1}^{b \land (w+1)} (2m-1) \, b! \, w! {l \choose b-m}{l \choose l-b-m+1}\\
		& \ \ \ \ + \sum\limits_{m=2}^{b \land w} 2m \, b! \, w! {l \choose b-m}{l \choose l-b-m} + 2 \, b! \, w! {l \choose b-1}{l \choose l-b-1}\\
		& \ \ \ \ + \sum\limits_{m=2}^{b \land w} b! \, w! {l \choose b-m}{l \choose l-b-m} + b! \, w! {l \choose b-1}{l \choose l-b-1}\\
		& \ \ \ \ + (\E[X_{1 \, 1}^4]-3) \, b! \, w! {l \choose b-1}{l \choose l-b-1}\\
		& = \sum\limits_{m=1}^{b \land (l-b+1)} (2m-1) \, b! \, w! {l \choose b-m}{l \choose l-b-m+1}\\
		& \ \ \ \ + \sum\limits_{m=1}^{b \land (l-b)} 2m \, b! \, w! {l \choose b-m}{l \choose l-b-m}\\
		& \ \ \ \ + \sum\limits_{m=1}^{b \land (l-b)} b! \, w! {l \choose b-m}{l \choose l-b-m} + (\E[X_{1 \, 1}^4]-3) \, 2 \, b! \, w! {l \choose b-1}{l \choose l-b-1}\\
		& = b! \, w! \sum\limits_{m=1}^{b \land (l-b+1)} (2m-1) {l \choose b-m}{l \choose b+m-1}\\
		& \ \ \ \ + b! \, w! \sum\limits_{m=1}^{b \land (l-b)} (2m+1) {l \choose b-m}{l \choose b+m}\\
		& \ \ \ \ + (\E[X_{1 \, 1}^4]-3) \, b! \, w! {l \choose b-1}{l \choose b+1} \ .
	\end{align*}
	This can further be simplified thanks to Peter Taylor's answer to our question on Math Overflow (see \cite{PeterTaylor}). Peter Taylor shows that the upper two sums are together
	\begin{align*}
		& \frac{b! w!}{2} \left( {2l \choose 2b} + (2b-1) {l \choose b}^2 \right) \ , 
	\end{align*}
	which means for $b<l$ we get
	\begin{align*}
		& \sum\limits_{\substack{G_{\left< \bm{i},\bm{k} \right>} \in \mathcal{C}_{[l],2l} \\ B(G_{\left< \bm{i},\bm{k} \right>}) = [b]}} \cW(G_{\left< \bm{i},\bm{k} \right>}) = \underbrace{\frac{b! w!}{2} \left( {2l \choose 2b} + (2b-1) {l \choose b}^2 \right)}_{=: A_1(l,b)} + (\E[X_{1 \, 1}^4]-3) \, \underbrace{b! \, w! {l \choose b-1}{l \choose b+1}}_{=: A_2(l,b)} \ .
	\end{align*}
	We have thus shown Equality (\ref{Eq_MainResult11}), which with (\ref{Eq_MainResult100}) and (\ref{Eq_MainResult110}) proves Theorem \ref{Thm_MeanExpansion}.
\end{proof}

\section{Proving Theorem \ref{Thm_CovExpansion}}\label{SectionMainRes2}

\subsection{Definition (Double-circuit graph)}\label{DefDoubleCircuit}
For any vertex set $V_r$ we call a pair $(\bm{i}_1,\bm{i}_2)$ of sequences $\bm{i}_1 \in V_r^{N_1}$, $\bm{i}_2 \in V_r^{N_2}$ with $\bm{\{i}_1\bm{\}} \cup \bm{\{i}_2\bm{\}} = V_r$ a \textit{double-route} over $V_r$ with lengths $N_1,N_2 \in \N$.\\
The pair of graphs $(G_{\bm{i}_1},G_{\bm{i}_2})$, with both $G_{\bm{i}_1},G_{\bm{i}_2}$ as in Definition \ref{DefCircuitGraph}, is called a \textit{double-circuit graph} on $V_r$ of lengths $N_1,N_2$. We interpret the pair as a single 'combined' directed graph defined by $(V_r,([N_1] \times \{1\}) \sqcup ([N_2] \times \{2\}), f_{(\bm{i}_1,\bm{i}_2)})$, where
\begin{align*}
	& f_{(\bm{i}_1,\bm{i}_2)}((\cdot,1)) := f_{G_1} = f_{\bm{i}_1} : [N_1] \rightarrow V_r \times V_r \ \ ; \ \ f_{(\bm{i}_1,\bm{i}_2)}((\cdot,2)) := f_{G_2} = f_{\bm{i}_2} : [N_2] \rightarrow V_r \times V_r \ .
\end{align*}
Let $\mathcal{C}^2_{V_{l},N_1,N_2}$ denote the set of all such double-circuit graphs.\\
\\
The coloring of a double-circuit graph is then defined such that vertices are white, iff they are white in both $G_{\bm{i}_1}$ and in $G_{\bm{i}_2}$, i.e.
\begin{align*}
	& B((G_{\bm{i}_1},G_{\bm{i}_2})) := B(G_{\bm{i}_1}) \cup B(G_{\bm{i}_2}) \ .
\end{align*}
Since a balanced leaf of $(G_{\bm{i}_1},G_{\bm{i}_2})$ is always either in  $(G_{\bm{i}_1}$ or $G_{\bm{i}_2})$, we can canonically define the seed graph $\operatorname{S}((G_{\bm{i}_1},G_{\bm{i}_2}))$ by iterative removal of balanced leaves. In fact Lemma (\ref{SproutingDoubleCondition}) is a trivial consequence of this property.

\subsection{Definition (Double ring-type graphs)}\label{DefDoubleRings}
For $l_0 \geq 3$ a double-circuit graph $(G_{\bm{i}_1},G_{\bm{i}_2}) \in \mathcal{C}^2_{V_{l_0},l_0,l_0}$ will be called of \textit{one-directional double ring-type} with ring-length $l_0$, if $\bm{i}_1$ and $\bm{i}_2$ are the same modulo starting position.\\
On the other hand for any $l_0 \in \N$ a double-circuit graph $(G_{\bm{i}_1},G_{\bm{i}_2}) \in \mathcal{C}^2_{V_{l_0},E_{l_0,1},E_{l_0,2}}$ will be called of \textit{two-directional double ring-type} with ring-length $l_0$, if $\bm{i}_1$ and the backward-sequence $((\bm{i}_2)_{l_0+1-i})_{i \leq l_0}$ are equal modulo starting position.\\
\\
Both of these definitions already imply that each vertex also has in- and out-degree $1$, that both graphs are exhaustive and that $G_{\bm{i}_1},G_{\bm{i}_2} \in \mathcal{C}_{V_{l_0},l_0}$. We could thus equally well have chosen $(G_{\bm{i}_1},G_{\bm{i}_2})$ from $\mathcal{C}_{V_{l_0},l_0} \times \mathcal{C}_{V_{l_0},l_0}$.\\
\\
For even $l_0 \geq 4$ and any $(G_{\bm{i}_1},G_{\bm{i}_2}) \in \mathcal{C}_{V_{l_0},l_0} \times \mathcal{C}_{V_{l_0},l_0}$ we say
$$
(G_{\bm{i}_1},G_{\bm{i}_2}) \in \operatorname{Double-1-d-Ring}_{l_0} \ ,
$$
if $(G_{\bm{i}_1},G_{\bm{i}_2})$ is of one-directional double ring-type with ring-length $l_0$ and the starting positions of $\bm{i}_1$ and $\bm{i}_2$ are an even number of steps apart. Analogously for even $l_0 \in \N$ we say
$$
(G_{\bm{i}_1},G_{\bm{i}_2}) \in \operatorname{Double-2-d-Ring}_{l_0} \ ,
$$
if $(G_{\bm{i}_1},G_{\bm{i}_2})$ is of two-directional double ring-type with ring-length $l_0$ and the starting positions of $\bm{i}_1$ and $\bm{i}_2$ are also an even number of steps apart.\\
\\
\begin{minipage}{0.45\textwidth} 
	\centering
	\begin{tikzpicture}[node distance=\d and \d,>=stealth',auto, every place/.style={draw}]
		\node [black] (v1) {$v_1$};
		
		\path (v1) edge [in=150,out=120,loop,blue] node[el,above] {$1$} (v1);
		\path (v1) edge [in=100,out=70,loop,red] node[el,above] {$1$} (v1);
	\end{tikzpicture}\\
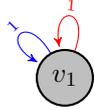
	$(\textcolor{blue}{G_{(v_1)}},\textcolor{red}{G_{(v_1)}})$, the only\\
	two-dir. double ring-type graph\\
	with ring-length $1$
\end{minipage}
\begin{minipage}{0.45\textwidth}
	\centering
	\begin{tikzpicture}[node distance=\d and \d,>=stealth',auto, every place/.style={draw}]
		\node [black] (v1) {$v_1$};
		\node [black] (v2) [right=of v1] {$v_2$};
		\node [black] (v3) [below=of v1] {$v_3$};
		\node [black] (v4) [below=of v2] {$v_4$};
		
		\path[->]
		(v1) edge [bend right=10,blue] node[el,below] {$1$} (v2)
		(v2) edge [bend right=10,blue] node[el,below] {$2$} (v4)
		(v4) edge [bend right=10,blue] node[el,above] {$3$} (v3)
		(v3) edge [bend right=10,blue] node[el,below] {$4$} (v1)
		(v1) edge [bend right=-10,red] node[el,above] {$2$} (v2)
		(v2) edge [bend right=-10,red] node[el,above] {$3$} (v4)
		(v4) edge [bend right=-10,red] node[el,below] {$4$} (v3)
		(v3) edge [bend right=-10,red] node[el,above] {$1$} (v1)
		;
		
	\end{tikzpicture}\\
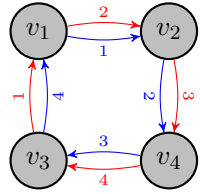
	$(\textcolor{blue}{G_{(v_1,v_2,v_4,v_3)}},\textcolor{red}{G_{(v_1,v_2,v_4,v_3)}})$, a\\
	one-dir. double ring-type graph\\
	\underline{not} in $\operatorname{Double-1-d-Ring}_{4}$
\end{minipage}
\\
\\
\\
\begin{minipage}{0.45\textwidth}
	\centering
	\begin{tikzpicture}[node distance=\d and \d,>=stealth',auto, every place/.style={draw}]
		\node [black] (v1) {$v_1$};
		\node [white] (v2) [right=of v1] {$v_2$};
		\node [white] (v3) [below=of v1] {$v_3$};
		\node [black] (v4) [below=of v2] {$v_4$};
		
		\path[->]
		(v1) edge [bend right=10,blue] node[el,below] {$1$} (v2)
		(v2) edge [bend right=10,blue] node[el,below] {$2$} (v4)
		(v4) edge [bend right=10,blue] node[el,above] {$3$} (v3)
		(v3) edge [bend right=10,blue] node[el,below] {$4$} (v1)
		(v1) edge [bend right=-10,red] node[el,above] {$1$} (v2)
		(v2) edge [bend right=-10,red] node[el,above] {$2$} (v4)
		(v4) edge [bend right=-10,red] node[el,below] {$3$} (v3)
		(v3) edge [bend right=-10,red] node[el,above] {$4$} (v1)
		;
		
	\end{tikzpicture}\\
	$(\textcolor{blue}{G_{(v_1,v_2,v_4,v_3)}},\textcolor{red}{G_{(v_1,v_2,v_4,v_3)}})$\\
	from $\operatorname{Double-1-d-Ring}_{4}$
\end{minipage}
\begin{minipage}{0.45\textwidth} 
	\centering
	\begin{tikzpicture}[node distance=\d and \d,>=stealth',auto, every place/.style={draw}]
		\node [white] (v1) {$v_1$};
		\node [black] (v2) [right=of v1] {$v_4$};
		\node [white] (v3) [below right=of v2] {$v_3$};
		\node [black] (v4) [below left=of v3] {$v_2$};
		\node [white] (v5) [left=of v4] {$v_6$};
		\node [black] (v6) [above left=of v5] {$v_5$};
		
		\path[->]
		(v1) edge [bend right=10,blue] node[el,below] {$2$} (v2)
		(v2) edge [bend right=10,red] node[el,above] {$1$} (v1)
		
		(v2) edge [bend right=10,blue] node[el,below] {$3$} (v3)
		(v3) edge [bend right=10,red] node[el,above] {$6$} (v2)
		
		(v3) edge [bend right=10,blue] node[el,above] {$4$} (v4)
		(v4) edge [bend right=10,red] node[el,below] {$5$} (v3)
		
		(v4) edge [bend right=10,blue] node[el,above] {$5$} (v5)
		(v5) edge [bend right=10,red] node[el,below] {$4$} (v4)
		
		(v5) edge [bend right=10,blue] node[el,above] {$6$} (v6)
		(v6) edge [bend right=10,red] node[el,below] {$3$} (v5)
		
		(v6) edge [bend right=10,blue] node[el,below] {$1$} (v1)
		(v1) edge [bend right=10,red] node[el,above] {$2$} (v6)
		
		;
		
	\end{tikzpicture}\\
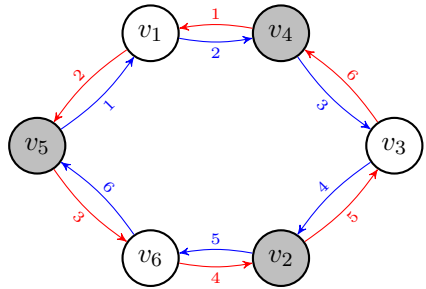
	$(\textcolor{blue}{G_{(v_5,v_1,v_4,v_3,v_2,v_6)}},\textcolor{red}{G_{(v_4,v_1,v_5,v_6,v_2,v_3)}})$\\
	from $\operatorname{Double-2-d-Ring}_{6}$
\end{minipage}

\subsection{Lemma (Characterizing sprouts of double ring-type graphs)}\label{SproutingDoubleCondition}
For any $(G_1,G_2) \in \mathcal{C}^2_{V_l,2l_1,2l_2}$, where $l=l_1+l_2$, and $l_0 > 2$ the following two statements are equivalent.
\begin{itemize}
	\item[a)] $\operatorname{S}((G_1,G_2)) \in \operatorname{Double-1-d-Ring}_{l_0}$
	
	\item[b)] $(\operatorname{S}(G_1),\operatorname{S}(G_2)) \in \operatorname{Double-1-d-Ring}_{l_0}$ and the sets
	\begin{align*}
		& V'_1 := V(G_1) \setminus V(\operatorname{S}(G_1)) \ \ \text{ and } \ \ V'_2 := V(G_2) \setminus V(\operatorname{S}(G_2))
	\end{align*}
	are disjoint.
\end{itemize}
The same holds for $\operatorname{Double-2-d-Ring}_{l_0}$ instead of $\operatorname{Double-1-d-Ring}_{l_0}$.
\begin{proof}
	This follows directly from the arguments at the end of Definition (\ref{DefDoubleCircuit}).
\end{proof}

\subsection{Lemma (Coloring of double ring-type graphs)}\label{DoubleRingColoring}
For even $l_0 \in \N$ and any $(G_{\bm{i}_1},G_{\bm{i}_2}) \in \mathcal{C}_{V_{l_0},l_0} \times \mathcal{C}_{V_{l_0},l_0,2}$ with
\begin{align*}
	& (G_{\bm{i}_1},G_{\bm{i}_2}) \in \operatorname{Double-1-d-Ring}_{l_0} \ \ \text{ or } \ \ (G_{\bm{i}_1},G_{\bm{i}_2}) \in \operatorname{Double-2-d-Ring}_{l_0}
\end{align*}
the vertices will alternate between black and white along the ring-structure of $(G_{\bm{i}_1},G_{\bm{i}_2})$.
\begin{proof}
	This trivially follows from the assumptions that $l_0$ is even and that the starting positions of the routes $\bm{i}_1$ and $\bm{i}_2$ are an even number of steps apart.
\end{proof}

\subsection{Lemma (Counting sprouts of double ring-type graphs)}\label{CountingSproutingDoubleRings}
For even $l_0 \in \N$ and any $b'_1,b'_2,w'_1,w'_2 \geq 0$ define $l_1 := \frac{l_0}{2}+b'_1+w'_1$, $l_2 := \frac{l_0}{2}+b'_2+w'_2$ as well as $l=l_1+l_2$. Let $B$ be a subset of $V_l$ such that $b := \#B = \frac{l_0}{2} + b'_1+b'_2$, then the number of $(G_1 := G_{\bm{i}_1},G_2 := G_{\bm{i}_2}) \in \mathcal{C}^2_{V_l,2l_1,2l_2}$ with $B(G) = B$, $\operatorname{S}((G_1,G_2)) \in \operatorname{Double-1-d-Ring}_{l_0}$ and the properties
\begin{align}\label{Eq_CountingSproutingDoubleRings1}
	& b'_1 = \#\overbrace{(V(G_1) \cap B) \setminus V(\operatorname{S}(G_1))}^{=: B'_1} \ \ ; \ \ b'_2 = \#\overbrace{(V(G_2) \cap B) \setminus V(\operatorname{S}(G_2))}^{=: B'_2} \nonumber\\
	& w'_1 = \#\underbrace{(V(G_1) \setminus B) \setminus V(\operatorname{S}(G_1))}_{=: W'_1} \ \ ; \ \ w'_2 = \#\underbrace{(V(G_2) \setminus B) \setminus V(\operatorname{S}(G_2))}_{=: W'_2}
\end{align}
is given by
\begin{align*}
	& \mathbbm{1}_{l_0 \geq 4} \, \frac{l_0}{2} \, b! \, w! \, {l_1 \choose b'_1} {l_1 \choose w'_1} {l_2 \choose b'_2} {l_2 \choose w'_2} \ ,
\end{align*}
where $w=l-b$. Analogously the number of $(G_1,G_2) \in \mathcal{C}^2_{V_l,2l_1,2l_2}$ with $B(G) = B$, $\operatorname{S}((G_1,G_2)) \in \operatorname{Double-2-d-Ring}_{l_0}$ and (\ref{Eq_CountingSproutingDoubleRings1}) is also given by
\begin{align*}
	& \frac{l_0}{2} \, b! \, w! {l_1 \choose b'_1} {l_1 \choose w'_1} {l_2 \choose b'_2} {l_2 \choose w'_2} \ .
\end{align*}\
\\
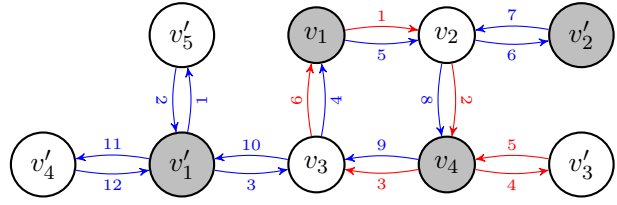
\begin{minipage}{0.95\textwidth}
	\centering
	\begin{tikzpicture}[node distance=\d and \d,>=stealth',auto, every place/.style={draw}]
		\node [black] (v1) {$v_1$};
		\node [white] (v2) [right=of v1] {$v_2$};
		\node [white] (v3) [below=of v1] {$v_3$};
		\node [black] (v4) [below=of v2] {$v_4$};
		\node [black] (vv1) [left=of v3] {$v_1'$};
		\node [black] (vv2) [right=of v2] {$v_2'$};
		\node [white] (vv3) [right=of v4] {$v_3'$};
		\node [white] (vv4) [left=of vv1] {$v_4'$};
		\node [white] (vv5) [left=of v1] {$v_5'$};
		
		\path[->]
		(v1) edge [bend right=10,blue] node[el,below] {$5$} (v2)
		(v2) edge [bend right=10,blue] node[el,below] {$8$} (v4)
		(v4) edge [bend right=10,blue] node[el,above] {$9$} (v3)
		(v3) edge [bend right=10,blue] node[el,below] {$4$} (v1)
		
		(v2) edge [bend right=10,blue] node[el,below] {$6$} (vv2)
		(vv2) edge [bend right=10,blue] node[el,above] {$7$} (v2)
		
		(v3) edge [bend right=10,blue] node[el,above] {$10$} (vv1)
		(vv1) edge [bend right=10,blue] node[el,below] {$3$} (v3)
		
		(vv5) edge [bend right=10,blue] node[el,below] {$2$} (vv1)
		(vv1) edge [bend right=10,blue] node[el,below] {$1$} (vv5)
		
		(vv4) edge [bend right=10,blue] node[el,below] {$12$} (vv1)
		(vv1) edge [bend right=10,blue] node[el,above] {$11$} (vv4)

		(v1) edge [bend right=-10,red] node[el,above] {$1$} (v2)
		(v2) edge [bend right=-10,red] node[el,above] {$2$} (v4)
		(v4) edge [bend right=-10,red] node[el,below] {$3$} (v3)
		(v3) edge [bend right=-10,red] node[el,above] {$6$} (v1)
		
		(v4) edge [bend right=10,red] node[el,below] {$4$} (vv3)
		(vv3) edge [bend right=10,red] node[el,above] {$5$} (v4)
		;
		
	\end{tikzpicture}
	$\big( \textcolor{blue}{G_{(v'_1,v'_5,v'_1,v_3,v_1,v_2,v'_2,v_2,v_4,v_3,v'_1,v'_4)}}, \textcolor{red}{G_{(v_1,v_2,v_4,v'_3,v_4,v_3)}} \big)$from $\mathcal{C}_{V_4 \sqcup V'_{5}, 12,6}$\\
	\vspace{0.3cm}
	$S\big(\textcolor{blue}{G_{(...)}}, \textcolor{red}{G_{(...)}}\big) = \big(\textcolor{blue}{G_{(v_1,v_2,v_3,v_4)}}, \textcolor{red}{G_{(v_1,v_2,v_3,v_4)}}\big)$\\
	$S(\textcolor{blue}{G_{(...)}}) = \textcolor{blue}{G_{(v_1,v_2,v_3,v_4)}} \ ; \ S(\textcolor{red}{G_{(...)}}) = \textcolor{red}{G_{(v_1,v_2,v_3,v_4)}}$
\end{minipage}
\begin{proof}\
	\\
	It is easily seen that the tree structures sprouting from the seed graph will each either belong to $G_1$ or $G_2$ with no intermixing, meaning we can count the sprouts of $S(G_1)$ and of $S(G_2)$ separately. The sets $B_i'$ and $W_i'$ describe the sets of sprouted black and white vertices visited by (only) $G_i$, while $V(S(G_1)) = V(S(G_2))$ is the set of vertices in the ring structure. There are clearly ${b \choose b'_1 , b'_2, \frac{l_0}{2}}$ many ways to distribute the total number of $b$ black vertices among the sets $B_1',B_2', V(S(G_1)) \cap B$ and also ${w \choose w'_1,w'_2,\frac{l_0}{2}}$ ways to distribute the total number of $w$ white vertices among the sets $W_1',W_2', V(S(G_1)) \setminus B$. Further there are $\frac{l_0}{2}!$ choices each for the order of the black/white vertices in the ring structure (with starting position of $S(G_1)$ accounted for) and finally $\frac{l_0}{2}$ many choices for the starting position of $S(G_2)$. We have now fixed the seed graphs and the sets $B_i',W_i'$ and may use Proposition \ref{CountingSprouting} to see that there remain $\frac{(l_0+b_i'+w_i')!^2}{(l_0+b_i')! (l_0+w_i')!}$ choices for $G_i$ in both cases $i \in \{1,2\}$. Multiplication of all the listed factors leads to the wanted formula.\\
	\\
	The case $l_0=2$ must be viewed separately, since in this case both $G_1$ and $G_2$ are trees. Again we have ${b \choose b'_1 , b'_2, \frac{l_0}{2}} \, {w \choose w'_1,w'_2,\frac{l_0}{2}}$ many ways of assigning the vertices and now Lemma \ref{CountingColoredTrees} tells us that there are $l_1! {l_1-1 \choose b'_1}$ many choices for $G_1$ as well as $l_2! {l_2-1 \choose b'_2}$ many choices for $G_2$ with the correct vertex assignment. A symmetry argument yields that the portion of choices of $(G_1,G_2)$, which actually have a connection between $v_1$ and $v_2$, is $\frac{l_1}{(b'_1+1)(w'_1+1)}\frac{l_2}{(b'_2+1)(w'_2+1)}$. Multiplication of all the listed factors again yields the wanted formula.
\end{proof}

\subsection{Proposition (Covariance-weight of double graphs)}\label{WeightCov}
For any $l_1,l_2 \in \N$ define $l := l_1+l_2$ and a vertex set $V_l$. For $\bm{i},\bm{k} \in V_{l}^{l_1}$ and $\bm{j},\bm{m} \in V_{l}^{l_2}$ with
\begin{align*}
	& \#( \bm{\{i\}} \cup \bm{\{k\}} \cup \bm{\{j\}} \cup \bm{\{m\}}) \geq l
\end{align*}
we have
\begin{align*}
	& \overbrace{\E\left[ (X_{i_1 \, k_1} X_{i_2 \, k_1}) ... (X_{i_{l_1} \, k_{l_1}} X_{i_1 \, k_{l_1}}) \cdot (X_{j_1 \, m_1} X_{j_2 \, m_1}) ... (X_{j_{l_2} \, m_{l_2}} X_{j_1 \, m_{l_2}}) \right]}^{=: \cW_{\bm{i},\bm{k},\bm{j},\bm{m}}}\\
	& \ \ \ \ - \underbrace{\E\left[ (X_{i_1 \, k_1} X_{i_2 \, k_1}) ... (X_{i_{l_1} \, k_{l_1}} X_{i_1 \, k_{l_1}})\right]}_{=: \cW_{\bm{i},\bm{k}}} \, \underbrace{\E\left[(X_{j_1 \, m_1} X_{j_2 \, m_1}) ... (X_{j_{l_2} \, m_{l_2}} X_{j_1 \, m_{l_2}}) \right]}_{=: \cW_{\bm{j},\bm{m}}}\\
	& = \begin{cases}
		0 & \text{, if } \#( \bm{\{i\}} \cup \bm{\{k\}} \cup \bm{\{j\}} \cup \bm{\{m\}}) > l\\
		\E[X_{1 \, 1}^4]-1 & \text{, if $\operatorname{S}((G_{\left< \bm{i},\bm{k} \right>},G_{\left< \bm{j},\bm{m} \right>})) \in \operatorname{Double-2-d-Ring}_{2}$}\\
		1 & \text{, if $\operatorname{S}((G_{\left< \bm{i},\bm{k} \right>},G_{\left< \bm{j},\bm{m} \right>})) \in \operatorname{Double-1-d-Ring}_{l_0}$ for even $l_0 \geq 4$}\\
		1 & \text{, if $\operatorname{S}((G_{\left< \bm{i},\bm{k} \right>},G_{\left< \bm{j},\bm{m} \right>})) \in \operatorname{Double-2-d-Ring}_{l_0}$ for even $l_0 \geq 4$}\\
		0 & \text{else}
	\end{cases} \ .
\end{align*}
\begin{proof}\
	\\
	The proof of this proposition is tedious and contributes little to the understanding of our methods. The idea is the same as in Proposition \ref{WeightOfFusedTrees}, where we can additionally use the fact that we are only looking for cases where the joint weight $\cW_{\bm{i},\bm{k},\bm{j},\bm{m}}$ differs from the product of the single weights $\cW_{\bm{i},\bm{k}} \cW_{\bm{j},\bm{m}}$. The product form (\ref{Eq_Weight2}) - which can similarly be applied to $\cW_{\bm{i},\bm{k},\bm{j},\bm{m}}$ - implies that the graphs $G_{\langle \bm{i},\bm{k} \rangle}$ and $G_{\langle \bm{j},\bm{m} \rangle}$ must share at least one connection between vertices, which makes the case $\#( \bm{\{i\}} \cup \bm{\{k\}} \cup \bm{\{j\}} \cup \bm{\{m\}}) > l$ impossible. Lemma \ref{SproutingDoubleCondition} can be used to show show that in the remaining relevant cases we still have $S\big((G_1,G_2)\big) = \big(S(G_1),S(G_2)\big)$ and the rest follows analogously to Proposition \ref{WeightOfFusedTrees}.
\end{proof}

\subsection{Proof of Theorem \ref{Thm_CovExpansion}}\label{MainResult2}
\begin{proof}\
	\\
	As all formulas in the formulation of Theorem \ref{Thm_CovExpansion} are symmetric in $l_1$ and $l_2$, without loss of generality assume $l_1 \leq l_2$. We can express the covariance as a weighted sum over graphs with
	\begin{align}\label{Eq_CovSum1}
		& \operatorname{Cov}[\tr(\bm{S}_{p,n}^{l_1}), \tr(\bm{S}_{p,n}^{l_2})] = \frac{1}{n^{l_1+l_2}} \sum\limits_{\substack{\bm{i} \in [p]^{l_1}, \bm{k} \in [n]^{l_1} \\ \bm{j} \in [p]^{l_2}, \bm{m} \in [n]^{l_2}}} \big( \cW_{\bm{i},\bm{k},\bm{j},\bm{m}} - \cW_{\bm{i},\bm{k}} \, \cW_{\bm{j},\bm{m}} \big)
	\end{align}
	for $\cW_{\bm{i},\bm{k}}$ as in (\ref{Eq_TraceRepresentation1}) or (\ref{Eq_Weight2}). Similarly $\cW_{\bm{i},\bm{k},\bm{j},\bm{m}}$ is given by
	\begin{align}
		\cW_{\bm{i},\bm{k},\bm{j},\bm{m}} &= \E\big[(X_{i_1 \, k_1} X_{i_2 \, k_1}) \, (X_{i_2 \, k_2} X_{i_3 \, k_2}) \, \dots \, (X_{i_l \, k_l} X_{i_1\, k_l}) \nonumber\\
		& \hspace{1cm} \times (X_{j_1 \, m_1} X_{j_2 \, m_1}) \, (X_{j_2 \, m_2} X_{j_3 \, m_2}) \, \dots \, (X_{j_l \, m_l} X_{j_1\, m_l})\big] \nonumber\\
		& = \prod\limits_{(s,t) \in [p] \times [n]} \E\big[X_{s,t}^{A_{s,t}(\operatorname{R}(G_{\langle \bm{i},\bm{k} \rangle})) + A_{s,t}(\operatorname{R}(G_{\langle \bm{j},\bm{m} \rangle}))}\big] \ .
	\end{align}
	We split the sum in (\ref{Eq_CovSum1}) again by number of total vertices $r$ and number of black vertices $b$ to get
	\begin{align}\label{Eq_CovSum2}
		& \operatorname{Cov}[\tr(\bm{S}_{p,n}^{l_1}), \tr(\bm{S}_{p,n}^{l_2})] = \frac{1}{n^{l_1+l_2}} \sum\limits_{r=1}^{(2l_1+2l_2) \land n} \sum\limits_{b=1}^{r \land p} \sum\limits_{\substack{\bm{i} \in [p]^{l_1}, \bm{k} \in [n]^{l_1} \\ \bm{j} \in [p]^{l_2}, \bm{m} \in [n]^{l_2} \\ \#\bm{\{i\}} \cup \bm{\{k\}} \cup \bm{\{j\}} \cup \bm{\{m\}} = r \\ \#\bm{\{i\}} \cup \bm{\{j\}} = b}} \big( \underbrace{\cW_{\bm{i},\bm{k},\bm{j},\bm{m}} - \cW_{\bm{i},\bm{k}} \, \cW_{\bm{j},\bm{m}}}_{=: W(G_{\langle \bm{i},\bm{k} \rangle}, G_{\langle \bm{j},\bm{m} \rangle})} \big) \nonumber\\
		& = \sum\limits_{r=1}^{(2l_1+2l_2) \land n} \sum\limits_{b=1}^{r \land p} \frac{{p \choose b} {n-b \choose r-b}}{n^{l_1+l_2}} \sum\limits_{\substack{\bm{i} \in [b]^{l_1}, \bm{k} \in [r]^{l_1} \\ \bm{j} \in [b]^{l_2}, \bm{m} \in [r]^{l_2} \\ \#\bm{\{i\}} \cup \bm{\{k\}} \cup \bm{\{j\}} \cup \bm{\{m\}} = r \\ \#\bm{\{i\}} \cup \bm{\{j\}} = b}} W(G_{\langle \bm{i},\bm{k} \rangle}, G_{\langle \bm{j},\bm{m} \rangle}) \ .
	\end{align}
	The same argument, which yielded (\ref{Eq_TraceRepresentation4}), now yields
	\begin{align*}
		& \operatorname{Cov}[\tr(\bm{S}_{p,n}^{l_1}), \tr(\bm{S}_{p,n}^{l_2})] = \sum\limits_{r=1}^{(l_1+l_2+2) \land n} \sum\limits_{b=1}^{r \land p} \frac{{p \choose b} {n-b \choose r-b}}{n^{l_1+l_2}} \sum\limits_{\substack{\bm{i} \in [b]^{l_1}, \bm{k} \in [r]^{l_1} \\ \bm{j} \in [b]^{l_2}, \bm{m} \in [r]^{l_2} \\ \#\bm{\{i\}} \cup \bm{\{k\}} \cup \bm{\{j\}} \cup \bm{\{m\}} = r \\ \#\bm{\{i\}} \cup \bm{\{j\}} = b}} W(G_{\langle \bm{i},\bm{k} \rangle}, G_{\langle \bm{j},\bm{m} \rangle})\\
		& = \sum\limits_{b=1}^{r \land p} \Bigg[ \sum\limits_{r=l_1+l_2}^{l_1+l_2+2} \frac{{p \choose b} {n-b \choose r-b}}{n^{l_1+l_2}} \hspace{-0.3cm} \sum\limits_{\substack{\bm{i} \in [b]^{l_1}, \bm{k} \in [r]^{l_1} \\ \bm{j} \in [b]^{l_2}, \bm{m} \in [r]^{l_2} \\ \#\bm{\{i\}} \cup \bm{\{k\}} \cup \bm{\{j\}} \cup \bm{\{m\}} = r \\ \#\bm{\{i\}} \cup \bm{\{j\}} = b}} \hspace{-0.3cm} W(G_{\langle \bm{i},\bm{k} \rangle}, G_{\langle \bm{j},\bm{m} \rangle}) + \mathcal{O}\Big( \frac{p^b}{n^{b+1}} \Big) \Bigg] \ .
	\end{align*}
	In Proposition \ref{WeightCov} we had seen that for $r \in \{l_1+l_2+1,l_1+l_2+2\}$ the weight $W(G_{\langle \bm{i},\bm{k} \rangle}, G_{\langle \bm{j},\bm{m} \rangle})$ will always be zero, meaning we are left with
	\begin{align}\label{Eq_CovSum3}
		& \operatorname{Cov}[\tr(\bm{S}_{p,n}^{l_1}), \tr(\bm{S}_{p,n}^{l_2})] = \sum\limits_{b=1}^{r \land p} \Bigg[ \frac{{p \choose b} {n-b \choose l_1+l_2-b}}{n^{l_1+l_2}} \hspace{-0.3cm} \sum\limits_{\substack{\bm{i} \in [b]^{l_1}, \bm{k} \in [l_1+l_2]^{l_1} \\ \bm{j} \in [b]^{l_2}, \bm{m} \in [l_1+l_2]^{l_2} \\ \#\bm{\{i\}} \cup \bm{\{k\}} \cup \bm{\{j\}} \cup \bm{\{m\}} = r \\ \#\bm{\{i\}} \cup \bm{\{j\}} = b}} \hspace{-0.3cm} W(G_{\langle \bm{i},\bm{k} \rangle}, G_{\langle \bm{j},\bm{m} \rangle}) + \mathcal{O}\Big( \frac{p^b}{n^{b+1}} \Big) \Bigg] \ .
	\end{align}
	It remains to show
	\begin{align*}
		& \sum\limits_{\substack{\bm{i} \in [b]^{l_1} , \, \bm{k} \in [l_1+l_2]^{l_1} \\ \bm{j} \in [b]^{l_2} , \, \bm{m} \in [l_1+l_2]^{l_2} \\ ( \bm{\{i\}} \cup \bm{\{k\}} \cup \bm{\{j\}} \cup \bm{\{m\}}) = [l_1+l_2] \\ \bm{\{i\}} \cup \bm{\{j\}} = [b]}} \hspace{-0.5cm} \big( \cW_{\bm{i},\bm{k},\bm{j},\bm{m}} - \cW_{\bm{i},\bm{k}} \, \cW_{\bm{j},\bm{m}} \big) = C_1(l_1,l_2,b) + (\E[X_{1 \, 1}^4]-3) C_2(l_1,l_2,b)
	\end{align*}
	for all $b \leq l_1+l_2$. We see by Proposition \ref{WeightCov} that
	\begin{align*}
		& \sum\limits_{\substack{\bm{i} \in [b]^{l_1} , \, \bm{k} \in [l_1+l_2]^{l_1} \\ \bm{j} \in [b]^{l_2} , \, \bm{m} \in [l_1+l_2]^{l_2} \\ ( \bm{\{i\}} \cup \bm{\{k\}} \cup \bm{\{j\}} \cup \bm{\{m\}}) = [l_1+l_2] \\ \bm{\{i\}} \cup \bm{\{j\}} = [b]}} \hspace{-0.5cm} \big( \cW_{\bm{i},\bm{k},\bm{j},\bm{m}} - \cW_{\bm{i},\bm{k}} \, \cW_{\bm{j},\bm{m}} \big) = \hspace{-0.5cm} \sum\limits_{\substack{(G_{\left< \bm{i},\bm{k} \right>},G_{\left< \bm{j},\bm{m} \right>}) \in \mathcal{C}^2_{[l_1+l_2],[l_1],[l_2]} \\ B((G_{\left< \bm{i},\bm{k} \right>},G_{\left< \bm{j},\bm{m} \right>})) = [b] }} \hspace{-0.5cm} \big( \cW_{\bm{i},\bm{k},\bm{j},\bm{m}} - \cW_{\bm{i},\bm{k}} \, \cW_{\bm{j},\bm{m}} \big)\\
		& = (\E[X_{1 \, 1}^4]-1) \#\{(G_1,G_2) \in \mathcal{C}^2_{[l_1+l_2],[l_1],[l_2]} \mid B((G_1,G_2)) = [b] ,\\
		& \hspace{4cm} (\operatorname{S}(G_1),\operatorname{S}(G_2)) \in \operatorname{Double-1-d-Ring}_{2}\}\\
		& \hspace{2cm} + \#\{(G_1,G_2) \in \mathcal{C}^2_{[l_1+l_2],[l_1],[l_2]} \mid B((G_1,G_2)) = [b] ,\\
		& \hspace{4cm} (\operatorname{S}(G_1),\operatorname{S}(G_2)) \in \operatorname{Double-1-d-Ring}_{l_0} \text{ for even } l_0 \geq 4 \}\\
		& \hspace{2cm} +\#\{(G_1,G_2) \in \mathcal{C}^2_{[l_1+l_2],[l_1],[l_2]} \mid B((G_1,G_2)) = [b] ,\\
		& \hspace{4cm} (\operatorname{S}(G_1),\operatorname{S}(G_2)) \in \operatorname{Double-2-d-Ring}_{l_0} \text{ for even } l_0 \geq 4 \} \ .
	\end{align*}
	Next sum over all choices of $b'_1,b'_2,w'_1,w'_2$ from Lemma \ref{CountingSproutingDoubleRings}.
	\begin{align*}
		& \sum\limits_{\substack{\bm{i} \in [b]^{l_1} , \, \bm{k} \in [l_1+l_2]^{l_1} \\ \bm{j} \in [b]^{l_2} , \, \bm{m} \in [l_1+l_2]^{l_2} \\ ( \bm{\{i\}} \cup \bm{\{k\}} \cup \bm{\{j\}} \cup \bm{\{m\}}) = [l_1+l_2] \\ \bm{\{i\}} \cup \bm{\{j\}} = [b]}} \hspace{-0.5cm} \big( \cW_{\bm{i},\bm{k},\bm{j},\bm{m}} - \cW_{\bm{i},\bm{k}} \, \cW_{\bm{j},\bm{m}} \big)\\
		& = (\E[X_{1 \, 1}^4]-1) \sum\limits_{b'_1=0 \lor (b-l_2)}^{(l_1-1) \land (b-1)} \#\{(G_1,G_2) \in \mathcal{C}^2_{[l_1+l_2],[l_1],[l_2]} \mid B((G_1,G_2)) = [b] ,\\
		& \hspace{4cm} \text{(\ref{Eq_CountingSproutingDoubleRings1}) holds}, \, (\operatorname{S}(G_1),\operatorname{S}(G_2)) \in \operatorname{Double-1-d-Ring}_{2}\}\\
		& \hspace{1cm} + \sum\limits_{\substack{ l_0 \text{ even} \\ 4 \leq l_0 \leq \underbrace{2l_1 \land 2l_2}_{=2l_1}}} \sum\limits_{b'_1=0 \lor (b-l_2)}^{(l_1-\frac{l_0}{2}) \land (b-\frac{l_0}{2})} \#\{(G_1,G_2) \in \mathcal{C}^2_{[l_1+l_2],[l_1],[l_2]} \mid B((G_1,G_2)) = [b] ,\\
		& \hspace{4cm} \text{(\ref{Eq_CountingSproutingDoubleRings1}) holds}, \, (\operatorname{S}(G_1),\operatorname{S}(G_2)) \in \operatorname{Double-1-d-Ring}_{l_0}\}\\
		& \hspace{1cm} + \sum\limits_{\substack{ l_0 \text{ even} \\ 4 \leq l_0 \leq \underbrace{2l_1 \land 2l_2}_{=2l_1}}} \sum\limits_{b'_1=0 \lor (b-l_2)}^{(l_1-\frac{l_0}{2}) \land (b-\frac{l_0}{2})} \#\{(G_1,G_2) \in \mathcal{C}^2_{[l_1+l_2],[l_1],[l_2]} \mid B((G_1,G_2)) = [b] ,\\
		& \hspace{4cm} \text{(\ref{Eq_CountingSproutingDoubleRings1}) holds}, \, (\operatorname{S}(G_1),\operatorname{S}(G_2)) \in \operatorname{Double-2-d-Ring}_{l_0}\}\\
		& \overset{\text{\ref{CountingSproutingDoubleRings}}}{=} (\E[X_{1 \, 1}^4]-1) \, b! \, w! \sum\limits_{b'_1=0 \lor (b-l_2)}^{(l_1 \land b) - 1} {l_1 \choose b'_1} {l_1 \choose w'_1} {l_2 \choose b'_2} {l_2 \choose w'_2}\\
		& \hspace{2cm} + b! \, w! \sum\limits_{\substack{ l_0 \text{ even} \\ 4 \leq l_0 \leq 2l_1}} \sum\limits_{b'_1=0 \lor (b-l_2)}^{(l_1\land b)-\frac{l_0}{2}} \frac{l_0}{2} \, {l_1 \choose b'_1} {l_1 \choose w'_1} {l_2 \choose b'_2} {l_2 \choose w'_2}\\
		& \hspace{2cm} + b! \, w! \sum\limits_{\substack{ l_0 \text{ even} \\ 4 \leq l_0 \leq 2l_1}} \sum\limits_{b'_1=0 \lor (b-l_2)}^{(l_1\land b)-\frac{l_0}{2}} \frac{l_0}{2} \, {l_1 \choose b'_1} {l_1 \choose w'_1} {l_2 \choose b'_2} {l_2 \choose w'_2} \ .
	\end{align*}
	The three equalities
	\begin{align}\label{Eq_BlackWhiteEqualities1}
		l_1 = \frac{l_0}{2} + b'_1 + w'_1 \ \ , \ \ l_2 = \frac{l_0}{2} + b'_2 + w'_2 \ \ , \ \ b = \frac{l_0}{2} + b'_1 + b'_2 \ ,
	\end{align}
	where the left hand sides and $l_0$ are fixed, show that $b'_2,w'_1,w'_2$ are already uniquely determined by $b'_1$. More precisely we have
	\begin{align}\label{Eq_BlackWhiteEqualities2}
		& w'_1 = l_1-\frac{l_0}{2}-b'_1 \ \ , \ \ w'_2 = l_2 - b + b'_1 \ \ , \ \ b'_2 = b-\frac{l_0}{2}-b'_1 \ .
	\end{align}
	In addition to plugging in these equalities we also add zero twice in the form of
	$$
	\pm b! \, w! \sum\limits_{b'_1=0 \lor (b-l_2)}^{(l_1 \land b) - 1} {l_1 \choose b'_1} {l_1 \choose w'_1} {l_2 \choose b'_2} {l_2 \choose w'_2} \ .
	$$
	As in the proof of Theorem \ref{Thm_MeanExpansion} this heuristically counteracts the exclusion of the case $l_0=2$ from the other two sums. The above formula becomes
	\begin{align*}
		& \sum\limits_{\substack{\bm{i} \in [b]^{l_1} , \, \bm{k} \in [l_1+l_2]^{l_1} \\ \bm{j} \in [b]^{l_2} , \, \bm{m} \in [l_1+l_2]^{l_2} \\ ( \bm{\{i\}} \cup \bm{\{k\}} \cup \bm{\{j\}} \cup \bm{\{m\}}) = [l_1+l_2] \\ \bm{\{i\}} \cup \bm{\{j\}} = [b]}} \hspace{-0.5cm} \big( \cW_{\bm{i},\bm{k},\bm{j},\bm{m}} - \cW_{\bm{i},\bm{k}} \, \cW_{\bm{j},\bm{m}} \big)\\
		& = (\E[X_{1 \, 1}^4]-3) \, b! \, w! \sum\limits_{b'_1=0 \lor (b-l_2)}^{(l_1 \land b) - 1} {l_1 \choose b'_1} {l_1 \choose l_1-\frac{l_0}{2}-b'_1} {l_2 \choose b-\frac{l_0}{2}-b'_1} {l_2 \choose l_2 - b + b'_1}\\
		& \hspace{2cm} + 2 \, b! \, w! \sum\limits_{\substack{ l_0 \text{ even} \\ l_0 \leq 2l_1}} \sum\limits_{b'_1=0 \lor (b-l_2)}^{(l_1\land b)-\frac{l_0}{2}} \frac{l_0}{2} \, {l_1 \choose b'_1} {l_1 \choose l_1-\frac{l_0}{2}-b'_1} {l_2 \choose b-\frac{l_0}{2}-b'_1} {l_2 \choose l_2 - b + b'_1}\\
		& = (\E[X_{1 \, 1}^4]-3) \, b! \, w! \sum\limits_{b'_1=0 \lor (b-l_2)}^{(l_1 \land b) - 1} {l_1 \choose b'_1} {l_1 \choose l_1-1-b'_1} {l_2 \choose b-1-b'_1} {l_2 \choose l_2 - b + b'_1}\\
		& \hspace{2cm} + 2 \, b! \, w! \sum\limits_{m=1}^{l_1} \sum\limits_{b'_1=0 \lor (b-l_2)}^{(l_1\land b)-m} m \, {l_1 \choose b'_1} {l_1 \choose l_1-m-b'_1} {l_2 \choose b-m-b'_1} {l_2 \choose l_2 - b + b'_1} \ ,
	\end{align*}
	where $w = l_1+l_2-b$ denotes the total number of white vertices. Using the fact that the generalized binomial coefficient ${n \choose k}$ is zero when $k>n$, this expression must be equal to
	\begin{align*}
		& (\E[X_{1 \, 1}^4]-3) \, b! \, w! \sum\limits_{k=0}^{b - 1} {l_1 \choose k} {l_1 \choose k+1} {l_2 \choose b-1-k} {l_2 \choose b - k}\\
		& \hspace{2cm} + 2 \, b! \, w! \sum\limits_{k=0}^{b} {l_1 \choose k} {l_2 \choose b - k} \sum\limits_{m=0}^{b-k} m {l_1 \choose k+m} {l_2 \choose b-m-k} \ .
	\end{align*}
\end{proof}

\newpage
\appendix
\section{Adaptability of the methods}\label{Section_Adaptability}
With some minor tweaks to the proofs of Theorems \ref{Thm_MeanExpansion} and \ref{Thm_CovExpansion} one can deduce the following generalizations.
\begin{theorem}[Complex entries]\label{ThmComplexEntries}\
	\\
	Let $(X_{i \, j})_{i,j \in \N}$ be iid complex random variables with $\E[X_{1 \, 1}] = 0$ and $\E[|X_{1 \, 1}|^2] = 1$. Let $\bm{S}_{n,p} = \frac{1}{n} \bm{X}_{p,n} \bm{X}_{p,n}^*$, where $\bm{X}_{p,n} = (X_{i \, j})_{i\leq p, j \leq n}$. For any $l \in \N$ assume $\E[|X_{1,1}|^{2l}]< \infty$, then
	\begin{align}\label{Eq_ComplexResult1}
		\E[\tr(\bm{S}_{p,n}^l)]& \nonumber\\
		= \sum\limits_{b=1}^{l\land p} \Bigg[ & \frac{{p \choose b}{n-b \choose l+1-b}}{n^l} \, l! {l-1 \choose b-1} \nonumber\\
		& \hspace{0cm} + \frac{{p \choose b}{n-b \choose l-b}}{n^l} \Big\{ \tilde{A}_1(l,b) + \big(\E[|X_{1 \, 1}|^4]-2-|\E[X_{1 \, 1}^2]|^2\big) \, A_2(l,b) \Big\} +  \mathcal{O}\left(\frac{p^b}{n^{b+1}}\right) \Bigg] \ ,
	\end{align}
	where $A_2(l,b) := b! \, (l-b)! {l \choose b-1}{l \choose b+1}$ remains unchanged with regards to Theorem \ref{Thm_MeanExpansion} and
	\begin{align*}
		& \tilde{A}_1(l,b) := b! \, (l-b)! \sum\limits_{m=1}^{b} (2m-1) {l \choose b-m}{l \choose b+m-1}\\
		& \hspace{2cm} + b! \, (l-b)! \sum\limits_{m=1}^{b} \big(2m+|\E[X_{1 \, 1}^2]|^{2m}\big) {l \choose b-m}{l \choose b+m} \ .
	\end{align*}
	Further, for any $l_1,l_2 \in \N$ assume $\E[|X_{1 \, 1}|^{2l_1+2l_2}] < \infty$, then
	\begin{align}\label{Eq_ComplexResult2}
		& \operatorname{Cov}[\tr(\bm{S}_{p,n}^{l_1}), \tr(\bm{S}_{p,n}^{l_2})] \nonumber\\
		& = \sum\limits_{b=1}^{l_1+l_2 \land p} \Bigg[ \frac{{p \choose b} {n-b \choose l_1+l_2-b}}{n^{l_1+l_2}} \Big\{ \tilde{C}_1(l_1,l_2,b) + (\E[|X_{1 \, 1}|^4]-2-|\E[X_{1 \, 1}^2]|^{2}) C_2(l_1,l_2,b) \Big\} + \mathcal{O}\left(\frac{p^{b}}{n^{b+1}}\right) \Bigg] \ ,
	\end{align}
	where
	\begin{align*}
		& \tilde{C}_1(l_1,l_2,b) := b! \, (l_1+l_2-b)! \sum\limits_{k=0}^{b} {l_1 \choose k} {l_2 \choose b - k} \sum\limits_{m=0}^{b-k} m \big( 1 + |\E[X_{1 \, 1}^2]|^{4m} \big) {l_1 \choose k+m} {l_2 \choose b-m-k}
	\end{align*}
	and $C_2(l_1,l_2,b)$ remains unchanged with regards to Theorem \ref{Thm_CovExpansion}.
\end{theorem}
\begin{proof}\
	\\
	We first prove (\ref{Eq_ComplexResult1}) analogously to Theorem \ref{Thm_MeanExpansion}. The new weight is
	\begin{align}\label{Eq_CompleyWeight1}
		\tilde{W}_{\bm{i},\bm{j}} & = \E\big[(X_{i_1 \, j_1} \ol{X_{i_2 \, j_1}}) \, (X_{i_2 \, j_2} \ol{X_{i_3 \, j_2}}) \, \dots \, (X_{i_l \, j_l} \ol{X_{i_1\, j_l}}) \big] \nonumber\\
		& = \prod\limits_{(s,t) \in p \times n} \E\Big[X_{s,t}^{e_{s,t}(G_{\langle \bm{i},\bm{j} \rangle})} \ol{X_{s,t}}^{\ol{e}_{s,t}(G_{\langle \bm{i},\bm{j} \rangle})}\Big] \ ,
	\end{align}
	where $e_{s,t}(G_{\langle \bm{i},\bm{j} \rangle})$ denotes the number of edges in $G_{\langle \bm{i},\bm{j} \rangle}$ from $s$ to $t$, which do not change their direction in the reversed graph $\operatorname{R}(G_{\langle \bm{i},\bm{j} \rangle})$, and $\ol{e}_{s,t}(G_{\langle \bm{i},\bm{j} \rangle})$ describes the number of edges from $s$ to $t$, which do change direction. Most of the arguments used in the proof of Theorem \ref{Thm_MeanExpansion} (see Section \ref{MainResult1}) remain unchanged, only the weights in Proposition \ref{WeightOfFusedTrees} change, while thankfully no new cases of seed graphs with non-zero weight arise. The new weights can be seen to be
	\begin{align*}
		& \cW(G_{\left< \bm{i}, \bm{k} \right>}) = \begin{cases}
			\E[|X_{1 \, 1}|^4] & \text{, if $\operatorname{S}(G_{\left< \bm{i}, \bm{k} \right>}) \in \operatorname{2-d-Ring}_2$}\\
			1 & \text{, if $\operatorname{S}(G_{\left< \bm{i}, \bm{k} \right>}) \in \operatorname{2-d-Ring}_{l_0}$ for some $l_0 \neq 2$}\\
			\big|\E[X_{1 \, 1}^2]\big|^{l_0} & \text{, if $\operatorname{S}(G_{\left< \bm{i}, \bm{k} \right>}) \in \operatorname{1-d-Ring}_{l_0}$ for even $l_0 \geq 4$}\\
			0 & \text{ else}
		\end{cases} \ .
	\end{align*}
	This changes the evaluation of $\sum\limits_{\substack{G_{\left< \bm{i},\bm{k} \right>} \in \mathcal{C}_{[l],2l} \\ B(G_{\left< \bm{i},\bm{k} \right>}) = [b]}} \cW(G_{\left< \bm{i},\bm{k} \right>})$ from the proof of Theorem \ref{Thm_MeanExpansion} to
	\begin{align*}
		\sum\limits_{\substack{G_{\left< \bm{i},\bm{k} \right>} \in \mathcal{C}_{[l],2l} \\ B(G_{\left< \bm{i},\bm{k} \right>}) = [b]}} \cW(G_{\left< \bm{i},\bm{k} \right>}) & = \sum\limits_{\substack{m=1 \\ m \leq w+1}}^{b} \#\{G \in \mathcal{C}_{[l],2l} \mid \operatorname{S}(G) \in \operatorname{2-d-Ring}_{2m-1}, \, B(G) = [b] \}\\
		& \ \ \ \ + \sum\limits_{\substack{m=2 \\ m \leq w}}^{b} \#\{G \in \mathcal{C}_{[l],2l} \mid \operatorname{S}(G) \in \operatorname{2-d-Ring}_{2m}, \, B(G) = [b] \}\\
		& \ \ \ \ + 2 \, \#\{G \in \mathcal{C}_{[l],2l} \mid \operatorname{S}(G) \in \operatorname{2-d-Ring}_{2} , \, B(G) = [b] \}\\
		& \ \ \ \ + \sum\limits_{\substack{m=2 \\ m \leq w}}^{b} |\E[X_{1 \, 1}^2]|^{2m} \#\{G \in \mathcal{C}_{[l],2l} \mid \operatorname{S}(G) \in \operatorname{1-d-Ring}_{2m} , \, B(G) = [b] \}\\
		& \ \ \ \ + |\E[X_{1 \, 1}^2]|^2 \#\{G \in \mathcal{C}_{[l],2l} \mid \operatorname{S}(G) \in \operatorname{2-d-Ring}_{2} , \, B(G) = [b] \}\\
		& \ \ \ \ + (\E[|X_{1 \, 1}|^4]-2-|\E[X_{1 \, 1}^2]|^2)\\
		& \hspace{3cm} \times \#\{G \in \mathcal{C}_{[l],2l} \mid \operatorname{S}(G) \in \operatorname{2-d-Ring}_{2} , \, B(G) = [b] \} \ ,
	\end{align*}
	which again with Corollaries \ref{SproutingCorollary1} and \ref{SproutingCorollary2} yields
	\begin{align*}
		\sum\limits_{\substack{G_{\left< \bm{i},\bm{k} \right>} \in \mathcal{C}_{[l],2l} \\ B(G_{\left< \bm{i},\bm{k} \right>}) = [b]}} \cW(G_{\left< \bm{i},\bm{k} \right>}) & = b! \, w! \sum\limits_{m=1}^{b \land (l-b+1)} \big(2m-1\big) {l \choose b-m}{l \choose b+m-1}\\
		& \ \ \ \ + b! \, w! \sum\limits_{m=1}^{b \land (l-b)} \big(2m+|\E[X_{1 \, 1}^2]|^{2m}\big) {l \choose b-m}{l \choose b+m}\\
		& \ \ \ \ + \big(\E[|X_{1 \, 1}|^4]-2-|\E[X_{1 \, 1}^2]|^2\big) \, b! \, w! {l \choose b-1}{l \choose b+1} \ .
	\end{align*}
	We have thus proven (\ref{Eq_ComplexResult1}) and employ similar changes to the proof of Theorem \ref{Thm_CovExpansion} (see Section \ref{MainResult2}) to show (\ref{Eq_ComplexResult2}).\\
	\\
	Similarly to (\ref{Eq_CompleyWeight1}) the new weight $\tilde{\cW}_{\bm{i},\bm{k},\bm{j},\bm{m}}$ is given by
	\begin{align*}
		& \tilde{\cW}_{\bm{i},\bm{k},\bm{j},\bm{m}} = \prod\limits_{(s,t) \in p \times n} \E\Big[X_{s,t}^{e_{s,t}(G_{\langle \bm{i},\bm{k} \rangle}) + e_{s,t}(G_{\langle \bm{j},\bm{m} \rangle})} \ol{X_{s,t}}^{\ol{e}_{s,t}(G_{\langle \bm{i},\bm{k} \rangle}) + \ol{e}_{s,t}(G_{\langle \bm{j},\bm{m} \rangle})}\Big]
	\end{align*}
	and the result of Proposition \ref{WeightCov} changes to
	\begin{align*}
		& \cW_{\bm{i},\bm{k},\bm{j},\bm{m}} - \cW_{\bm{i},\bm{k}} \cW_{\bm{j},\bm{m}}\\
		& = \begin{cases}
			0 & \text{, if } \#( \bm{\{i\}} \cup \bm{\{k\}} \cup \bm{\{j\}} \cup \bm{\{m\}}) > l\\
			\E[|X_{1 \, 1}|^4]-1 & \text{, if $\operatorname{S}((G_{\left< \bm{i},\bm{k} \right>},G_{\left< \bm{j},\bm{m} \right>})) \in \operatorname{Double-2-d-Ring}_{2}$}\\
			|\E[X_{1 \, 1}^2]|^{2l_0} & \text{, if $\operatorname{S}((G_{\left< \bm{i},\bm{k} \right>},G_{\left< \bm{j},\bm{m} \right>})) \in \operatorname{Double-1-d-Ring}_{l_0}$ for even $l_0 \geq 4$}\\
			1 & \text{, if $\operatorname{S}((G_{\left< \bm{i},\bm{k} \right>},G_{\left< \bm{j},\bm{m} \right>})) \in \operatorname{Double-2-d-Ring}_{l_0}$ for even $l_0 \geq 4$}\\
			0 & \text{else}
		\end{cases} \ .
	\end{align*}
	As in the proof of Theorem \ref{Thm_CovExpansion} this yields
	\begin{align*}
		& \sum\limits_{\substack{\bm{i} \in [b]^{l_1} , \, \bm{k} \in [l_1+l_2]^{l_1} \\ \bm{j} \in [b]^{l_2} , \, \bm{m} \in [l_1+l_2]^{l_2} \\ ( \bm{\{i\}} \cup \bm{\{k\}} \cup \bm{\{j\}} \cup \bm{\{m\}}) = [l_1+l_2] \\ \bm{\{i\}} \cup \bm{\{j\}} = [b]}} \hspace{-0.5cm} \big( \tilde{\cW}_{\bm{i},\bm{k},\bm{j},\bm{m}} - \tilde{\cW}_{\bm{i},\bm{k}} \, \tilde{\cW}_{\bm{j},\bm{m}} \big)\\
		& = (\E[|X_{1 \, 1}|^4]-1) \#\{(G_1,G_2) \in \mathcal{C}^2_{[l_1+l_2],[l_1],[l_2]} \mid B((G_1,G_2)) = [b] ,\\
		& \hspace{4cm} (\operatorname{S}(G_1),\operatorname{S}(G_2)) \in \operatorname{Double-1-d-Ring}_{2}\}\\
		& \hspace{2cm} + \sum\limits_{\substack{l_0 \text{ even} \\ 4\leq l_0 \leq 2l_1}} |\E[X_{1 \, 1}^2]|^{2l_0} \#\{(G_1,G_2) \in \mathcal{C}^2_{[l_1+l_2],[l_1],[l_2]} \mid B((G_1,G_2)) = [b] ,\\
		& \hspace{6cm} (\operatorname{S}(G_1),\operatorname{S}(G_2)) \in \operatorname{Double-1-d-Ring}_{l_0} \}\\
		& \hspace{2cm} +\#\{(G_1,G_2) \in \mathcal{C}^2_{[l_1+l_2],[l_1],[l_2]} \mid B((G_1,G_2)) = [b] ,\\
		& \hspace{4cm} (\operatorname{S}(G_1),\operatorname{S}(G_2)) \in \operatorname{Double-2-d-Ring}_{l_0} \text{ for even } l_0 \geq 4 \}\\
		& = (\E[|X_{1 \, 1}|^4]-1) \, b! \, w! \sum\limits_{b'_1=0 \lor (b-l2)}^{(l_1 \land b) - 1} {l_1 \choose b'_1} {l_1 \choose w'_1} {l_2 \choose b'_2} {l_2 \choose w'_2}\\
		& \hspace{2cm} + b! \, w! \sum\limits_{\substack{ l_0 \text{ even} \\ 4 \leq l_0 \leq 2l_1}} |\E[X_{1 \, 1}^2]|^{2l_0} \sum\limits_{b'_1=0 \lor (b-l2)}^{(l_1\land b)-\frac{l_0}{2}} \frac{l_0}{2} \, {l_1 \choose b'_1} {l_1 \choose w'_1} {l_2 \choose b'_2} {l_2 \choose w'_2}\\
		& \hspace{2cm} + b! \, w! \sum\limits_{\substack{ l_0 \text{ even} \\ 4 \leq l_0 \leq 2l_1}} \sum\limits_{b'_1=0 \lor (b-l2)}^{(l_1\land b)-\frac{l_0}{2}} \frac{l_0}{2} \, {l_1 \choose b'_1} {l_1 \choose w'_1} {l_2 \choose b'_2} {l_2 \choose w'_2}\\
		& = (\E[|X_{1 \, 1}|^4]-2-|\E[X_{1 \, 1}^2]|^{2}) \, \overbrace{b! \, w! \sum\limits_{k=0}^{b - 1} {l_1 \choose k} {l_1 \choose k+1} {l_2 \choose b-1-k} {l_2 \choose b - k}}^{= C_2(l_1l_2,b)}\\
		& \hspace{2cm} + \underbrace{b! \, w! \sum\limits_{k=0}^{b} {l_1 \choose k} {l_2 \choose b - k} \sum\limits_{m=0}^{b-k} m \big( 1 + |\E[X_{1 \, 1}^2]|^{4m} \big) {l_1 \choose k+m} {l_2 \choose b-m-k}}_{= \tilde{C}_1(l_1,l_2,b)} \ .
	\end{align*}
\end{proof}

\begin{theorem}[Covariance with re-sampled rows]\label{Thm_ReSample}\
	\\
	For $p<n$ let $\bm{X} = (X_{i \, j})_{i \leq p, j \leq n}$ be a $(p \times n)$ matrix with iid centered entries with $\E[X_{i\,j}^2]=1$. Let $\tilde{\bm{X}}$ have the same distribution as $\bm{X}$ with $\tX_{i\,j} = X_{i\,j}$ for $i\leq s$ and $(\tX_{i\,j})_{s<i\leq p, j \leq n}$ independent of $(X_{i\,j})_{i,j \in \N}$. Heuristically, $\tilde{\bm{X}}$ is the result of re-sampling the first $s$ rows of $\bm{X}$. For any $l_1,l_2 \in \N$ assume $\E[X_{1 \, 1}^{2l_1+2l_2}] < \infty$, then for the sample covariance matrices $\bm{S}_{p,n} = \frac{1}{n} \bm{X} \bm{X}^T$ and $\tilde{\bm{S}}_{p,n} = \frac{1}{n} \tilde{\bm{X}} \tilde{\bm{X}}^T$ we have
	\begin{align*}
		& \operatorname{Cov}[\tr(\bm{S}_{p,n}^{l_1}), \tr(\tilde{\bm{S}}_{p,n}^{l_2})]\\
		& = \sum\limits_{b=1}^{(l_1+l_2) \land p} \Bigg[ \sum\limits_{t=1}^{b \land s} \frac{{s \choose t} {p-t \choose b-t} {n-b \choose l_1+l_2-b}}{{b \choose t} n^{l_1+l_2}} \tilde{C}_1(l_1,l_2,b,t)\\
		& \hspace{2cm} + \frac{s {p-1 \choose b-1} {n-b \choose l_1+l_2-b}}{b n^{l_1+l_2}} (\E[X_{1 \, 1}^4]-3) C_2(l_1,l_2,b) + \mathcal{O}\Big( \frac{p^b}{n^{b+1}} \Big) \Bigg] \ ,
	\end{align*}
	where
	\begin{align*}
		& \tilde{C}_1(l_1,l_2,b,t) := 2 b! (l_1+l_2-b)! \sum\limits_{k=0}^{b-t} t \, {l_1 \choose k} {l_1 \choose k+t} {l_2 \choose b-t-k} {l_2 \choose b-k}
	\end{align*}
	and $C_2(l_1,l_2,b)$ is as in the formulation of Theorem \ref{Thm_CovExpansion}.
\end{theorem}
\begin{proof}\
	\\
	By symmetry we can without loss of generality assume $l_1 \leq l_2$. Similar to the proof of Theorem \ref{Thm_CovExpansion} we have
	\begin{align}\label{Eq_ReSample1}
		& \operatorname{Cov}[\tr(\bm{S}_{p,n}^{l_1}), \tr(\tilde{\bm{S}}_{p,n}^{l_2})] \nonumber\\
		& = \sum\limits_{b=1}^{r \land p} \Bigg[ \sum\limits_{r=l_1+l_2}^{l_1+l_2+2} \frac{{p \choose b} {n-b \choose r-b}}{n^{l_1+l_2}} \hspace{-0.3cm} \sum\limits_{\substack{\bm{i} \in [b]^{l_1}, \bm{k} \in [r]^{l_1} \\ \bm{j} \in [b]^{l_2}, \bm{m} \in [r]^{l_2} \\ \#\bm{\{i\}} \cup \bm{\{k\}} \cup \bm{\{j\}} \cup \bm{\{m\}} = r \\ \#\bm{\{i\}} \cup \bm{\{j\}} = b}} \hspace{-0.3cm} \big( \tilde{\cW}_{\bm{i},\bm{k},\bm{j},\bm{m}} - \cW_{\bm{i},\bm{k}} \cW_{\bm{j},\bm{m}} \big) + \mathcal{O}\Big( \frac{p^b}{n^{b+1}} \Big) \Bigg] \ ,
	\end{align}
	where t singular weights $\cW_{\bm{i},\bm{k}}$ and $\cW_{\bm{j},\bm{m}}$ remain unchanged with regards to Theorem \ref{Thm_CovExpansion}, while the joint weight has the form
	\begin{align*}
		\tilde{\cW}_{\bm{i},\bm{k},\bm{j},\bm{m}} &= \E\left[ (X_{i_1 \, k_1} X_{i_2 \, k_1}) ... (X_{i_{l_1} \, k_{l_1}} X_{i_1 \, k_{l_1}}) \cdot (\tX_{j_1 \, m_1} \tX_{j_2 \, m_1}) ... (\tX_{j_{l_2} \, m_{l_2}} \tX_{j_1 \, m_{l_2}}) \right]\\
		& = \prod\limits_{(s,t)\in [p] \times [n]} \E\bigg[X_{s,t}^{A_{s,t}(\operatorname{R}(G_{\left< \bm{i}, \bm{k} \right>}))} \tX_{s,t}^{A_{s,t}(\operatorname{R}(G_{\left< \bm{j}, \bm{m} \right>}))}\bigg] \ .
	\end{align*}
	The result of Proposition \ref{WeightCov} changes to
	\begin{align*}
		& \tilde{\cW}_{\bm{i},\bm{k},\bm{j},\bm{m}} - \cW_{\bm{i},\bm{k}} \cW_{\bm{j},\bm{m}}\\
		& = \mathbbm{1}_{\bm{\{i\}} \cap V(S(G_{\left< \bm{i},\bm{k} \right>})) \subset [s]} \begin{cases}
			0 & \text{, if } \#( \bm{\{i\}} \cup \bm{\{k\}} \cup \bm{\{j\}} \cup \bm{\{m\}}) > l\\
			\E[X_{1 \, 1}^4]-1 & \text{, if $\operatorname{S}((G_{\left< \bm{i},\bm{k} \right>},G_{\left< \bm{j},\bm{m} \right>})) \in \operatorname{Double-2-d-Ring}_{2}$}\\
			1 & \text{, if $\operatorname{S}((G_{\left< \bm{i},\bm{k} \right>},G_{\left< \bm{j},\bm{m} \right>})) \in \operatorname{Double-1-d-Ring}_{l_0}$ for even $l_0 \geq 4$}\\
			1 & \text{, if $\operatorname{S}((G_{\left< \bm{i},\bm{k} \right>},G_{\left< \bm{j},\bm{m} \right>})) \in \operatorname{Double-2-d-Ring}_{l_0}$ for even $l_0 \geq 4$}\\
			0 & \text{else}
		\end{cases}\\
		& = \mathbbm{1}_{\bm{\{i\}} \cap V(S(G_{\left< \bm{i},\bm{k} \right>})) \subset [s]} \big( \cW_{\bm{i},\bm{k},\bm{j},\bm{m}} - \cW_{\bm{i},\bm{k}} \cW_{\bm{j},\bm{m}} \big) \ ,
	\end{align*}
	since the ring-type seed graphs will only contribute, if all black vertices in the ring structure are from $[s]$, and $V(S(G_{\left< \bm{i},\bm{k} \right>})) = V(S(G_{\left< \bm{j},\bm{m} \right>})) = V\big( S(G_{\left< \bm{i},\bm{k} \right>},G_{\left< \bm{j},\bm{m} \right>}) \big)$ for all seed graphs $S(G_{\left< \bm{i},\bm{k} \right>},G_{\left< \bm{j},\bm{m} \right>})$ with non-zero weight, which follows from Lemma \ref{SproutingDoubleCondition}. Consequently the equality (\ref{Eq_ReSample1}) becomes
	\begin{align*}
		& \operatorname{Cov}[\tr(\bm{S}_{p,n}^{l_1}), \tr(\tilde{\bm{S}}_{p,n}^{l_2})] \nonumber\\
		& = \sum\limits_{b=1}^{(l_1+l_2) \land p} \Bigg[ \frac{{p \choose b} {n-b \choose l_1+l_2-b}}{n^{l_1+l_2}} \hspace{-0.3cm} \sum\limits_{\substack{\bm{i} \in [b]^{l_1}, \bm{k} \in [l_1+l_2]^{l_1} \\ \bm{j} \in [b]^{l_2}, \bm{m} \in [l_1+l_2]^{l_2} \\ \bm{\{i\}} \cup \bm{\{k\}} \cup \bm{\{j\}} \cup \bm{\{m\}} = [l_1+l_2] \\ \bm{\{i\}} \cup \bm{\{j\}} = [b]}} \hspace{-0.3cm} \mathbbm{1}_{[b] \cap V(S(G_{\left< \bm{i},\bm{k} \right>})) \subset [s]} \big( \cW_{\bm{i},\bm{k},\bm{j},\bm{m}} - \cW_{\bm{i},\bm{k}} \cW_{\bm{j},\bm{m}} \big) + \mathcal{O}\Big( \frac{p^b}{n^{b+1}} \Big) \Bigg]\\
		& = \sum\limits_{b=1}^{(l_1+l_2) \land p} \Bigg[ \frac{{p \choose b} {n-b \choose l_1+l_2-b}}{n^{l_1+l_2}} \hspace{-0.3cm} \sum\limits_{\substack{\bm{i} \in [b]^{l_1}, \bm{k} \in [l_1+l_2]^{l_1} \\ \bm{j} \in [b]^{l_2}, \bm{m} \in [l_1+l_2]^{l_2} \\ \bm{\{i\}} \cup \bm{\{k\}} \cup \bm{\{j\}} \cup \bm{\{m\}} = [l_1+l_2] \\ \bm{\{i\}} \cup \bm{\{j\}} = [b]}} \hspace{-0.3cm} \mathbbm{1}_{\bm{\{i\}} \cap \bm{\{j\}} \subset [s]} \big( \cW_{\bm{i},\bm{k},\bm{j},\bm{m}} - \cW_{\bm{i},\bm{k}} \cW_{\bm{j},\bm{m}} \big) + \mathcal{O}\Big( \frac{p^b}{n^{b+1}} \Big) \Bigg]\\
		& = \sum\limits_{b=1}^{(l_1+l_2) \land p} \Bigg[ \sum\limits_{t=1}^{b \land s} \frac{{s \choose t} {p-t \choose b-t} {n-b \choose l_1+l_2-b}}{n^{l_1+l_2}} \hspace{-0.3cm} \sum\limits_{\substack{\bm{i} \in [b]^{l_1}, \bm{k} \in [l_1+l_2]^{l_1} \\ \bm{j} \in [b]^{l_2}, \bm{m} \in [l_1+l_2]^{l_2} \\ \bm{\{i\}} \cup \bm{\{k\}} \cup \bm{\{j\}} \cup \bm{\{m\}} = [l_1+l_2] \\ \bm{\{i\}} \cup \bm{\{j\}} = [b] \\ \bm{\{i\}} \cap \bm{\{j\}} = [t]}} \hspace{-0.3cm} \big( \cW_{\bm{i},\bm{k},\bm{j},\bm{m}} - \cW_{\bm{i},\bm{k}} \cW_{\bm{j},\bm{m}} \big) + \mathcal{O}\Big( \frac{p^b}{n^{b+1}} \Big) \Bigg] \ .
	\end{align*}
	For easier calculation later on we de-specify the set $\bm{\{i\}} \cap \bm{\{j\}} \subset [b]$, which gives us an additional factor $\frac{1}{{b \choose t}}$, since we were only looking at one of the ${b \choose t}$ many possible subsets of $[b]$ with $t$ elements. We arrive at
	\begin{align*}
		& \operatorname{Cov}[\tr(\bm{S}_{p,n}^{l_1}), \tr(\tilde{\bm{S}}_{p,n}^{l_2})]\\
		& = \sum\limits_{b=1}^{(l_1+l_2) \land p} \Bigg[ \sum\limits_{t=1}^{b \land s} \frac{{s \choose t} {p-t \choose b-t} {n-b \choose l_1+l_2-b}}{{b \choose t} n^{l_1+l_2}} \hspace{-0.3cm} \sum\limits_{\substack{\bm{i} \in [b]^{l_1}, \bm{k} \in [l_1+l_2]^{l_1} \\ \bm{j} \in [b]^{l_2}, \bm{m} \in [l_1+l_2]^{l_2} \\ \bm{\{i\}} \cup \bm{\{k\}} \cup \bm{\{j\}} \cup \bm{\{m\}} = [l_1+l_2] \\ \bm{\{i\}} \cup \bm{\{j\}} = [b] \\ \#\bm{\{i\}} \cap \bm{\{j\}} = t}} \hspace{-0.3cm} \big( \cW_{\bm{i},\bm{k},\bm{j},\bm{m}} - \cW_{\bm{i},\bm{k}} \cW_{\bm{j},\bm{m}} \big) + \mathcal{O}\Big( \frac{p^b}{n^{b+1}} \Big) \Bigg] \ .
	\end{align*}
	As in the proof of Theorem \ref{Thm_CovExpansion} we can use Proposition \ref{WeightCov} to see
	\begin{align*}
		& \sum\limits_{\substack{\bm{i} \in [b]^{l_1} , \, \bm{k} \in [l_1+l_2]^{l_1} \\ \bm{j} \in [b]^{l_2} , \, \bm{m} \in [l_1+l_2]^{l_2} \\ ( \bm{\{i\}} \cup \bm{\{k\}} \cup \bm{\{j\}} \cup \bm{\{m\}}) = [l_1+l_2] \\ \bm{\{i\}} \cup \bm{\{j\}} = [b] \\ \bm{\{i\}} \cap \bm{\{j\}} = [t]}} \hspace{-0.5cm} \big( \cW_{\bm{i},\bm{k},\bm{j},\bm{m}} - \cW_{\bm{i},\bm{k}} \, \cW_{\bm{j},\bm{m}} \big) = \hspace{-0.5cm} \sum\limits_{\substack{(G_{\left< \bm{i},\bm{k} \right>},G_{\left< \bm{j},\bm{m} \right>}) \in \mathcal{C}^2_{[l_1+l_2],[l_1],[l_2]} \\ B((G_{\left< \bm{i},\bm{k} \right>},G_{\left< \bm{j},\bm{m} \right>})) = [b] \\ [b] \cap V(S(G_{\left< \bm{i},\bm{k} \right>})) = [t]}} \hspace{-0.5cm} \big( \cW_{\bm{i},\bm{k},\bm{j},\bm{m}} - \cW_{\bm{i},\bm{k}} \, \cW_{\bm{j},\bm{m}} \big)\\
		& = (\E[X_{1 \, 1}^4]-1) \#\{(G_1,G_2) \in \mathcal{C}^2_{[l_1+l_2],[l_1],[l_2]} \mid B((G_1,G_2)) = [b] , \#[b] \cap V(S(G_{\left< \bm{i},\bm{k} \right>})) = t ,\\
		& \hspace{4cm} (\operatorname{S}(G_1),\operatorname{S}(G_2)) \in \operatorname{Double-1-d-Ring}_{2}\}\\
		& \hspace{2cm} + \#\{(G_1,G_2) \in \mathcal{C}^2_{[l_1+l_2],[l_1],[l_2]} \mid B((G_1,G_2)) = [b] , \#[b] \cap V(S(G_{\left< \bm{i},\bm{k} \right>})) = t ,\\
		& \hspace{4cm} (\operatorname{S}(G_1),\operatorname{S}(G_2)) \in \operatorname{Double-1-d-Ring}_{l_0} \text{ for even } l_0 \geq 4 \}\\
		& \hspace{2cm} +\#\{(G_1,G_2) \in \mathcal{C}^2_{[l_1+l_2],[l_1],[l_2]} \mid B((G_1,G_2)) = [b] , \#[b] \cap V(S(G_{\left< \bm{i},\bm{k} \right>})) = t ,\\
		& \hspace{4cm} (\operatorname{S}(G_1),\operatorname{S}(G_2)) \in \operatorname{Double-2-d-Ring}_{l_0} \text{ for even } l_0 \geq 4 \} \ .
	\end{align*}
	By Lemma \ref{DoubleRingColoring} we know that the number of black vertices in the ring structure of a Double-Ring type graph of even length $l_0$ is $t$, iff $l_0=2t$. If $t=1$, the above formula becomes
	\begin{align*}
		& \sum\limits_{\substack{\bm{i} \in [b]^{l_1} , \, \bm{k} \in [l_1+l_2]^{l_1} \\ \bm{j} \in [b]^{l_2} , \, \bm{m} \in [l_1+l_2]^{l_2} \\ ( \bm{\{i\}} \cup \bm{\{k\}} \cup \bm{\{j\}} \cup \bm{\{m\}}) = [l_1+l_2] \\ \bm{\{i\}} \cup \bm{\{j\}} = [b] \\ \bm{\{i\}} \cap \bm{\{j\}} = [t]}} \hspace{-0.5cm} \big( \cW_{\bm{i},\bm{k},\bm{j},\bm{m}} - \cW_{\bm{i},\bm{k}} \, \cW_{\bm{j},\bm{m}} \big) =\\
		& = (\E[X_{1 \, 1}^4]-1) \#\{(G_1,G_2) \in \mathcal{C}^2_{[l_1+l_2],[l_1],[l_2]} \mid B((G_1,G_2)) = [b] ,\\
		& \hspace{4cm} (\operatorname{S}(G_1),\operatorname{S}(G_2)) \in \operatorname{Double-1-d-Ring}_{2}\}\\
		& = (\E[X_{1 \, 1}^4]-1) \sum\limits_{b_1'=0\lor(b-l_2)}^{(l_1-1)\land(b-1)} \#\{(G_1,G_2) \in \mathcal{C}^2_{[l_1+l_2],[l_1],[l_2]} \mid B((G_1,G_2)) = [b] ,\\
		& \hspace{4cm} \text{(\ref{Eq_CountingSproutingDoubleRings1}) holds}, (\operatorname{S}(G_1),\operatorname{S}(G_2)) \in \operatorname{Double-1-d-Ring}_{2}\}\\
		& \overset{\text{\ref{CountingSproutingDoubleRings}}}{=} (\E[X_{1 \, 1}^4]-1) \sum\limits_{b_1'=0\lor(b-l_2)}^{(l_1-1)\land(b-1)} b! w! {l_1 \choose b_1'} {l_1 \choose w_1'} {l_2 \choose b_2'} {l_2 \choose w_2'}\\
		& \overset{\text{(\ref{Eq_BlackWhiteEqualities2})}}{=} (\E[X_{1 \, 1}^4]-1) b! w! \sum\limits_{b_1'=0}^{b-1} {l_1 \choose b_1'} {l_1 \choose b_1'+1} {l_2 \choose b-1-b_1'} {l_2 \choose b-b_1'}
	\end{align*}
	and for $t > 1$ we analogously have
	\begin{align*}
		& \sum\limits_{\substack{\bm{i} \in [b]^{l_1} , \, \bm{k} \in [l_1+l_2]^{l_1} \\ \bm{j} \in [b]^{l_2} , \, \bm{m} \in [l_1+l_2]^{l_2} \\ ( \bm{\{i\}} \cup \bm{\{k\}} \cup \bm{\{j\}} \cup \bm{\{m\}}) = [l_1+l_2] \\ \bm{\{i\}} \cup \bm{\{j\}} = [b] \\ \bm{\{i\}} \cap \bm{\{j\}} = [t]}} \hspace{-0.5cm} \big( \cW_{\bm{i},\bm{k},\bm{j},\bm{m}} - \cW_{\bm{i},\bm{k}} \, \cW_{\bm{j},\bm{m}} \big)\\
		& = \#\{(G_1,G_2) \in \mathcal{C}^2_{[l_1+l_2],[l_1],[l_2]} \mid B((G_1,G_2)) = [b] ,\\
		& \hspace{4cm} (\operatorname{S}(G_1),\operatorname{S}(G_2)) \in \operatorname{Double-1-d-Ring}_{2t} \}\\
		& \hspace{2cm} +\#\{(G_1,G_2) \in \mathcal{C}^2_{[l_1+l_2],[l_1],[l_2]} \mid B((G_1,G_2)) = [b] ,\\
		& \hspace{4cm} (\operatorname{S}(G_1),\operatorname{S}(G_2)) \in \operatorname{Double-2-d-Ring}_{2t} \}\\
		& = \sum\limits_{b_1'=0\lor(b-l_2)}^{(l_1-t) \land (b-t)} \#\{(G_1,G_2) \in \mathcal{C}^2_{[l_1+l_2],[l_1],[l_2]} \mid B((G_1,G_2)) = [b] ,\\
		& \hspace{4cm} \text{(\ref{Eq_CountingSproutingDoubleRings1}) holds}, (\operatorname{S}(G_1),\operatorname{S}(G_2)) \in \operatorname{Double-1-d-Ring}_{2t} \}\\
		& \hspace{1cm} + \sum\limits_{b_1'=0\lor(b-l_2)}^{(l_1-t) \land (b-t)} \#\{(G_1,G_2) \in \mathcal{C}^2_{[l_1+l_2],[l_1],[l_2]} \mid B((G_1,G_2)) = [b] ,\\
		& \hspace{4cm} \text{(\ref{Eq_CountingSproutingDoubleRings1}) holds}, (\operatorname{S}(G_1),\operatorname{S}(G_2)) \in \operatorname{Double-2-d-Ring}_{2t} \}\\
		& \overset{\text{\ref{CountingSproutingDoubleRings}}}{=} \sum\limits_{b_1'=0\lor(b-l_2)}^{(l_1-t) \land (b-t)} t \, b! \, w! \, {l_1 \choose b'_1} {l_1 \choose w'_1} {l_2 \choose b'_2} {l_2 \choose w'_2}\\
		& \hspace{1cm} + \sum\limits_{b_1'=0\lor(b-l_2)}^{(l_1-t) \land (b-t)} t \, b! \, w! \, {l_1 \choose b'_1} {l_1 \choose w'_1} {l_2 \choose b'_2} {l_2 \choose w'_2}\\
		& \overset{\text{(\ref{Eq_BlackWhiteEqualities2})}}{=} 2 b! w! \sum\limits_{b_1'=0}^{b-t} t \, {l_1 \choose b_1'} {l_1 \choose b_1'+t} {l_2 \choose b-t-b_1'} {l_2 \choose b-b_1'}
	\end{align*}
	Gathering our results we see
	\begin{align*}
		& \operatorname{Cov}[\tr(\bm{S}_{p,n}^{l_1}), \tr(\tilde{\bm{S}}_{p,n}^{l_2})]\\
		& = \sum\limits_{b=1}^{(l_1+l_2) \land p} \Bigg[ \sum\limits_{t=1}^{b \land s} \frac{{s \choose t} {p-t \choose b-t} {n-b \choose l_1+l_2-b}}{{b \choose t} n^{l_1+l_2}} \hspace{-0.3cm} \sum\limits_{\substack{\bm{i} \in [b]^{l_1}, \bm{k} \in [l_1+l_2]^{l_1} \\ \bm{j} \in [b]^{l_2}, \bm{m} \in [l_1+l_2]^{l_2} \\ \bm{\{i\}} \cup \bm{\{k\}} \cup \bm{\{j\}} \cup \bm{\{m\}} = [l_1+l_2] \\ \bm{\{i\}} \cup \bm{\{j\}} = [b] \\ \#\bm{\{i\}} \cap \bm{\{j\}} = t}} \hspace{-0.3cm} \big( \cW_{\bm{i},\bm{k},\bm{j},\bm{m}} - \cW_{\bm{i},\bm{k}} \cW_{\bm{j},\bm{m}} \big) + \mathcal{O}\Big( \frac{p^b}{n^{b+1}} \Big) \Bigg]\\
		& = \sum\limits_{b=1}^{(l_1+l_2) \land p} \Bigg[ \sum\limits_{t=1}^{b \land s} \frac{{s \choose t} {p-t \choose b-t} {n-b \choose l_1+l_2-b}}{{b \choose t} n^{l_1+l_2}} 2 b! (l_1+l_2-b)! \sum\limits_{k=0}^{b-t} t \, {l_1 \choose k} {l_1 \choose k+t} {l_2 \choose b-t-k} {l_2 \choose b-k}\\
		& \hspace{2cm} + \frac{s {p-1 \choose b-1} {n-b \choose l_1+l_2-b}}{b n^{l_1+l_2}} (\E[X_{1 \, 1}^4]-3) b! (l_1+l_2-b)! \sum\limits_{k=0}^{b-1} {l_1 \choose k} {l_1 \choose k+1} {l_2 \choose b-1-k} {l_2 \choose b-k}\\
		& \hspace{2cm} + \mathcal{O}\Big( \frac{p^b}{n^{b+1}} \Big) \Bigg] \ .
	\end{align*}
\end{proof}

\section{A technical lemma on bipartite trees}\label{Section_TechnicalLemma}
Let $K_{b+1,w+1}$ denote the fully connected bipartite graph with $b+1$ many vertices on the left hand side and $w+1$ many vertices on the right hand side. The spanning trees of $K_{b+1,w+1}$ are clearly all bipartite trees with connections to all the vertices in $K_{b+1,w+1}$.

\subsection{Lemma (Spanning trees of $K_{b+1,w+1}$ with the edge $(a_{1},c_{1})$)}\label{BipartiteSpanningTrees}
For any $b,w \geq 0$ and given $d_1,...,d_{b+1},e_1,...,e_{w+1} \in \N$ with the properties
\begin{align}\label{Eq_BipartiteSpanningTrees1}
	& d_1+...+d_{b+1} = b+w+1 = e_1+...+e_{w+1}
\end{align}
let $S$ denote the set of spanning trees $t$ of $K_{b+1,w+1}$ such that $t$ contains the edge $(a_{1},c_{1})$ and
\begin{align}\label{Eq_BipartiteSpanningTrees2}
	& \forall i \leq b+1 : \ \deg_t(a_i) = d_i\\
	& \forall j \leq w+1 : \ \deg_t(c_j) = e_j \ .
\end{align}
The cardinality of $S$ is
\begin{align*}
	& \#S = {b \choose e_1-1,...,e_{w+1}-1} \, {w \choose d_1-1,...,d_{b+1}-1}\\
	& \hspace{4cm} \times \begin{cases}
		1-\frac{(b-e_{1}+1)(w-d_{1}+1)}{bw} & \text{, if } w,b>0\\
		1 & \text{, if } w=0 \text{ or } b=0
	\end{cases} \ .
\end{align*}
By symmetry the same holds if we prescribe any edge $(a_{k},c_{l})$ instead of $(a_{1},c_{1})$. We only need to replace $d_{1}$ and $e_{1}$ in our formulas with $d_k$ and $e_{l}$.
\begin{proof}\
	\\
	We split this proof into three cases.
	\begin{itemize}
		\item[1)] Suppose $b=0$ or $w=0$:\\
		In this case one side of $K_{b+1,w+1}$ has only a single vertex, meaning there is only one possible spanning tree. This spanning tree contains the edge $(a_{1},c_{1})$ and thus the formula holds.
		
		\item[2)] Suppose either $d_{1}=1$ or $e_{1} = 1$, while $b,w>0$:\\
		Without loss of generality assume $d_{1}=1$, then $a_{1}$ must be a leaf. Let $c_j$ denote its neighbor. The tree is only counted, if $j=1$. Removal of the leaf $a_{1}$ (and relabeling the vertices $a_2,...,a_{b+1}$ into $a_1,...,a_b$) then defines a bijection between the set of all spanning trees $t$ of $K_{b+1,w+1}$ with (\ref{Eq_BipartiteSpanningTrees2}), where the leaf $a_{1}$ is connected to $c_{1}$, and the set of all spanning trees $\widehat{t}$ of $K_{b,w+1}$ with
		\begin{align*}
			& \forall i \leq b : \ \deg_{\widehat{t}}(a_i) = d_{i+1}\\
			& \forall 1 < j \leq w+1 : \ \deg_{\widehat{t}}(c_j) = e_j
		\end{align*}
		and $\deg_{\widehat{t}}(c_{1}) = e_{1}-1$. By (2.2) of \cite{MoonTrees} there are
		\begin{align*}
			& {b-1 \choose e_1-2,e_2-1,...,,e_{w+1}-1} {w \choose d_1-1,...,d_{b+1}-1}\\
			& = \frac{e_{1}-1}{b} {b \choose e_1-1,...,e_{w+1}-1} {w \choose \underbrace{d_{1}-1}_{=0},d_2-1,...,d_{b+1}-1}\\
			& = \bigg( 1 - \frac{b-e_{1}+1}{b} \underbrace{\frac{w-d_{1}+1}{w}}_{=1} \bigg) {b \choose e_1-1,...,e_{w+1}-1} {w \choose d_1-1,...,d_{b+1}-1}
		\end{align*}
		such trees.
		
		\item[3)] Suppose both $d_{1}$ and $e_{1}$ are larger than $1$ and $b,w>0$:\\
		Without loss of generality assume $b \geq w$, it is then easily seen that at least one element of $d_1,...,d_{b+1}$ must have value $1$ and by assumption this element is not $d_1$. Without loss of generality assume $d_{b+1}=1$ and let $\bm{S}_q$ denote the set of all spanning trees $t$ of $K_{b+1,w+1}$ with (\ref{Eq_BipartiteSpanningTrees2}), where the edge $(a_1,c_1)$ exists and $c_q$ is the only neighbor of $a_{b+1}$, which is a leaf since $d_{b+1}=1$.\\
		\\
		In order to use an inductive argument over the total number of vertices $N=w+b+2$, we first need to show that removing the vertex $a_{b+1}$ can only land us in the cases (2) or (3) and not in case (1). We could only land in case (1), if $b=1$ and thus also $w=1$ hold, as we had assumed both $b,w>0$ and $b \geq w$. However for $b=1=w$, we can not have $d_1 > 1$ and $w_1 > 1$ simultaneously, since then (\ref{Eq_BipartiteSpanningTrees1}) could not hold. It follows that we can not land in case (1) by removing the leaf $a_{b+1}$.\\
		\\
		We now inductively prove the formula
		\begin{align*}
			& \#S = {b \choose e_1-1,...,e_{w+1}-1} \, {w \choose d_1-1,...,d_{b+1}-1} \bigg(1-\frac{(b-e_{1}+1)(w-d_{1}+1)}{bw} \bigg) \ ,
		\end{align*}
		which is consistent with case (2), but not case (1). The above argument allows us to use case (2) as the start to the induction (effectively assume $b\geq 2$) and we only need to do the inductive step.\\
		\\
		For each $q \leq w+1$ we similarly as in (2) have a bijection between $\bm{S}_q$ and the set of spanning trees $\hat{t}$ of $K_{b,w+1}$ with connection $(a_1,c_1)$,
		\begin{align*}
			& \forall i \leq b: \, \deg_{\hat{t}}(a_i) = d_i\\
			& \forall j \leq w+1, \, j \neq q: \, \deg_{\hat{t}}(c_j) = e_j
		\end{align*}
		and $\deg_{\hat{t}}(c_q)=e_q-1$ by removing the leaf $a_{b+1}$. By inductive assumption we thus know the cardinality of $\bm{S}_q$ to be
		{\small
			\begin{align*}
				& \#\bm{S}_q = {b-1 \choose e_1-1,...,e_q-2,...,e_{w+1}-1} \, {w \choose d_1-1,...,d_{b}-1} \bigg(1-\frac{(b-e_{1})(w-d_{1}+1)}{(b-1)w} \bigg)
			\end{align*}
		}
		for $q \neq 1$ and
		{\small
			\begin{align*}
				& \#\bm{S}_1 = {b-1 \choose e_1-2,,...,e_{w+1}-1} \, {w \choose d_1-1,...,d_{b}-1} \bigg(1-\frac{(b-e_{1}+1)(w-d_{1}+1)}{(b-1)w} \bigg) \ .
			\end{align*}
		}
		It follows that
		{\small
			\begin{align*}
				& \#S = \sum\limits_{q=1}^{w+1} \#\bm{S}_q = \bm{S}_1 + \sum\limits_{q=2}^{w+1} \#\bm{S}_q\\
				& = {b-1 \choose e_1-2,,...,e_{w+1}-1} \, {w \choose d_1-1,...,d_{b}-1} \bigg(1-\frac{(b-e_{1}+1)(w-d_{1}+1)}{(b-1)w} \bigg)\\
				& \hspace{1cm} + \sum\limits_{q=2}^{w+1} {b-1 \choose e_1-1,...,e_q-2,...,e_{w+1}-1} \, {w \choose d_1-1,...,d_{b}-1}\\
				& \hspace{4cm} \times \bigg(1-\frac{(b-e_{1})(w-d_{1}+1)}{(b-1)w} \bigg)\\
				& = \frac{e_1-1}{b} {b \choose e_1-1,,...,e_{w+1}-1} \, {w \choose d_1-1,...,d_{b}-1,\underbrace{d_{b+1}-1}_{=0}} \bigg(1-\frac{(b-e_{1}+1)(w-d_{1}+1)}{(b-1)w} \bigg)\\
				& \hspace{1cm} + \sum\limits_{q=2}^{w+1}\frac{e_q-1}{b} {b \choose e_1-1,...,e_{w+1}-1} \, {w \choose d_1-1,...,d_{b}-1,\underbrace{d_{b+1}-1}_{=0}}\\
				& \hspace{4cm} \times \bigg(1-\frac{(b-e_{1})(w-d_{1}+1)}{(b-1)w} \bigg)\\
				& = {b \choose e_1-1,...,e_{w+1}-1} \, {w \choose d_1-1,...,d_{b+1}-1}\\
				& \hspace{1cm} \times \Bigg[ \frac{e_1-1}{b} \bigg(1-\frac{(b-e_{1}+1)(w-d_{1}+1)}{(b-1)w} \bigg) + \sum\limits_{q=2}^{w+1}\frac{e_q-1}{b} \bigg(1-\frac{(b-e_{1})(w-d_{1}+1)}{(b-1)w} \bigg) \Bigg] \ .
			\end{align*}
		}\
		\\
		By (\ref{Eq_BipartiteSpanningTrees1}) the sum $\sum\limits_{q=2}^{w+1} e_q$ must be equal to $b+w+1-e_1$ and the expression in the square bracket is calculated to be $\bigg(1-\frac{(b-e_{1}+1)(w-d_{1}+1)}{bw} \bigg)$. \qedhere
	\end{itemize}
\end{proof}

\newpage
\section*{List of symbols}\label{ListOfSymbols}
\addcontentsline{toc}{section}{List of symbols}
\begin{itemize}
	
	\item[] $A(G)$ \tabto{1.5cm} Adjacency matrix: a matrix with one row and column per vertex of $G$,\\
	\tabto{1.5cm} where $A(G)_{v_i,v_j}$ (slight abuse of notation) counts the number of edges\\
	\tabto{1.5cm} from $v_i$ to $v_j$ (If $G$ is undirected, edges are counted in both directions.)
	
	\item[] $\E[X_{1 \, 1}^4]$ \tabto{1.5cm} fourth moment of the entries $(X_{i\, j})_{i,j}$ of $\bm{X}$, i.e. $\E[X_{1 \, 1}^4] = \E[X_{1 \, 1}^4]$ (see \ref{MainResult1})
	
	\item[] $B(G)$ \tabto{1.5cm} set of black colored vertices in $G$ (see \ref{DefBlackWhiteCol} and \ref{DefDoubleCircuit})
	
	\item[] $\operatorname{Cov}$ \tabto{1.5cm} covariance of two random variables
	
	\item[] $\mathcal{C}_{V_r,N}$ \tabto{1.5cm} set of circuit graphs (see \ref{DefCircuitGraph})
	
	\item[] $\mathcal{C}^2_{V_r,N_1,N_2}$ \tabto{1.7cm} set of double-circuit graphs (see \ref{DefDoubleCircuit})
	
	\item[] $\deg$ \tabto{1.5cm} Degree: the number of edges of a vertex in an undirected (multi-)graph
	
	\item[] $\E$ \tabto{1.5cm} mean of a random variable
	
	\item[] $e_k$ \tabto{1.5cm} an edge of a directed multigraph (see \ref{DefDirectedMultigraph})
	
	\item[] $E_N$ \tabto{1.5cm} linearly oriented edge set of a directed multigraph (see \ref{DefDirectedMultigraph})
	
	\item[] $f_G$ \tabto{1.5cm} map from $E_N$ to $V_r \times V_r$, which defines the graph $G$ (see \ref{DefDirectedMultigraph})
	
	\item[] $G_{\bm{i}}$ \tabto{1.5cm} uniquely defined circuit graph to the route $\bm{i}$ (see \ref{DefCircuitGraph})
	
	\item[] $\mathcal{G}_{r,N}$ \tabto{1.5cm} set of directed multigraphs $G=(V_r,E_N,f)$ (see \ref{DefDirectedMultigraph})
	
	\item[] $\operatorname{head}$ \tabto{1.5cm} a graph-dependent map assigning the termination vertex of an edge\\
	\tabto{1.5cm}(see \ref{DefDirectedMultigraph})
	
	\item[] $\bm{i}$ \tabto{1.5cm} route of a circuit graph (see \ref{DefCircuitGraph})
	
	\item[] $\operatorname{indeg}$ \tabto{1.5cm} In-Degree: number of directed edges with a given vertex as their head
	
	\item[] $K_{b,w}$ \tabto{1.5cm} complete bipartite graph with $b$ left-hand-vertices and $w$ right-hand-vertices
	
	\item[] $\operatorname{outdeg}$ \tabto{1.5cm} Out-Degree: number of directed edges with a given vertex as their tail
	
	\item[] $\operatorname{R}(G)$ \tabto{1.5cm} reversed graph of $G$ (see \ref{DefReversed})
	
	\item[] $\bm{S}_{p,n}$ \tabto{1.5cm} sample covariance matrix (see page 1)
	
	\item[] $\operatorname{S}(G)$ \tabto{1.5cm} seed graph of $G$ (see \ref{DefSeed} and \ref{DefDoubleCircuit})
	
	\item[] $\operatorname{tail}$ \tabto{1.5cm} a graph-dependent map assigning the origin vertex of an edge (see \ref{DefDirectedMultigraph})
	
	\item[] $\mathcal{T}_{V_{l+1}}$ \tabto{1.5cm} set of balanced trees with vertex set $V_{l+1}$ (see \ref{DefBalancedTree})
	
	\item[] $\operatorname{tr}$ \tabto{1.5cm} trace of a matrix
	
	\item[] $\operatorname{U}(G)$ \tabto{1.5cm} undirected simple graph constructed by replacing undirected connections\\
	\tabto{1.5cm} of $G$ with edges (see \ref{DefUndirectedConnection})
	
	\item[] $v_i$ \tabto{1.5cm} a vertex of a directed multigraph (see \ref{DefDirectedMultigraph})
	
	\item[] $V_r$ \tabto{1.5cm} ordered vertex set of a directed multigraph with $r$ elements (see \ref{DefDirectedMultigraph})
	
	\item[] $V(G)$ \tabto{1.5cm} the set of visited vertices of $G$ (see \ref{DefVisitedVertex})
	
	\item[] $\bm{X}_{p,n}$ \tabto{1.5cm} random data matrix with $n$ data-points and $p$ features\\
	\tabto{1.5cm} (see the introduction)
	
	\item[] $X_{i,j}$ \tabto{1.5cm} entry of the random data matrix $X_{p,n}$ (see the introduction)
	
	\item[] $\land$ \tabto{1.5cm} gives the minimum of two numbers, i.e. $a \land b := \min(a,b)$
	
	\item[] $\lor$ \tabto{1.5cm} gives the maximum of two numbers, i.e. $a \lor b := \max(a,b)$
	
	\item[] $\bm{\{ \ \}}$ \tabto{1.5cm} set of vertices occurring in a sequence, i.e. $\bm{\{i\}} = \{i_1,...,i_N\}$ (see \ref{DefCircuitGraph})
	
	\item[] $[ \ ]$ \tabto{1.5cm} set of positive integers up to a given integer, i.e. $[N] = \{1,...,N\}$
	
	\item[] $\lfloor \ \rfloor$ \tabto{1.5cm} gives the nearest lower (or equal) whole number
	
	\item[] $\lceil \ \rceil$ \tabto{1.5cm} gives the nearest higher (or equal) whole number
	
	\item[] $\left< \ , \ \right>$ \tabto{1.5cm} zipped sequence of two sequences of equal length,\\
	\tabto{1.5cm} i.e. $\left< \bm{i}, \bm{k} \right> = (i_1,k_1,...,i_N,k_N)$ (see Section \ref{MethodOverview})
	
	\item[] $\operatorname{1-d-Ring}_{l_0}$ \tabto{3.5cm} set of one-dir. ring-type graphs of length $l_0$ (see \ref{DefRing})
	
	\item[] $\operatorname{2-d-Ring}_{l_0}$ \tabto{3.5cm} set of two-dir. ring-type graphs of length $l_0$ (see \ref{DefRing})
	
	\item[] $\operatorname{Double-1-d-Ring}_{l_0}$ \tabto{3.5cm} set of one-dir. double ring-type graphs of length $l_0$ (see \ref{DefDoubleRings})
	
	\item[] $\operatorname{Double-2-d-Ring}_{l_0}$ \tabto{3.5cm} set of two-dir. double ring-type graphs of length $l_0$ (see \ref{DefDoubleRings})
	
\end{itemize}

\newpage
\section*{Acknowledgments}\label{Declarations}
This work was supported by the DFG Research Unit 5381. It would not have been possible without the backing and enthusiasm of my supervisor Angelika Rohde.\\
\\
We thank Peter Taylor for proving a particular combinatorial identity for us on Mathoverflow (see \cite{PeterTaylor}). This identity allowed us to simplify the formula from Theorem \ref{Thm_MeanExpansion} by showing consistency of our results with an earlier, more specialized result of Bai and Silverstein.\\
\\
\\
A special thank you goes to Clemens Brüser for supplying the essential idea to the proof of Lemma \ref{BipartiteSpanningTrees} and for being a great colleague who always had an ear open for problems of graph theoretical nature. Best of luck in Dresden, you will be missed.\\
\\
\\

\section*{Declarations}
\addcontentsline{toc}{section}{Declarations}

The author is employed under the research unit 5381 of the DFG (Deutsche Forschungsgemeinschaft) at the University of Freiburg. Otherwise there was no funding relevant to this research.\\
\\
The author has no relevant financial or non-financial interests to disclose. The author has no competing interests to declare that are relevant to the content of this article. The author certifies that they have no affiliations with or involvement in any organization or entity with any financial interest or non-financial interest in the subject matter or materials discussed in this manuscript. The author has no financial or proprietary interests in any material discussed in this article.\\
\\

\subsection*{Data Availability}
Data sharing is not applicable to this article as no datasets are analysed.

\end{document}